\documentclass[11pt,a4paper,twoside]{article}
\usepackage{enumerate}

\usepackage{latexsym}
\usepackage{amsmath}
\usepackage{amssymb}
\usepackage{amsthm}
\usepackage{chicago}
\usepackage{amsbsy}
\usepackage{multicol}
\usepackage{color}
\usepackage{graphicx}

\usepackage[utf8]{inputenc}
\usepackage{amsfonts}
\usepackage{calligra}
\usepackage{calrsfs}
\usepackage[left=2cm,right=2cm,top=2cm,bottom=2cm]{geometry}
\usepackage[mathscr]{euscript}
\DeclareMathAlphabet{\matAthcalligra}{T1}{calligra}{m}{n}

\newcommand{\SPdp}{{\mathbb S}_+^{d-1}}
\newcommand{\pintegers}{{\mathbb N}_+}

\newcommand{\reals}{{\mathbb R}}
\newcommand{\realsd}{{\mathbb R}^d}

\newcommand{\SPd}{{\mathbb S}^{d-1}}
\newcommand{\halmos}{\vspace{3mm} \hfill \mbox{$\Box$}}
\newcommand{\DD}{\mathbb D}
\newcommand{\cG}{\mathfrak G}
\newcommand{\bfDV}{{\boldsymbol \delta}_v}

\newcommand{\D}{\mathcal{D}}
\newcommand{\Rdo}{\R^d_+\setminus\{0\}}
\newcommand{\Pstat}{\widehat{\Prob}^\alpha}
\newcommand{\Pstattheta}{\widehat{\Prob}^\theta}
\newcommand{\Estat}{\widehat{\E}^\alpha}

\allowdisplaybreaks

\newtheorem{thm}{Theorem}[section]
\newtheorem{cor}[thm]{Corollary}
\newtheorem{lemma}[thm]{Lemma}
\newtheorem{prop}[thm]{Proposition}
\newtheorem{sublemma}[thm]{Sublemma}

\theoremstyle{definition}
\newtheorem{defin}[thm]{Definition}
\newtheorem{rem}[thm]{Remark}

\numberwithin{equation}{section}

\usepackage{color}

\newcommand{\E}{\mathbb{E}}
\newcommand{\Prob}{\mathbb{P}}

\newcommand{\N}{\mathbb{N}}
\newcommand{\Z}{\mathbb{Z}}

\newcommand{\R}{\mathbb{R}}

\newcommand{\F}{\mathcal{F}}

\newcommand{\supp}{\mathrm{supp}\,  }

\renewcommand{\epsilon}{\varepsilon}
\renewcommand{\rho}{\varrho}

\newcommand{\1}[1][]{\mathbf{1}_{#1}}
\newcommand{\norm}[1]{\ensuremath{\left\| {#1} \right\|}}

\newcommand{\abs}[1]{\ensuremath{\left| {#1} \right|}}
\newcommand{\skalar}[1]{\langle #1 \rangle}

\newcommand{\mf}[1]{\mathfrak{#1}}
\newcommand{\M}{\mf{M}}

\newcommand{\Step}[1][]{\textsc{Step #1}}

\newcommand{\lm}{\lambda}

\newcommand{\Rd}{\R^d}

\newcommand{\Rdp}{{\R^d_>}}

\newcommand{\est}[1][\theta]{r^*_{#1}}

\newcommand{\es}[1][\theta]{r_{#1}}

\renewcommand{\P}[2][]{\ensuremath{\mathbb{P}_{#1} \left( {#2} \right)}}

\newcommand{\wt}{\widetilde}

\renewcommand{\Rdp}{\R^d_+}

\title{\bf  Large excursions and conditioned laws for 
recursive sequences generated by random matrices}

\author{\sc Jeffrey F.\ Collamore and Sebastian Mentemeier$^1$}

\date{\it University of Copenhagen and TU Dortmund}

\begin{document}







\maketitle

\footnotetext[1]{Corresponding author. 
\par \hspace*{.1cm} {\it AMS 2010 subject classifications.}  
Primary 60K15, 60F10. Secondary: 60J05, 60G70, 60G17.
\par  \hspace*{.1cm} {\it Keywords and phrases.}  Random recurrence equations, stochastic fixed point equations, products of random matrices,
Markov chain theory in general state space, nonlinear Markov renewal theory, large deviations, first passage times, conditional limit theorems,
extreme value theory.}

\begin{abstract}
We determine the large exceedance probabilities and large exceedance paths for the matrix recursive sequence
$V_n = M_n V_{n-1} + Q_n, \: n=1,2,\ldots,$
where $\{M_n\}$ is an i.i.d.\ sequence of $d \times d$ random matrices and $\{ Q_n\}$ is an i.i.d.\ sequence of random vectors, both with
nonnegative entries.
Early work on this problem dates to Kesten's (1973) seminal paper, motivated by an application to multi-type branching processes.  Other applications arise in financial time series modeling (connected to the study of the GARCH($p,q$) processes) and in physics, and this recursive sequence has also been the focus of extensive work in the recent probability literature. 
In this work, we characterize the distribution of the first passage time $T_u^A := \inf \{n:  V_n \in u A \}$, where $A$ is a subset
of the nonnegative quadrant in $\realsd$, showing that $T_u^A/u^\alpha$ converges to an exponential law.  In the process, we also revisit and refine
Kesten's classical estimate, showing that if $V$ has the stationary distribution of $\{ V_n \}$,
then  ${\mathbb P} \left( V \in uA \right) \sim C_A u^{-\alpha}$ as $u \to \infty$, 
providing, most importantly, a new characterization
of the constant $C_A$.  Finally, we describe the large exceedance paths via two conditioned limit laws.  In the first, we show that conditioned
on a large exceedance, the process $\{ V_n\}$ follows an exponentially-shifted Markov random walk, which we identify,
thereby generalizing results for classical random walk
to matrix recursive sequences.  In the second,
we characterize  the empirical distribution of $\{ \log |V_n| - \log |V_{n-1}| \}$  prior to a large exceedance, showing 
that this  distribution converges to the stationary law of the exponentially-shifted Markov random walk.
\end{abstract}

%
%
%
%

\section{Introduction}
The goal of this paper is to describe the extremal behavior, tail asymptotics, and conditioned path properties of the matrix  recursive sequence
\begin{equation} \label{intro1}
V_n = M_n V_{n-1} + Q_n, \quad n=1,2,\ldots, \qquad V_0 = v \in \realsd_+,
\end{equation}
where $\{ M_n \}$ is an i.i.d.\ sequence of $d \times d$ random matrices with nonnegative entries, $\{ Q_n \}$ is an i.i.d.\ sequence of nonnegative random vectors,
and $\realsd_+:=[0,\infty)^d$ denotes the nonnegative quadrant in $d$-dimensional Euclidean space.

Motivated by branching processes in random environments with immigration, as considered in Solomon \citeyear{FS72,FS75},
the matrix recursive sequence \eqref{intro1} was originally studied in the fundamental paper
of \shortciteN{HK73}. 
If $\E\big[ \log \norm{M_1} + \log|Q_1|\big] < \infty$ and the upper Lyapunov exponent is negative, i.e.
$$ \lim_{n \to \infty} \frac{1}{n} \log \norm{M_1 \cdots M_n} <0,$$ 
then it is readily verified that the law of
$$ V ~:=~ \sum_{k=1}^\infty M_1 \cdots M_{k-1} Q_{k}$$ is the unique stationary distribution for the Markov chain $\{V_n\}$ (see e.g. \shortciteN{Bougerol1992}).   Then under a Cram\'er-type condition stating that $\lambda(\alpha)=1$ for some $\alpha>0$, where $$\lambda(\theta) := \lim_{n \to \infty}  \left({\mathbb E} \left[ \| M_n \cdots M_1\|^\theta \right] \right)^{1/n}$$
and $\| \cdot \|$ denotes
operator norm, Kesten studied ${\mathbb P} \left( \big< w, V \big> > u \right)$ as $u \to \infty$ for $w \in \realsd_+$. It is shown in \shortciteN{HK73} that under appropriate moment and irreducibility conditions,
\begin{equation} \label{intro2}
{\mathbb P} \left( \big< w, V \big> > u \right) \sim {\cal C} u^{-\alpha} \quad \mbox{\rm as} \quad u \to \infty
\end{equation}
for a certain constant ${\cal C}$. 
Historically, this estimate resolved a conjecture by Spitzer,  verifying that $V$ lies in the domain of attraction of a stable law. 

Recently, there has been a renewed interest in Kesten's estimate.  The asymptotics in \eqref{intro2} have been shown to characterize
 the stationary tail decay in the 
GARCH($p,q$) financial time series models or, similarly, the ARMA($p,q$) processes with random coefficients; cf.\ \shortciteN{LHHRSRCV89}, \shortciteN{TM03}.
The process \eqref{intro1} is also relevant for the study of random walk in random environment (cf., e.g., \shortciteN{Kesten1975}, \shortciteN{Wang2013}),
and in a variety of other problems related to branching processes and Mandelbrot cascades; cf.\ \shortciteN{Gui1990}, \shortciteN{QL00},
\citeN{BDMbook}
and references therein.
Furthermore, in recent years, the scope of Kesten's method has broadened to include
more general fixed point equations in $\reals$; namely
equations of the form
\begin{equation} \label{introSFPE}
V \stackrel{d}{=} F(V),
\end{equation}
where $F:  \reals \to \reals$ is a random function independent of $V$, and $F(v) \approx Mv$ for large $v$, where $M$ is a random variable in
$\reals$; cf.\  \shortciteN{CG91},
\shortciteN{MM11}, Collamore and Vidyashankar (2013a,\,b), \shortciteN{GA16}. 
[Here, $\stackrel{d}{=}$ denotes equality in distribution.] Moreover, generalizations to Markov-dependent recursive sequences
(satisfying different assumptions from the processes we consider here)
have been obtained by \shortciteN{AR07} and \shortciteN{JC09}.

It is natural to ask whether this theory may be extended to reveal more refined path properties of the
process $\{ V_n \}$.  In particular, the behavior of $\{ V_n \}$ over large excursions
may essentially be inferred from that of the Markov random walk $\{ (X_n,S_n):  n=0,1,\ldots \}$, where 
\begin{equation} \label{intro2a}
X_n = \frac{M_n \cdots M_1 v}{|M_n \cdots M_1v |}, \qquad S_n = \log |M_n \cdots M_1 v|,
\end{equation}
and $|\cdot|$ denotes a norm in $\realsd$
(thus, $\{ X_n \}$ describes the directional component of the matrix product $M_n \cdots M_1 v$, and $\{ S_n \}$ describes its radial growth).
Note that $e^{S_n}X_n$ corresponds with $V_n$ when $Q=0$.
While the rough equivalence between $\{ V_n \}$ and $\{ e^{S_n}X_n \}$ has been utilized by numerous authors, including \shortciteN{HK73},
the correspondence between these processes has typically only been employed to obtain estimates
such as \eqref{intro2}, and not to characterize more detailed path properties.  In contrast, our approach will be to 
{\it quantify} this discrepancy using Markov nonlinear renewal theory,
as developed by Melfi (1992, 1994), yielding---after
accounting for the small-time behavior---that
the process $\{ V_n \}$ is closely approximated by $\{ e^{S_n} X_n \}$ in a manner which we characterize mathematically.
Consequently, it is natural to expect that, over a large excursion, the random walk structure inherent in
$\{ (X_n, S_n) \}$ may be exploited to yield deeper characteristics of the process $\{ V_n \}$ which mimic known properties 
of Markov random walk.
Following this approach, we shall reexamine Kesten's estimate, then extend the approach to obtain related asymptotic results relevant in extreme value theory,
and, ultimately,  derive path estimates conditioned on a large excursion, showing quantitatively that 
the path of $\{V_n\}$
 to a large exceedance roughly
follows that of $\{e^{S_n} X_n\}$ in an $\alpha$-shifted measure---also known as the exponentially tilted measure or Esscher transform---generalized to the setting of {\em Markov} random walk.

We start by revisiting \eqref{intro2}, establishing that, for an arbitrary set $A \subset \realsd_+$ satisfying certain regularity constraints,
\begin{equation} \label{intro3}
{\mathbb P} \left( V \in u A \right) \sim \frac{C}{\lambda^\prime(\alpha)} {\mathfrak L}_\alpha(A) \, u^{-\alpha}  \quad \mbox{\rm as} \quad u \to \infty.
\end{equation}
for a universal constant $C$ and a  measure ${\mathfrak L}_\alpha$.  
In particular, the constant $C$ is now {\it explicitly} identified as the $\alpha^{\footnotesize {\rm th}}$ moment of a certain power series 
derived from $\{ (M_n,Q_n) \}$
and the time-reversed products of $\{ M_n \}$; see \eqref{thm2A-eq1} and \eqref{thm2A-eq2}
below.  The formula we obtain can be viewed as a multidimensional extension of the main result
in  \shortciteN{JCAV13}.   (For a related one-dimensional estimate, see also \shortciteN{NECSOZ09}.)  Roughly speaking, 
$C$ is obtained
by studying the ratio $|V_n|/|M_n \cdots M_1 v|$ along a large excursion, and
intuitively, the constant arises by comparing the growth of $ |V_n|$ to the product $|M_n \cdots M_1 v|$ along paths
where both of these processes {\it diverge}; see the discussion in Remark \ref{remark:hitting times Vn} below. 

  Moreover, from \eqref{intro3}
we immediately conclude that $V$ is  {\em multivariate } regularly varying, as  could only be deduced  from \eqref{intro2}, based on the current literature, with the help of the Cram\'{e}r-Wold device; see \shortciteN{BBRDTM02} and \shortciteN{BL2009}. We 
emphasize that this additional step is not needed in our approach. 

Following a similar approach, we then examine
 the extremal behavior of $\{ V_n \}$.
Specifically, letting $A$ be a subset of the nonnegative quadrant and setting $T_u^A = \inf \{ n:  V_n \in uA \}$, we study the growth  rate
of $T_u^A$ as $u \to \infty$.   We show that
\begin{equation} \label{intro4}
\lim_{u \to \infty} {\mathbb P} \left( \frac{T_u^A}{u^{\alpha}} 
 \le z \, \bigg| \, V_0 = v \right) = e^{-K_A z}, \quad  z \ge 0,
\end{equation} 
where $\alpha$ is given as in \eqref{intro2} and $K_A$ is a constant which we also characterize, relating this constant explicitly to the
prefactor appearing in the asymptotic expression, as $u \to \infty$, for the
hitting probability of the set $uA$ of  $\{e^{S_n} X_n\}$ and to the constant $C$.  
As a special case,  setting $A=\{ x:  |x| > 1 \}$, we can then conclude that $\{ |V_n| \}$ belongs to the 
maximum domain of attraction of the Fr\'{e}chet distribution.
However, it should be emphasized that \eqref{intro4} is actually a stronger result, yielding the directional dependence of $\{ V_n \}$ and suggesting
a natural extension of classical extreme value theory to this multidimensional setting.
Note that in our estimate, $T_u^A < \infty$ a.s.; this contrasts from the asymptotics one obtains for perpetuity sequences
(i.e., the backward sequences corresponding to \eqref{intro1} in $\reals$).  In that setting, one 
obtains asymptotics which partly mimic those of random walk; cf.\  \citeN{DBJCEDJZ16}, which may be compared with
random walk estimates such as \shortciteN{SL84} or
\shortciteN{JC98}.  In contrast, \eqref{intro4} is qualitatively similar to {\it reflected} random walk,
and \eqref{intro4} can be viewed as an extension, to our setting, of a classical result due to \shortciteN{DI72}.
(For maximal segmental sums
of random walks, closely related estimates have also been provided by Dembo and Karlin (1991a,\:b), \shortciteN{SKAD92},
and  \shortciteN{ADSKOZ94}.)
Our result also sharpens earlier work, largely restricted to one-dimensional recursions, due to \shortciteN{LHHRSRCV89},
Perfekt (1994, 1997), and \citeN{BDMbook}; cf.\ Remark \ref{RkPerfekt} below.

The key to establishing  \eqref{intro3} and \eqref{intro4} is a proposition, where we study the behavior of $\{V_n\}$ over {\it cycles}
emanating from, and then returning to, a given set ${\mathbb D} \subset \realsd_+$.  Drawing an analogy with reflected random walk, 
these returns to ${\mathbb D}$ play the role of Iglehart's (1972) returns of a reflected random walk to the origin.
Letting $\tau$ denote the first return time to ${\mathbb D}$, then for any suitable function $g$ and any $m \in \{1,2,\ldots\}$,  we  consider 
in Proposition \ref{prop1} the limit behavior, as $u \to \infty$, of 
\[
u^\alpha {\mathbb E} \left[ g\left(\frac{V_{T_u^A}}{u}, \ldots, \frac{V_{T_u^A+m}}{u}\right) {\bf 1}_{\{T_u^A < \tau\}} \, \Big| \, V_0=v \right] .
\]
If $g =1$, then this represents the rescaled probability that $\{V_n\}$ enters the set $uA$ before returning to $\DD$.
Moreover, for general $g$, we show that the post-$T_u^A$-process behaves as  $\{e^{S_n} X_n\}$, but starting with the {\em stationary overjump distribution}.
This idea is then extended in the final section of the article
 to include the path behavior {\it prior} to time $T_u^A$ more explicitly,  drawing a close analogy to the behavior of the process $\{e^{S_n} X_n\}$ in the  $\alpha$-shifted
measure.

 The motivation of our concluding results is to establish an extension of a well-known estimate for random walk; namely, that a negative-drift random
 walk satisfying a Cram\'er-type condition and {\it conditioned to stay positive} behaves as its associate  (i.e., the random walk in the $\alpha$-shifted measure); cf.\ \shortciteN{WF71}, Section XII.6.(d); \shortciteN{Bertoin.Doney:1994}.
 Similarly, a negative-drift random walk {\it conditioned to achieve a high barrier at level} $u$ will also 
 converge to its associate as the level $u \to \infty$; cf.\ \shortciteN{SA82}.  Thus, it is natural to expect that, as $u \to \infty$,
 the process $\{ V_n \}$ will behave analogous to   $\{e^{S_n}X_n\}$ in the $\alpha$-shifted measure,
  and in Section \ref{Sec2.4}, we 
make this idea precise.  
 As a special case, we then consider the empirical law of $\{ \log |V_n| - \log |V_{n-1}| \}$ conditioned on $\{T_u^A < \tau \}$.
 We show that for any suitable continuous function $g$,
\begin{equation}\label{intro7}
\E \left[ \left. \bigg| \frac{1}{T_u^A} \sum_{k=1}^{T_u^A} g\left( \log \bigg(  \frac{|V_k|}{|V_{k-1}|} \bigg) \right) 
  ~-~ \widehat{\E}^\alpha \left[g( S_1) \right]  \bigg|  \, \right| \, V_0=v, \: T_u^A < \tau \right] ~\to ~ 0
 \end{equation}
as $u \to \infty$, where $\widehat{{\mathbb E}}^\alpha [\cdot]$ denotes expectation, under stationarity, in the $\alpha$-shifted
measure, and $S_1$ is given as in \eqref{intro2a}.  Thus,  the empirical law of $\left\{ \left(\log|V_n| - \log|V_{n-1}| \right) \right\}$ converges weakly
in ${\mathbb P} \left( \cdot \, | T_u^A < \tau \right)$-probability to the distribution,
under stationarity, of $S_1$.

We emphasize that we shall develop our limit theorems {\it without} the assumption that the Markov chain $\{ V_n \}$
 is Harris recurrent, and thus---while we shall occasionally draw upon the theory of Harris recurrent chains---our appoach
will differ markedly from the more classical approach outlined, for example, in \shortciteN{PNEN87a}.
Indeed, the assumption of Harris recurrence is rather unnatural in our setting.  Instead, we circumvent
this requirement by introducing a smoothing technique,
where the sequence
$\{ Q_{kn} \}$ is ``smoothed" for some $k \in \{1,2,\ldots\}$, thereby ensuring that the resulting
process is Harris recurrent, yet the effect of this smoothing is negligible in an asymptotic limit.
As this technique  could be adapted to 
 other recursive sequences satisfying a stochastic fixed point equation of the form \eqref{introSFPE}---where the assumption of Harris recurrence is also restrictive---this method could potentially be 
 of some general interest.   Indeed, rather than assuming Harris recurrence, 
 we shall rely throughout the article on the recently-developed theory of \shortciteN{Guivarch2012}, which exploits spectral gap properties on special function spaces for matrix products under weak regularity conditions. While the theory of Guivarc'h and Le Page is developed for invertible matrices, a formulation for matrices with nonnegative entries, as we consider here, is given in \shortciteN{Buraczewski.etal:2014}.
 We now turn to a precise statement of our main results.

%
%
%
%

\setcounter{equation}{0}
\section{Statement of results}
\label{sect:2}

\subsection{Preliminaries:  notation, assumptions, and background}
Let $d \ge 1$ and let $M$ be a  $d \times d$ random matrix whose entries are nonnegative a.s., and
let $Q$ be a random vector in $\realsd$ with nonnegative entries a.s.
Let $\mu$ denote the probability law of $(M,Q)$, and let $\mu_M, \mu_Q$ denote the marginal laws of $M$ and $Q$, respectively. Now let $\{ ( M_n, Q_n):  n \in \pintegers \}$ be a sequence of i.i.d.\ copies of $(M, Q)$ where, here and in the following, 
$\pintegers := \{1,2,\ldots \}$ denotes the positive integers.
We assume throughout the article that $\{ (M_i,Q_i):  i=1,\ldots, n \}$ are adapted to a given filtration $\{{\cal F}_i:  i = 1,2,\ldots\}$. 

Then the aim of this paper is to study the extremal properties of the stochastic recursive sequence defined by
\begin{equation}\label{RDE} 
V_n = M_{n} V_{n-1} + Q_n, \quad n=1,2,\ldots, \quad V_0\sim \nu,
 \end{equation}
for a given initial distribution $\nu$, where, unless specifically noted, $\nu$ is concentrated at a deterministic point $v \in \realsd$.

Next we introduce some additional notation, as follows.
Let $\realsd$ be endowed with the
scalar product $\big< \cdot, \cdot \big>$ and canonical orthonormal basis $\{ e_i:  i=1, \ldots, d \}$.  
The nonnegative cone in $\realsd$ is defined by
\[
\realsd_+ = \left\{ x \in \realsd:  \big< x, e_i \big>  \ge 0 \right\}.
\]

Let $| \cdot |$ denote a norm in $\realsd$, and assume throughout the article that $| \cdot |$ is monotone, i.e., if $x, \, y \in {\mathbb R}^d_+$ satisfy $y - x \in {\mathbb R}^d_+$,
then $|x| \le |y|$. 
Let $\SPd := \{ x \in \Rd \, : \, \abs{x}=1\}$ denote the unit sphere and
${\mathbb S}_+^{d-1} := \realsd_+ \cap \SPd$.   For any $x \in \realsd$, let $\widetilde{x}$ denote its projection onto the unit sphere, that is,
\[
\widetilde{x} \equiv \left( x \right)^\sim := |x|^{-1} x, \quad x \in \realsd.
\]
For any subspace of $\realsd$, let ${\cal B}({\cal S})$ denote the collection of Borel sets on ${\cal S}$; and
for any $E \in {\cal B}({\cal S})$, let $E^\circ$, $\bar{E}$, $E^c$, and $\partial E$ denote the interior, closure, complement, 
and boundary of $E$, respectively.  For any $r>0$ and $y \in \Rd$, let $B_r(y)=\{x \in \Rd \, : \, |x-y| < r\}$.  For any measure $\nu$ on ${\cal S} \subset \realsd$, let $\supp \nu$ denote the support of $\nu$.
Also, denote the set of bounded continuous real-valued functions on a space $E$  by ${\cal C}_b(E)$, equipped with the supremum norm, 
namely $|f|_\infty:=\sup \{|f(x)| \, : \, x \in E\}$.

Let ${\mathfrak M}$ denote the collection of $d\times d$ matrices having nonnegative 
coefficients, and let $\| {\mathfrak m} \|$ denote the operator norm, i.e.,
\[
\| {\mathfrak m} \| := \sup_{x \in {\mathbb S}^{d-1}} | {\mathfrak m} x |, \quad {\mathfrak m} \in {\mathfrak M}.
\]

\textbf{\textit{Allowable and positively regular matrices.}}
We say that a matrix ${\mathfrak m} \in {\mathfrak M}$ is {\it allowable} if it has no zero row or column.
Moreover, if the coefficients of a given matrix ${\mathfrak m} \in {\mathfrak M}$ are {\it strictly} positive, 
then we write ${\mathfrak m} \succ 0$ and say that ${\mathfrak m}$ is {\it positively regular}.  With a slight abuse of notation, we write
\[
{\mathfrak M}^\circ = \left\{ {\mathfrak m} \in {\mathfrak M}:  {\mathfrak m} \succ 0 \right\}.
\]
{\it As a standing assumption, we shall always assume that there exists an $n \in \pintegers$ such that 
${\mathfrak N} := \inf \left\{ n \in \pintegers:  M_n \cdots M_1 \succ 0 \right\}$ is finite a.s.; that is, ultimately, the product $M_n \cdots M_1$ is positively regular
with probability one.}
This assumption will be subsumed in the stronger Hypothesis $(H_1)$, given below (as can be inferred from \shortciteN{Hennion1997}, Lemma 3.1 or \shortciteN{Buraczewski.etal:2014}, Lemma 6.3).

Under this standing assumption, the elements of the vectors $\{ V_n \}$ communicate, leading to a common polynomial decay rate for the exceedance probabilities,
regardless of direction, while the pre-factor $C_A := (C/\lambda^\prime(\alpha)) {\mathfrak L}_\alpha(A)$ in \eqref{intro3}, or the constant $K_A$ in \eqref{intro4}, will be directionally-dependent.
\\[-.2cm]

\textbf{\textit{Non-arithmetic distributions for random matrices.}}   We will need a generalization of the notion of a non-arithmetic distribution to the setting of random
matrices.  First consider a more general framework, where $\{ (X_n,S_n):  n=0,1,\ldots \}$ is a Markov random walk, i.e.\ $\{(X_n, S_n-S_{n-1}) \}$ is a Markov chain
with a transition kernel which only depends on the state of the driving chain $\{X_n\}$. 
The most satisfactory generalization of an arithmetic distribution in this setting
is due to \shortciteN{VS84}.  In his formulation, the Markov random walk $\{ (X_n,S_n):  n=0,1,\ldots \}$ is {\it arithmetic}
if there exists a $t >0$, $\theta \in [0,2\pi)$, and a ``shift-function" $\vartheta$ on $\reals$ such that
\[
{\mathbb E} \left[ \exp \left\{ i t S_1 - i \theta + i\left( \vartheta(X_1) - \vartheta(X_0) \right) \right\} \right] = 1,
\]
and {\it non-arithmetic} otherwise.  Clearly, if $\{ S_n \}$ denotes the sums of an i.i.d.\ sequence of random variables (rather than a Markov-dependent sequence), then we may take $\vartheta = 0$,
and the above condition is equivalent to the requirement that the support of $S_1$ lies in an arithmetic progression; that is, the distribution of $S_1$ is arithmetic in the classical sense.

This condition is not easily verified in the matrix setting, so it is more natural to adopt a requirement akin to that of  \shortciteN{HK73}.  Namely,
let $\Gamma_M$ denote the smallest closed subsemigroup of $\M$ which contains $\supp  \mu_M$. 

\begin{defin} \label{nonarith}
We say that $\mu_M$ is {\it non-arithmetic} if the additive group generated by
\[
\left\{ \log \norm{\mathfrak m}:  {\mathfrak m} \in \Gamma_M \cap {\mathfrak M}^\circ \right\}
\]
is dense in $\reals$.
\end{defin}

It is shown in \shortciteN{BM2013}, Lemma 2.7, that this condition implies that of \shortciteN{VS84}.
It is also worth observing that, alternatively, we could replace $\log \| {\mathfrak m} \|$ with the Frobenius eigenvalue of ${\mathfrak m}$ in Definition \ref{nonarith}; thus, we see that our definition is, indeed, in
agreement with the one given in \shortciteN{HK73}. 

We are now prepared to introduce our basic assumptions on the distribution function $\mu_M$ of $M$.\\[-.1cm]

\noindent
{\bf Hypothesis \boldmath{($H_1$)}.}   $\mu_M$ is non-arithmetic, and
$\mu_M \{ {\mathfrak m}:   {\mathfrak m}$ is allowable$\} = 1$.\\[-.2cm]

Since we will employ Markov renewal theory, the appearance of a non-arithmetic assumption is natural. 
The further requirement that the matrix $M$ is allowable $\mu_M$-a.s. is also standard and appears in numerous related works in the literature (cf., e.g., 
\shortciteN{Kesten1984}, \shortciteN{Hennion1997}, \shortciteN{Buraczewski.etal:2014}).
In comparison, \shortciteN{HK73} assumes that the rows of $M$ are nonzero a.s., but does not assume that the columns of $M$ are also nonzero a.s. (as we assume by requiring
that the matrices are ``allowable").
This further requirement will guarantee the uniqueness of the invariant measures and functions in Lemma \ref{prop:ip} below.\\[-.2cm]

We now turn to certain moment conditions which will be imposed on the pair $(M,Q)$.  Set
\[
{\mathfrak D} = \left\{ \theta \ge 0 \: : \: \int_{\mathfrak M} \norm{\mathfrak m}^\theta \mu_M (d{\mathfrak m}) < \infty \right\} = \left\{ \theta \ge 0 \: : \: {\mathbb E} \left[ \norm{M}^\theta \right]
  < \infty \right\}.
\]
Let ${\mathfrak m}^T$ denote the transpose of ${\mathfrak m}$.  Then
for any $\theta \in {\mathfrak D}$, define:
\begin{align*}
P_\theta f(x) &= {\mathbb E} \left[ \abs{M x}^\theta f\big( \widetilde{M x} \big) \right],  \quad f \in {\cal C}_b(\SPdp);\\[.2cm]
P_\theta^* f(x) &= {\mathbb E} \left[ \abs{M^T x}^\theta f\big( \widetilde{M^T x} \big) \right],  \quad f \in {\cal C}_b(\SPdp);
\end{align*}
and
\[
\lm(\theta) = \lim_{n \to \infty} \left( {\mathbb E} \left[ \norm{M_n \cdots M_1}^\theta  \right] \right)^{1/n}; \qquad
\Lambda(\theta) = \log \lambda(\theta).
\]
\begin{lemma}\label{prop:ip}
Assume $\theta \in {\mathfrak D}$ and $\mu_M \{ {\mathfrak m}:   {\mathfrak m}$ is allowable$\} = 1$.
Then $\lambda(\theta)$ 
is the spectral radius of $P_\theta$, and there is a unique probability measure $l_\theta$ on $\SPdp$ and a unique, strictly positive function $\es \in {\cal C}_b \big( \SPdp \big)$ with  $\int_{{\mathbb S}^{d-1}_+} r_\theta(x) d l_\theta(x) =1$ such that 
\begin{equation}
l_\theta P_\theta = \lm(\theta) l_\theta, \qquad P_\theta \es = \lambda(\theta) \es.
 \end{equation}
Moreover, the function $\es$ is $\max\{\theta,1\}$-H\"older continuous; thus, in particular, $\es$ is bounded from above and below by finite positive constants.

Similarly, the spectral radius of $P^\ast_\theta$ equals $\lambda(\theta)$, and there is a unique probaility measure $l_\theta^*$ on $\SPdp$ and a unique, strictly positive function $\est[\theta]$ such that
$$
l_\theta^* P_\theta^* = \lm(\theta) l_\theta^*,\quad P_\theta^* \est = \lambda(\theta) \est,
\quad \mbox{\rm and} \quad \int_{{\mathbb S}^{d-1}_+} r^*_\theta(x) d l_\theta^*(x) =1.$$
Moreover,
\begin{equation}\label{eq:relation rs and lstar}
\est(x) ~=~ c \int_{\SPdp} \skalar{x,y} ^\theta l_\theta(dy) \quad \text{\rm for all } x \in \SPdp,
\end{equation}
where $c =  \left( \int_{{\mathbb S}^{d-1}_+ \times {\mathbb S}^{d-1}_+} \big< x, y \big>^\alpha l_\theta^\ast(dx) l_\theta(dy) \right)^{-1}$.
Furthermore, \eqref{eq:relation rs and lstar} also holds if  $r_\theta^\ast$ and $l_\theta$ are replaced with
$\es$ and $l_\theta^*$, respectively.
\end{lemma}

In the above lemma, we have written $l_\theta P_\theta$ for the application of the adjoint operator $P_\theta^\prime$ to the measure $l_\theta$, i.e. $l_\theta P_\theta$ is the unique measure satisfying
$$ \int_{\SPdp} f(x) \big(l_\theta P_\theta)(dx) ~=~ \int_{\SPdp} \big(P_\theta f(x)) \, l_\theta(dx)$$
for all $f \in {\cal C}_b(\SPdp)$.
The proof of Lemma 2.1 can be found in \shortciteN{Buraczewski.etal:2014}, Proposition 3.1; see also \shortciteN{Guivarch2012}, Theorem 2.16 for an
analogous result in the setting of invertible matrices.
For some related results for  Harris recurrent Markov chains;
see, for example, \shortciteN{EN84}, \shortciteN{PNEN87a}, or \shortciteN{AM2012}.

For any allowable matrix $\mathfrak{m}$, now define
$$ \textbf{\textit{i}}(\mathfrak{m}) ~:=~ \inf_{x \in \SPdp} \abs{\mathfrak{m} x}.$$

\noindent
{\bf Hypothesis \boldmath{($H_2$)}.}  There exists an $\alpha >0$ such that $\lambda(\alpha) =1$, and the following moment conditions hold:
\[
{\mathbb E} \Big[ \| M \|^\alpha \max \{\abs{\log \norm{\mathfrak{m}}}, \abs{\log \textbf{\textit{i}}(\mathfrak{m})} \} \Big] < \infty; \quad \mbox{and} \quad
{\mathbb E} \big[ | Q |^\alpha \big] < \infty.
\]\\[-.6cm]

Once again, $(H_2)$ is quite standard, also in the one-dimensional setting; cf.\ \shortciteN{CG91}.  In comparison with our assumptions,
\shortciteN{HK73} requires the slightly weaker condition $\E \big[ \norm{M}^\alpha \, \log^+ \norm{M} \big] < \infty$, rather than
\[
{\mathbb E} \Big[ \| M \|^\alpha \max \{\abs{\log \norm{\mathfrak{m}}}, \abs{\log \textbf{\textit{i}}(\mathfrak{m})} \} \Big] < \infty.
\]
However, from our modest strengthening of his condition, we will be able to identify the limit in the  Furstenberg-Kesten theorem, as given below (in
an extended form) in Lemma \ref{lem:slln}.\\[-.2cm]

\textbf{\textit{The shifted distribution. }}  We shall utilize the constant $\alpha$ in $(H_2)$ to employ a change of measure.
Now in the one-dimensional setting, it is natural to apply this change of measure to the first member of the pair $(\log M,Q) \subset \reals^2$ (cf.\ \shortciteN{JCAV13}, \shortciteN{JCGDAV14}),
and our main objective here is to extend this idea to the multidimensional framework.

Adopting the setting of random matrices and proceeding more formally, let ${\mathfrak m}$ be any allowable nonnegative matrix; and for any $\theta \in {\mathfrak D}$ and any 
$n \in \pintegers$,
introduce the density function
\[
p_n^{\theta} (x,{\mathfrak m}) = \frac{\abs{{\mathfrak m} x}^\theta}{\left(\lm(\theta)\right)^n} \frac{r_\theta \big( \widetilde{{\mathfrak m} x} \big)}{\es(x)}, \quad\quad x \in \SPdp.
\]
Note by an application of Lemma \ref{prop:ip} that
\[
 \int p_n^{\theta}(x, {\mathfrak m}_n \cdots {\mathfrak m}_1)\ \mu^{\otimes
n} \left( \{d{\mathfrak m}_i,dq_i\}_{i=1}^n \right) =1, \quad x \in \SPdp, \quad \theta \in {\mathfrak D}.
\]
Moreover, the system of probability measures $\mu_{n,x}^{\theta} = p_n^{\theta}(x, \cdot)
\mu^{\otimes n}$ is a projective system; hence by the Kolmogorov extension
theorem, there exists a unique probability measure ${\mathbb P}_x^{\theta}$ on $({\cal M} \times
\realsd_+)^{\pintegers}$ having marginals $\mu_{n,x}^{\theta}$.  When the random variables $\left\{ (M_n, Q_n):  n=1,2,\ldots \right\}$
are generated by the measure ${\mathbb P}_x^{\theta}$ rather than the true underlying probability measure, we shall write ${\mathbb E}_x^{\theta} [\cdot]$.  
{\it We shall refer to this distribution as the ``$\theta$-shifted distribution."}

 It is worth observing that, although $\{ (M_n,Q_n): n=1,2,\ldots \}$ is assumed to be i.i.d.\ in the unshifted measure,  
this sequence will be Markov-dependent in the $\theta$-shifted measure for any $\theta >  0$.
However, defining \begin{equation}\label{def:eta}
\eta_\theta (E) = \int_E \es[\theta](x) l_\theta(dx),
\end{equation}
we see that the measure
\begin{equation} \label{stationary-theta}
 \Pstattheta ~:=~\int_{\SPdp} {\mathbb P}_x^\theta \: \eta_\theta(dx)
 \end{equation}
 is shift-invariant; i.e., the sequence $\{(M_n, Q_n)\}$ is stationary under $\Pstattheta$; cf.\ Section 3.1 of \shortciteN{Buraczewski.etal:2014}.
 This is an important observation, as it will allow us to apply the results of \shortciteN{Hennion1997} on products of random matrices; cf.\ Section 4 below.
In addition, by Lemma 6.2 of \shortciteN{Buraczewski.etal:2014}, we have that
 ${\mathbb P}_x^\theta \ll \Pstattheta$, for all $x \in \SPdp$. We will use this result 
 frequently to infer convergence $\Prob^\alpha_x$-a.s., for arbitrary $x \in \SPdp$, by proving $\Pstattheta$-a.s.\ convergence.
Furthermore, the finite-dimensional distributions of $\{ (M_i, Q_i):  i=0,\ldots,n \}$ under each ${\mathbb P}_x^\alpha$ are equivalent to their unshifted distributions, since the density function $p_n^\theta$ is strictly positive.  

\medskip

\textbf{\textit{The Markov random walk.}}  The process $\{ M_n \}$ induces a Markov chain on ${\mathbb S}_+^{d-1}$,
defined by setting
\begin{equation} \label{def-X}
X_n = \left( M_n \cdots M_1 X_0 \right)^\sim, 
\end{equation}
for some initial state $X_0 \in  {\mathbb S}^{d-1}_+$.
This process will play an important role in the sequel. In the $\theta$-shifted measure,
$\{ X_n \}$ has a unique stationary distribution given by $\eta_\theta$.
If we further define 
\begin{equation}\label{def-S}
S_n = \log | M_n \cdots M_1 X_0 |, \quad n=1,2,\ldots; \quad S_0 = 0;
\end{equation}
then $\{ (X_n,S_n):  n=0,1,\ldots \}$ forms a Markov random walk, which we will utilize frequently below due to the fact that,
in the $\alpha$-shifted measure,
$\{ (\tilde{V}_n,\log|V_n|-\log|V_{n-1}|)  \}$
closely resembles  $\{ (X_n,S_n-S_{n-1}) \}$ for large $n$.

\medskip

\textbf{\textit{Probability measures.}}  We introduce the following conventions to describe conditional probabilities which
depend specifically on the
 initial values of $X_0$ and $V_0$.  Write:
\begin{align*}
\Prob_v(\cdot)=\Prob(\cdot | V_0=v), \quad 
\Prob^\theta_x(X_0=x)=1, \quad 
\Prob^\theta_{x,v}(\cdot)=\Prob^\theta_x( \cdot | V_0=v),
\end{align*}
and use the same notation for expectations.  When conditioning on an initial distribution $V_0 \sim \gamma$, write:
\begin{align*}
{\mathbb P}_\gamma (\cdot) =\int_{\Rdo} {\mathbb P}_v(\cdot) \, \gamma(dv) = {\mathbb P} (\cdot | V_0 \sim \gamma), \quad
{\mathbb P}^\theta_\gamma(\cdot) = \int_{\Rdo} {\mathbb P}^\theta_{\tilde{v},v}(\cdot) \gamma(dv), \quad {\mathbb P}^\theta_{\bfDV}(\cdot) =  {\mathbb P}^\theta_{\tilde{v},v}(\cdot),
\end{align*}
and analogously for expectations.   Note that while working in the $\theta$-shifted measure, we must specify {\it both} $X_0$ and $V_0$
in these last equations, and we specifically take $X_0 = \widetilde{V}_0.$  The reasoning for this asymmetry comes from the observation that,
while the initial state $X_0$ does not influence the distribution of $\{ V_n \}$ in the original measure, this initial state
does affect the law of $\{M_n\}$ and hence that of $\{ V_n \}$ in the $\theta$-shifted measure and, thus, both of the initial states, $X_0$ and $V_0$, must be specified in the latter probabilities or expectations. 
Finally, we will often suppress the dependence on $(x,v)$ and write $\Prob^\alpha$-a.s.\ in place of $\Prob^\alpha_{x,v}$-a.s. for all $x \in \SPdp$, $v \in \Rdo$.
With $\{(X_n, S_n)\}$ defined as above, using that $\lambda(\alpha)=1$, we can rewrite the change of measure as follows:
for all $n \in \pintegers$, $x \in \SPdp$, and any bounded measurable function $f : \SPdp \times (\mathfrak{M} \times \Rdp)^n$, 
\begin{equation}
r_\alpha(x) \E_{x,v}^\alpha \Big[ \frac{e^{-\alpha S_n}}{r_\alpha(X_n)} f(X_0,V_0, M_1, Q_1, \ldots, M_n, Q_n) \Big] ~=~ \E \Big[ f(x,v, M_1, Q_1, \cdots, M_n, Q_n) \Big].
\end{equation}

\medskip

%

\textbf{\textit{Limit theory for the Markov random walk.}}
Now recall the theorem of \shortciteN{Furstenberg1960}, which states that, if $\E\big[ \log \norm{M_1} \big] < \infty$, then 
\begin{equation} \label{FurstKestThm}
 \lim_{n \to \infty} \frac{1}{n} \log \norm{M_n \cdots M_1} 
\end{equation}
converges to its ergodic average.  We will need a refinement of this result, developed in the context of the $\theta$-shifted measure, where $\theta \in {\mathfrak D}^\circ$.
Before stating this result, first recall that under $(H_1)$,  ${\mathfrak N} :=  \inf\big\{ n \in {\mathbb Z}_+ : \, M_n \cdots M_1 \succ 0 \big\} < \infty$; cf.\ \shortciteN{Hennion1997},
Lemma 3.1.

\begin{lemma}\label{lem:slln}
Assume that $(H_1)$ is satisfied and let
$\theta \in \mathfrak{D}^\circ$, and suppose that
\[
{\mathbb E} \left[ \| M \|^\theta \max \{\abs{\log \norm{\mathfrak{m}}}, \abs{\log \textbf{\textit{i}}(\mathfrak{m})} \} \right] < \infty.
\]
Then in the $\theta$-shifted measure, we have for all $x \in \SPdp$ that
\begin{equation} 
\lim_{n \to \infty} \frac{1}{n} \log \abs{M_n \cdots M_1 x} 
= \lim_{n \to \infty} \frac{1}{n} \log \norm{M_n \cdots M_1} ~
= \Lambda^\prime(\theta) ~=~ \hat{\E}^\theta [S_1] \qquad \Prob^\theta\text{\rm -a.s};
\end{equation}
and for all $x, \, y \in \SPdp$,
\begin{equation} \label{FurstKestThm-1}
\lim_{n \to \infty} \, \sup \bigg\{ \abs{\frac1n \1[{\{ {\mathfrak N} \le n\}}] \log \skalar{y, M_{n} \dots M_{1} x} - \Lambda^\prime(\theta)} \, : \, x,y \in \SPdp   \bigg\} ~=~0 \qquad \Prob^\theta\text{\rm -a.s.}
\end{equation}
\end{lemma}

\noindent
{\bf Proof.} See \shortciteN{Hennion1997}, Theorem 2 (where the uniformity is proved) and \shortciteN{Buraczewski.etal:2014}, Theorem 6.1 (where the limit is identified in the $\theta$-shifted measure). \halmos

Finally, when studying the empirical measure conditioned on a large exceedance, it will be helpful to compare with
the {\it unconditional} behavior of the corresponding $\alpha$-shifted Markov random walk. For this purpose, the following lemma is useful.

\begin{lemma}\label{lem:ergodic theorem}
Assume Hypotheses $(H_1)$ and $(H_2)$ are satisfied, and suppose that the transition kernel of $\{ (X_n,S_n):  n=0,1,\ldots \}$
follows the $\alpha$-shifted measure, ${\mathbb P}^\alpha$.
Then for any measurable function $g : \SPdp \times \SPdp \times \R \to \R$ and any initial state $X_0 = x \in \SPdp$,
$$ \lim_{n \to \infty} \frac{1}{n} \sum_{k=1}^n g\left(X_k, X_{k-1}, S_k - S_{k-1} \right) ~=~ \Estat \big[g(X_1, X_0, S_1) \big] \quad \Prob^\alpha_x\text{\rm -a.s.},$$
provided that the expectation on the right-hand side exists.
\end{lemma}

\noindent
{\bf Proof.}
This is proved as in the first part of \shortciteN{Buraczewski.etal:2014}, Theorem 6.1 [employing, in the notation of that article, the function $f(x,\omega):=g(a_1(\omega).x,x, \log|a_1(\omega)x|)$].
\halmos

%
%
%
%

\subsection{Tail estimates for \boldmath$\{V_n\}$}
We now turn to our first main result, where we revisit and extend Kesten's (1973) theorem, establishing an {\it explicit} expression
for the constant $C_A$ in \eqref{intro3}. 

 Let $\pi$ denote the stationary distribution of $\{ V_n \}$, which is given by the law of the random variable 
$$V ~:=~ \sum_{k=1}^\infty M_1 \cdots M_{k-1} Q_{k-1}.$$
Then it is well known that $V$ is finite a.s.\ under the hypotheses of this paper, and hence $\pi$ exists; cf.\ \shortciteN{HK73}.
[The necessary moment hypotheses follow from $(H_2)$, while the negativity of the upper Lyapunov exponent follows from Lemma \ref{lem:slln} and the convexity of $\Lambda$, which implies that $\Lambda'(0)<0$.]

Now fix a set ${\mathbb D} \subset \realsd_+,$ where $\pi(\DD)>0$, and
let $\pi_{\mathbb D}$ denote the stationary distribution of $\{V_n\}$ restricted to ${\mathbb D}$; that is, 
\begin{equation} \label{def-PI-D}
\pi_{\mathbb D}(E) = \frac{ \pi( E \cap {\mathbb D})}{ \pi ({\mathbb D})}, \quad E \in {\cal B}(\reals^d_+);
\end{equation}
and let $\tau$ denote the first  return time   of $\{V_n\}$ into ${\mathbb D}$; namely,
\[
\tau = \inf\{n  \in \pintegers:  V_n \in {\mathbb D} \}.
\]

Next, let 
 ${\vec 1} = (1, \ldots, 1)^T$, and define
\begin{equation} \label{def-upsilon}
Y_i   = \lim_{n \to \infty} \left( M_i^\top \cdots M_n^\top {\vec 1} \right)^\sim, \quad n = 1,2,\ldots.
\end{equation}
Note that if  $\theta \in {\mathfrak D}$, then the limit on the right-hand side exists ${\mathbb P}^\theta$-a.s., since this product constitutes a backward sequence
of an iterated function system and the maps $\{M_k\}$  act as contractions on $\SPdp$; cf.\ \shortciteN{Hennion1997}, Section 3.
Moreover, the law of $Y_i$ is given by
\[
\eta^\ast_\theta(E) := \int_E r_\theta^\ast(x) l^\ast_\theta(dx), \quad E \in {\cal B}({\mathbb S}^{d-1}_+),
\]
where $r^\ast_\theta$ and $l^\ast_\theta$ are given as in Lemma \ref{prop:ip} (cf.\ \shortciteN{Guivarch2012}, Theorem 3.2;
 \shortciteN{Buraczewski.etal:2014}, Proposition 3.1).

\medskip

\textbf{\textit{The condition \boldmath{$({\mathfrak K})$}.}}
Next recall that under $(H_1)$, the measure $\mu_M$ is non-arithmetic
and hence $M_n \cdots M_1$ is positively regular for sufficiently large $n$ w.p.1,
implying that for some positive integer $k$ and some $s>0$,
\begin{equation} \label{requirement on k}  
M_{k} \cdots M_{2} Q_1 \succ s {\vec 1}   \quad \mbox{with positive probability.}
\tag{${\mathfrak K}$}
\end{equation} 
Now if $k>1$, then it is natural to introduce the $k$-step process; namely, for all $k \in \pintegers$, set
\[
\widehat{M}_n := M_{kn} \cdots M_{k(n-1) + 1} \quad \mbox{\text and} \quad  \widehat{Q}_n = \sum_{i=k(n-1)+1}^{kn} M_{kn} \cdots M_{i+1} Q_i,
\]
and note as a consequence of these definitions that
\[
V_{kn} = \widehat{M}_n V_{k(n-1)} + \widehat{Q}_n, \quad n=1,2,\ldots,
\]
where $\widehat{Q}_n -  s {\vec 1} \succ 0 $ with positive probability.
It is worth observing here that the stationary distributions of $\{ V_{kn} \}$ and $\{ V_n \}$ are, of course, identical.

\medskip

Finally, denote by ${\cal C}_{0}\big(\realsd_+\setminus\{0\}\big)$ the set of bounded continuous functions on $\realsd_+\setminus\{0\}$ which are supported on $\R^d_+ \setminus B_r(0)$, for some $r>0$.

We are now prepared to state our first main result.

\begin{thm} \label{thm2A}  Assume that Hypotheses $(H_1)$ and $(H_2)$ are satisfied, and suppose that ${\mathbb D} =  B_r(0) \cap \left(\realsd_+ \setminus \{0\} \right)$,
where $r$ is sufficiently large
such that $\pi({\mathbb D}) > 0$.  Then if $f \in {\cal C}_{0}\big(\realsd_+\setminus\{0\}\big)$ 
and $k = 1$ in \eqref{requirement on k},  we have that
\begin{equation} \label{thm2A-eq1}
 \lim_{u \to \infty} u^\alpha {\mathbb E} \left[ f \left(\frac{V}{u}\right) \right]  = \frac{C}{\lambda^\prime(\alpha)} 
 \int_{{\mathbb S}^{d-1}_+ \times \reals} e^{-\alpha s} f(e^s x) l_\alpha(dx) ds,
  \end{equation}
where 
 \begin{equation} \label{thm2A-eq2}
C =  \int_{\mathbb D} r_\alpha(\widetilde{v}) {\mathbb E}^\alpha_{\bfDV} \left[ \left(|v|+ \sum_{i=1}^\infty \frac{\big< Y_i, \widetilde{Q}_i \big>}{\big<Y_i,X_i\big>}
    \: \frac{|Q_i|}{|M_i \cdots M_1 
 \widetilde{v}|} \right)^\alpha
 {\bf 1}_{\{ \tau = \infty\}} \right] \pi(dv).
  \end{equation}
  If $k>1$ in \eqref{requirement on k}, then the theorem still holds, but with respect to the $k$-step chain $\{ V_{kn} \}$ generated by $\{ (\widehat{M}_i, \widehat{Q}_i) \}$, rather than
the 1-step chain $\{ V_n \}$ generated by $\{ (M_i,Q_i)\}$.
\end{thm}

If $\{V_n\}$ is a Harris recurrent chain, then we may always take $k=1$; see Proposition \ref{thm2A-underH3} below.
Moreover, if $Q \succ 0$ with positive probability, then we may again take $k=1$.  Thus, the condition $k=1$ is seen to be an exceedingly weak requirement.

More generally, for the $k$-step chain, note that the stopping time $\tau$ in \eqref{thm2A-eq2} is then
taken with respect to that chain (rather than the 1-step chain), and the drift factor $\lambda^\prime(\alpha)$ in \eqref{thm2A-eq1} must be replaced with
the drift of the $k$-step chain, namely $k \lambda^\prime(\alpha)$; cf.\ Remark \ref{drift.factor} below.

\begin{rem}   \label{remark:hitting times Vn}
A more intuitive description of $C$ is obtained by setting $Z_n = V_n/|M_n \cdots M_1 V_0| $.  
Then in Lemma \ref{lem:convergence Z} below, we will show that $| Z_n| \to  |Z|$ a.s., where $|Z|$ represents the quantity
appearing in \eqref{thm2A-eq2}, i.e.,
\[
 \abs{Z} ~=~|v| + \sum_{i=1}^{\infty} \frac{\skalar{Y_i, \widetilde{Q}_i}}{\skalar{Y_i, X_i}} \frac{\abs{Q_i}}{\abs{M_i \cdots M_1 \widetilde{v}}}.
 \]
Thus, the constant $C$ is obtained by {\it comparing} the growth rate of $|V_n|$ to the growth rate of $|M_n \cdots M_1 \widetilde{v}|$
 in the $\alpha$-shifted measure, i.e.\ in a setting where these processes diverge.
\end{rem}

It is worth observing that there are other equivalent formulations to \eqref{thm2A-eq1}, 
as follows.

\begin{rem} \label{remark-Lalpha}
 Let ${\mathfrak L}_\alpha$ be the measure on $\R^d_+\setminus\{0\}$ defined by the equation
\[
  \int_{{\mathbb S}^{d-1}_+ \times \reals} e^{-\alpha s} f(e^s x) l_\alpha(dx) ds =  \int_{{\mathbb R}^d_+\backslash \{0\}} f(x) {\mathfrak L}_\alpha(dx).
\]
Then \eqref{thm2A-eq1} yields the vague convergence (of measures on $\R^d_+\setminus\{0\}$)
\begin{equation} \label{vague}
u^\alpha {\mathbb P} \left( \frac{V}{u} \in \cdot\right) \stackrel{v}{\longrightarrow} \frac{C}{\lambda^\prime(\alpha)} {\mathfrak L}_\alpha(\cdot)
\qquad \mbox{\rm as} \:\: u \to \infty.
\end{equation}
 [Here, $\Rd_+\setminus\{0\}$ is considered as a subset of the one-point compactification $(\Rd_+ \cup \{\infty\}) \setminus\{0\}$, where sets of the form $(x,\infty)^d$ are relatively compact. The test functions for vague convergence are those $f \in {\cal C}_{0}\big(\realsd_+\setminus\{0\}\big)$ for which $\lim_{|x| \to \infty} f(x):=f(\infty)$ exists.] 
Now for any measurable set $A \subset \realsd_+$ which is bounded away from zero and satisfies ${\mathfrak L}_\alpha(\partial A)=0$,
it follows from the Portmanteau theorem that
\begin{equation} \label{openset}
\lim_{u \to \infty} u^\alpha {\mathbb P} \left( V \in uA\right) = \frac{C}{\lambda^\prime(\alpha)} 
{\mathfrak L}_\alpha(A).
\end{equation}
Thus, in particular, Theorem \ref{thm2A} yields an estimate for ${\mathbb P}\left( V \in uA \right)$ for any open set $A \subset \realsd_+$ which is bounded 
away from the origin.

Furthermore, note that for any $t >0$ and any measurable $E \subset \SPdp$ with $l_\alpha({\partial E})=0$,
 the sets $E^t:=\{ x \in \R^d_+ \; : \; |x|>t, x/|x| \in E\}$   are ${\mathfrak L}_\alpha$-continuous.  Hence,
 for all $E \subset \SPdp$ with $l_\alpha(\partial E)=0$, \begin{equation} 
\lim_{u \to \infty} u^\alpha {\mathbb P} \left( \abs{V}> tu, \, \frac{V}{\abs{V}} \in   E\right) ~=~ \frac{C}{\alpha \lambda'(\alpha)} t^{-\alpha} l_\alpha(E).
\end{equation}
Thus, we infer the weak convergence 
\begin{equation} \label{thm2A-eq3b}
\lim_{u \to \infty} {\mathbb P} \left( \frac{V}{\abs{V}} \in \cdot \:\Big| \:  \abs{V} > u \right) \Rightarrow  l_\alpha(\cdot).
\end{equation}
In fact, it is easily seen that
\eqref{thm2A-eq3b} together with ${\mathbb P} \left( |V|>u \right) \sim \left(C/\alpha \lambda^\prime(\alpha) \right) u^{-\alpha}$
also yields \eqref{thm2A-eq1}; 
  i.e., these formulations are essentially {\it equivalent}.
For further information on multivariate regular variation and vague convergence, see \shortciteN{Resnick:2004}, Section 3.
\end{rem}

\vspace*{.2cm}

We conclude this section by comparing our theorem with some related results in the literature.  As already noted,
in contrast to \shortciteN{HK73}, the identification of the constant $C$ is {\it explicit} and we have also characterized the directional dependence, whereas \shortciteN{HK73}
considers ${\mathbb P} \left( \big< v, V \big> >u \right)$ for vectors $v \in \realsd$, which is a special case of \eqref{thm2A-eq1}.

In the one-dimensional setting, where \eqref{RDE} holds for $(M,Q) \subset (0,\infty) \times [0,\infty)$,  it was shown  in \shortciteN{JCAV13} that
\begin{equation} \label{JCAV-eq1}
{\mathbb P} \left( V > u \right) \sim \frac{C}{\alpha \lambda^\prime(\alpha)} u^{-\alpha} \quad \mbox{\rm as} \quad u \to \infty,
\end{equation}
where $\lambda(\alpha) = {\mathbb E} \left[ M^\alpha \right]$ for $\alpha$  chosen such that $\lambda(\alpha) =1$, and 
\begin{equation*} 
C = \frac{1}{{\mathbb E}_{\pi_{\mathbb D}}[\tau]} {\mathbb E}_{\pi_{\mathbb D}}^{\alpha} \left[ \left(V_0 + \sum_{i=1}^\infty \frac{|Q_i|}{|M_i \cdots M_1|} \right)^\alpha
 {\bf 1}_{\{ \tau = \infty\}}
 \right].
\end{equation*}
After an application of Lemma  \ref{lemma:hitting times Vn} (showing that $\pi(\DD) = \left( {\mathbb E}_{\pi_{\mathbb D}}[\tau] \right)^{-1}$),
this constant is easily seen to have the form described in \eqref{thm2A-eq2}, since in this simplified setting, $\left( \big< Y_i, \widetilde{Q}_i \big>/\big< Y_i, X_i \big> \right) = 1$.
[In \shortciteN{JCAV13}, the stopping time $\tau$ is taken to be a regeneration time of the process $\{ V_n \}$, but the result holds equally well with the return time $\tau$
in place of the regeneration time, utilizing Lemma \ref{regn-cycle} below.]
If it is further assumed that $\{ M_n \} \subset (0,\infty)$ and $\{ Q_n \} \subset (0,\infty)$ are independent, then an alternative characterization of the constant $C$ has also
been given by \shortciteN{NECSOZ09} using entirely different techniques.

\subsection{Extremal estimates for maxima and first passage times}
Our next objective is to study the large exceedance probability over a single cycle emanating from, and then returning
to, a given set ${\mathbb D} \subset \realsd_+\setminus\{0\}$ and, in this way, characterize the distribution of the first passage time
$T_u := \inf \left\{ n \in \pintegers:  |V_n | > u \right\}$ and, more generally,
\[
T_{u}^A := \inf\left\{ n \in \pintegers:  V_n \in u A \right\}, \quad\mbox{\rm where } A \subset \{ x \in \realsd_+:  |x| > 1 \}.
\]
[Equivalently, we could assume that $A \subset \realsd_+$ is supported on $\realsd_+ \setminus B_r(0)$ for some $r > 0$, and the proofs
would still hold with only minor change.]
First impose the following additional requirement on the set $A$.

\begin{defin}  
We say that a set $A \in {\cal B}(\realsd_+)$ is a {\it semi-cone} if 
$x \in \partial A \Longrightarrow \{ tx: t > 1 \} \subset A$; that is, the ray generated by any point on 
the boundary of $A$ is entirely contained within the set $A$.
\end{defin}

Now suppose that $A \subset \left\{ x \in \realsd_+:  |x| > 1 \right\}$ is a semi-cone, and define
\[
d_A(x) = \inf \left\{ t > 1:  tx  \in A \right\}, \quad x \in {\mathbb S}^{d-1}_+,
\]
and
 \begin{align} \label{def-S*}
& S_n^A ~=~ S_n -  \log d_A(X_n) , \quad n=0,1,2,\ldots;  \nonumber\\[.2cm]
&r_\alpha^A (x) ~=~ r_\alpha (x) \big(d_A(x) \big)^\alpha, \quad x \in {\mathbb S}^{d-1}_+; \qquad
 {\mathfrak P}_A = \left\{ x \in \SPdp:  d_A(x) < \infty \right\};
\end{align}
where $\{ S_n \}$ is defined as in \eqref{def-S}.

As a consequence of Kesten's renewal theorem, it will be shown in Lemma \ref{ruin.est} below that, if ${\mathfrak P}_A = {\mathbb S}^{d-1}_+$, then
\begin{equation} \label{Ruin}
{\mathbb P} \left( \left. M_n \cdots M_1 \widetilde{V}_0 \in uA, \:\:\mbox{\rm for some}\:\: n \in \pintegers \right| V_0 = v\right)
\sim r_\alpha(\widetilde{v}) D_A u^{-\alpha} \quad \mbox{as} \quad u \to \infty,
\end{equation}
where 
\begin{equation} \label{D_A}
D_A :=  \int_{{\mathbb S}^{d-1}_+ \times \reals_+} \: 
\frac{e^{-\alpha s}}{r^A_\alpha(x)} \rho^A(dx, ds)
\end{equation} 
for a measure $\rho^A$ which will be specified below in Section \ref{sect:nonlinearMRT}.
Essentially, \eqref{Ruin} is the {\it ruin estimate} for the Markov random walk $\{ (X_n,S_n^A):  n=0,1,\ldots \}$
under the initial state $X_0 =\tilde{v}$,
and $\rho^A$ corresponds to the
stationary excess distribution for this process.  Indeed,
if $A$ is a semi-cone and $d_A$ is continuous, then it follows immediately from the definitions
 that, on the left-hand side of \eqref{Ruin},
\[
M_n \cdots M_1 \widetilde{V}_0 \in uA \Longleftrightarrow e^{S_n} X_n \in uA \Longleftrightarrow S_n^A > \log u - \log d_A(v).
\]
Now if ${\mathfrak P}_A$ is strictly contained in ${\mathbb S}^{d-1}_+$, then  \eqref{Ruin} will still hold and this {\it defines} the constant $D_A$, although the identification of  $D_A$ is
less explicit in this case (i.e., there is no formula equivalent to \eqref{D_A}).  However, $D_A$ can nonetheless
be interpreted as the ruin constant for the Markov random walk; see Section \ref{sect:6} below.

Finally, let $C$ be defined as in \eqref{thm2A-eq2} and set 
\begin{equation} \label{C-v}
C(v) =  r_\alpha(\widetilde{v}) {\mathbb E}^\alpha_{\bfDV} \left[ \left(|v|+ \sum_{i=1}^\infty \frac{\big< Y_i, \widetilde{Q}_i \big>}{\big<Y_i,X_i\big>}
    \: \frac{|Q_i|}{|M_i \cdots M_1 
 \widetilde{v}|} \right)^\alpha
 {\bf 1}_{\{ \tau = \infty\}} \right].
  \end{equation}

\begin{thm} \label{thm2B} 
Suppose that Hypotheses $(H_1)$ and $(H_2)$ are satisfied and that ${\mathbb D} =  B_r(0) \cap \realsd_+$,
where $r$ is sufficiently large
such that $\pi({\mathbb D}) > 0$.  Assume that $A \in {\cal B}(\realsd_+ \setminus B_1(0))$ is a semi-cone
and the function
$d_A$  is continuous.   
Moreover, assume that $k=1$ in \eqref{requirement on k}.
Then for any probability measure $\gamma$ supported on $\Rd_+\setminus\{0\}$,
\begin{equation} \label{thm2B-eq1}
 \lim_{u \to \infty} u^\alpha {\mathbb P} \left( \left. T_u^A < \tau \, \right| \, V_0 \sim \gamma \right)  = D_A  \int_{\mathbb D} C(v) \gamma(dv).
   \end{equation}
Furthermore, the sequence $\{ T_u^A \}$ converges in distribution; more precisely, 
\begin{equation} \label{thm2B-eq2}
\lim_{u \to \infty} {\mathbb P}\left( \frac{T_u^A }{ u^{\alpha}} \le z \, \Big| \, V_0 = v \right) = 1 - e^{-K_A \, z}, \quad z \ge 0,
\end{equation}
where $K_A = C D_A$.
\end{thm}

We emphasize that the boundary of the set $A$ is allowed to have an unbounded distance to the origin, so that $A$ need not intersect
every ray emanating from the origin in $\realsd_+$.  Also, similarly to Theorem \ref{thm2A}, the assumption that $k=1$ is not necessary if $\{V_n\}$ is Harris recurrent, or if $Q \succ 0$ with positive probability.

\begin{rem} \label{RkPerfekt}
  Eq.\ \eqref{thm2B-eq2} generalizes various known results from extreme value theory relating to the recursive sequence
\eqref{RDE}.   For one-dimensional recursions, estimates have previously been given for the distribution of $\max \big\{ V_i:  1 \le i \le  n^{1/\alpha}u \big\}$ as $n \to \infty$; cf.\ \shortciteN{LHHRSRCV89}, Theorem 2.1.; \shortciteN{RP94}.
In the multidimensional setting, the only result of which we are aware is that of 
\shortciteN{RP97}, who studies the componentwise maxima, namely the process 
\[
\left( \max_{1 \le i \le n} V_1, \ldots, \max_{1 \le i \le n} V_d \right) \quad \mbox{as} \quad n \to \infty.
\]
Note that the componentwise maxima need not be achieved simultaneously; hence Perfekt's results do not coincide with ours.
Moreover, in all of these references, additional conditions are assumed
which we do not impose here; in particular, in their formulations it must be assumed that $V_0 \sim \pi$, so that the sequence $\{ V_n \}$ is stationary.

However, a main contribution of our theorem, beyond its generality, is the  {\it explicit} identification of the constants involved, namely
$K_A$, $D_A$, and $C$ (which are in fact computable, especially in the one-dimensional case).  
We emphasize that \eqref{thm2B-eq2} also allows for {\it general}\, sets $A$ and, thus, it suggests how classical extremal
estimates may be naturally extended to multidimensional problems having spatial dependence,
replacing maxima with first passage times.  Cf.\ \shortciteN{ADSKOZ94} for a similar type of estimate in a different
multidimensional setting.
\end{rem}

\begin{rem} \label{extremal index}
As a particular application of the previous theorem, we now determine the extremal index of $\{ | V_n| \}$.
Integrating with respect to the measure $\pi$ in  \eqref{thm2B-eq2},  we obtain 
\[
\lim_{u \to \infty} {\mathbb P}\left( \frac{T_u^A }{ u^{\alpha}} \le z \, \Big| \, V_0 \sim \pi \right) = 1 - e^{-K_A \, z}, \quad z \ge 0.
\]
Now set $A = \{ x:  |x|>1 \}$.  Then it easily follows with $u=n^{1/\alpha} w$ and $z=w^{-\alpha}$ that
\begin{equation} \label{extremal-1}
\lim_{n \to \infty} {\mathbb P} \left( \max_{1 \le i \le n} |V_i| \le n^{1/\alpha} w \, \Big| \, V_0 \sim \pi \right) = e^{-K_A w^{-\alpha}}.
\end{equation}
Moreover, for this choice of $A$, we obtain by Theorem \ref{thm2A} that
\begin{equation} \label{extremal-2}
\lim_{n \to \infty} n {\mathbb P} \left( |V| > n^{1/\alpha}w \right) = \frac{C}{\alpha \lambda^\prime(\alpha)} w^{-\alpha}.
\end{equation}
Then reasoning as in \shortciteN{RLHR88}, Section 2.2, we conclude from \eqref{extremal-1} and \eqref{extremal-2} that the extremal index
of  $\{ | V_n| \}$ is given by
\begin{equation} \label{extremal-3}
\Theta = \alpha \lambda^\prime(\alpha) D_A.
\end{equation}
For a related result in the one-dimensional setting, see \shortciteN{JCAV13}, Proposition 2.2.
\end{rem}

\begin{rem}  
If $\{ V_n \}$ is Harris recurrent, 
then from the proofs of Theorems \ref{thm2A} and \ref{thm2B}, it can be seen that these results also hold
if one replaces $\tau$ with the first regeneration time of $\{ V_n \}$, assuming that regeneration has occurred at time 0.
However, in this case, the constant in \eqref{thm2A-eq2} takes a slightly different form, namely,
 \begin{equation} \label{thm2A-eq2-regeneration}
C = \frac{1}{{\mathbb E}_{\pi_{\mathbb D}}[\tau]} {\mathbb E}^\alpha_{\nu} \left[ r_\alpha(\widetilde{V}_0) \left(|V_0|+ \sum_{i=1}^\infty \frac{\big< Y_i, \widetilde{Q}_i \big>}{\big<Y_i,X_i\big>}
    \: \frac{|Q_i|}{|M_i \cdots M_1 
 \widetilde{V}_0|} \right)^\alpha
 {\bf 1}_{\{ \tau = \infty\}}
 \right],
\end{equation}
where $\nu$ comes from the minorization; that is, $P^k(x,dy) \ge h(x) \nu(dy)$ for a suitable function $h$ and measure $\nu$, where $P$ is the transition kernel of $\{V_n\}$.
\end{rem}

\subsection{The conditional path under a large exceedance and its empirical law} \label{Sec2.4}

Finally, we consider the path behavior of $\{ V_n\}$ prior to a large exceedance; namely, the conditional law $\P{ \,\cdot \, |T_u^A < \tau }$ as $u \to \infty$, where $\tau$ is the return time to a set $\DD=B_r(0) \cap \Rd_+ \setminus\{0\}$ with $\pi(\DD)>0$.  

We first recall the classical problem.
Suppose that $S_n^{\xi}= \xi_1 + \cdots + \xi_n$ is a random walk on $\reals$ with i.i.d.\ innovations $\{ \xi_i \}  \sim G$ and negative mean, and suppose that Cram\`er's condition is statisfied; namely,
\[  
\int_\reals e^{\alpha  x} G(dx) = 1, \quad \mbox {for some } \alpha >0.
\]
Now if one conditions on ${\cal E}:=\{S_n^{\xi} \text{ hits $u$ before $0$}\}$, then it is well known that
 the law of $S_n^{\xi}$ behaves, as $u \to \infty$, like its associated random walk (\shortciteN{WF71}, Section
 XII.6.(d)); that is, as a random walk whose increments have  the $\alpha$-shifted distribution
\[
G^\alpha(dx) := e^{\alpha x} G(dx).
\] 
Similarly, in ruin theory, it is well known that the likely path to ruin follows a random walk with the $\alpha$-shifted distribution; cf.\ \shortciteN{SA82}.
For further information on {\em random walks conditioned to stay nonnegative}, we refer to \shortciteN{Bertoin.Doney:1994} 
(who also consider $\{S_n^{\xi}\}$ conditioned on ${\cal E}^\prime:=\{ S_k \ge 0 \text{ for all } k  \le n\}$, showing
that this leads to yet another limit law as $n \to \infty$).

Thus, by conditioning on $\{ T_u^A < \tau \}$ in our problem, it is natural to expect that the likely path to a large exceedance follows the $\alpha$-shifted distribution for the {\it Markov} random walk
$\{ (X_n,S_n)\}$. 
The purpose of this section is to make this intuition precise, studying, in the first theorem, the convergence  of finite-dimensional distributions of $\{ V_n \}$ under $\{ T_u^A < \tau \}$, and showing that these distributions converge
 to those of the process $\{e^{S_n}X_n\}$ in the $\alpha$-shifted measure.

Since $\{V_n\}$ is an affine recursion, one cannot expect that its behavior will mimic that of $\{e^{S_n}X_n\}$ over the
entire trajectory.  For this reason, we introduce an ``initial" level $\epsilon_u$ with $\epsilon_u=o(u)$ and $\epsilon_u \uparrow \infty$
as $u \to \infty$,
and study the trajectory of $\{V_n\}$ subsequent to its exceedance beyond the level $\epsilon_u$, showing
that the perturbation by the additive components $\{Q_n\}$ become asymptotically negligible.
Moreover, the  asymptotic ``overjump" distribution $V_{T_{\epsilon_u}}$ can be equated
 to the asymptotic overjump distribution of $\{(X_n, S_n)\}$, which we denote by $\rho$ and characterize in Subsection \ref{sect:nonlinearMRT}
 below.

\begin{thm} \label{thm:excursion}
Suppose that Hypotheses $(H_1)$ and $(H_2)$ are satisfied, and assume that the set $A$ is a semi-cone and the function $d_A$ is 
continuous and bounded on $\SPdp$. 
Let $m \in \pintegers$ and  $g : (\R^d_+)^{m+1} \to \R$ be $\theta$-H\"{o}lder continuous for some $\theta \le \min\{1,\alpha\}$.
Set 
\[
I_u = T_{\epsilon_u}, \quad\mbox{\rm where}\quad \epsilon_u = o(u) \quad \mbox{\rm and} \quad \epsilon_u \nearrow \infty \quad \mbox{\rm as}\quad u \to \infty.
\]
Then for all $v \in \R_+^d$, 
\begin{align}
 \lim_{u \to \infty} &\, \E_v\bigg[ g \bigg( \frac{V_{I_u}}{|V_{I_u}|}, \ldots, \frac{V_{I_u+m}}{|V_{I_u}|} \bigg) \, \bigg| \, 
    T_u^A < \tau \bigg] \nonumber\\[.1cm]
 & \hspace*{1.5cm} =~ \int_{{\mathbb S}^{d-1}_+ \times \reals_+}
 \E_x^\alpha \bigg[ g \big(X_0, e^{S_1} X_1, \ldots, e^{S_m}X_m)\bigg] \rho(dx,ds).
 \end{align}
 \end{thm}

 The class of $\theta$-H\"older continuous functions is a separating class, and thus we deduce, for all $m \in \pintegers$, the weak convergence  
 $$\Prob\left( \bigg( \frac{V_{I_u}}{|V_{I_u}|}, \ldots, \frac{V_{I_u+m}}{|V_{I_u}|} \bigg) \in \cdot \, \bigg| \, T_u^A < \tau \right) ~\Rightarrow~ \int_{{\mathbb S}^{d-1}_+ \times \reals_+} \Prob^\alpha_x \left(  \big(X_0, e^{S_1} X_1, \ldots, e^{S_m}X_m) \in \cdot \right) \rho(dx,ds),$$
for any initial distribution of $V_0$.
 
Note that if we could take $I_u=0$ in the above theorem, then we would obtain an asymptotic description for all paths of finite
length, and by the Kolmogorov extension theorem, we could then conclude---as is obtained in \shortciteN{WF71} or \shortciteN{Bertoin.Doney:1994}---that the
conditional path follows the $\alpha$-shifted distribution.  However, we cannot expect so strong a result here, since, as already noted, the process $\{ V_n \}$ 
is not homogeneous and---as a nonlinear renewal process---only resembles the $\alpha$-shifted distribution for sufficiently large $n$ (e.g., for $n \ge I_u$).

A stronger version of this theorem---allowing for paths of infinite length---will be proved below in Theorem \ref{prop8.1}.
 In this general setting, it is natural to consider a scaled process, normalized by a factor $a^n$ subsequent to time $I_u$.
This normalization is needed, since the distance between $\{V_n\}$ and $\{e^{S_n}X_n\}$, in fact, diverges in Theorem \ref{thm:excursion} as $m \to \infty$.
This type of result is consistent with related conditioned limit theorems from large deviation theory; cf.\ \shortciteN{ADOZ98}, Chapter 7 and references
therein.

Nevertheless, in our final result, we consider the complete path between time zero and a large excursion
terminating at time $T_u^A$. 
Specifically, we prove that the empirical law of the increments $\left\{\log |V_n| - \log |V_{n-1}| \right\}$ along a large excursion has the same limit law, as $u \to \infty$, as $\left\{S_n-S_{n-1} \right\}$ under ${\mathbb P}^\alpha$; cf.\ Lemma \ref{lem:ergodic theorem}.   

\begin{thm} \label{thm:empirical law}
Suppose that Hypotheses $(H_1)$ and $(H_2)$ are satisfied, and assume that the set $A$ is a semi-cone
and the function  $d_A$ is continuous and bounded on $\SPdp$.
Then for any $v \in \R_+^d$ and any bounded Lipschitz continuous function $g: \R \to \R$,
 \begin{equation}\label{empirical distribution function}
 \lim_{u \to \infty} \, \E_v \left[ \left. \bigg| \frac{1}{T_u^A} \sum_{k=1}^{T_u^A} g\left( \log \bigg(  \frac{|V_k|}{|V_{k-1}|} \bigg) \right)  ~-~ \Estat \left[g( S_1) \right] \bigg|  \, \right|   T_u^A < \tau \right] ~=~ 0.
 \end{equation}
\end{thm}

Thus the empirical law of $\left\{ \left(\log|V_k| - \log|V_{k-1}| \right) \right\}$ converges weakly, in $\Prob_v(\cdot | T_u^A < \tau)$-probability, to 
$\widehat{\Prob^\alpha}\big(S_1 \in \cdot\big)$.
  
 Finally we remark that---under a different formulation from ours---conditioned limit theorems have also been studied recently in
  \shortciteN{AJJS14}.  In contrast, they consider path behavior conditioned on a large {\it initial} value,
  whereas we study {\it stopped} processes and obtain an entirely distinct characterization in terms of the $\alpha$-shifted 
  Markov random walk.

\subsection{Structure of the paper}
The  organization of the rest of the paper is as follows. Section \ref{sect:3} contains results about the processes $\{V_n\}$ and $\{(X_n, S_n)\}$ which will be needed for the proofs of the main theorems.
Specifically, in Subsection \ref{subsect:3.1} we quantify recurrence properties of the Markov chain $\{V_n\}$. Subsection \ref{sect:convergenceZ} contains an essential result characterizing the asymptotic ratio 
$|V_n|/e^{S_n}$, which then forms the basis for our application of Melfi's nonlinear Markov renewal theory in Subsection \ref{sect:nonlinearMRT}. In Section \ref{sect:4}, we provide a precise
description of how the distribution of the post-$T_u$-process $\{V_{T_u+k} \, : 0 \le k \le m\}$ relates to that of  $\{e^{S_k}{X_k} \, : \, 0 \le k \le m\}$, for any finite $m$. This characterization 
is then used in the proofs of both Theorem \ref{thm2A} and Theorem \ref{thm2B}, given in Sections \ref{sect:5} and \ref{sect:6}, respectively.  Finally, in Section \ref{sect:7}, we prove a stronger version
of Theorem \ref{thm:excursion} and provide the proof of Theorem \ref{thm:empirical law}.

%
%
%
%

\setcounter{equation}{0}
\section{Background:  Markov chain theory and Markov nonlinear renewal theory} \label{sect:3}

In the present section, we present several results from Markov chain theory and nonlinear renewal theory which will be needed for the proofs
of the main theorems.   
After a quick review of the primary results of this section (most importantly, Lemma \ref{lem:geometric return}, Hypothesis $(H_3)$, Lemma \ref{lem:convergence Z}, Theorem \ref{prop:asymptotic overjump distribution}, and Theorem \ref{thm:Melfi}), the reader may wish to proceed to Sections \ref{sect:4}-\ref{sect:7}, where the main results
of the paper are proved, referring back to Section \ref{sect:3}  as necessary.

\subsection{Markov chain theory for \boldmath{$\{V_n\}$}}\label{subsect:3.1}

Recall that $\pi$ denotes the stationary distribution of $\{V_n\}$. In this subsection, we show that the return times to $\pi$-positive sets which are neighborhoods of the origin
have exponential moments. We further introduce the supplementary Hypothesis $(H_3)$, which will be used in some proofs in the initial steps, 
although this hypothesis will ultimately be removed in the proofs of our main theorems.

\begin{lemma} \label{lm-drift}
Assume that $(H_1)$ and $(H_2)$ are satisfied.
Then for any  $0 < \theta < \min\{1,\alpha\}$, 
there exist positive constants $t <1$ and $L < \infty$ such that  for ${\DD}^\dagger := \{v \in \realsd_+:  |v| \le L \}$,
\begin{equation} \label{drift.2}
{\mathbb E} \left[ \left. |V_n|^\theta r_\theta(\widetilde{V}_n) \: \right| \: {\cal F}_{n-1} \right] \le t |V_{n-1}|^\theta r_\theta (\wt{V}_{n-1}), \quad \mbox{\rm for all} \:\: V_{n-1} \in 
\realsd_+ \setminus {\DD}^\dagger.
\end{equation}
In particular, for $\tau^\dagger := \inf \{ n \in \pintegers \, : \, V_n \in {\DD}^\dagger\}$,  there exists a finite 
constant $B$ such that
\begin{equation} \label{drift.2a}
{\mathbb E} \left[ \left. |V_n|^\theta  {\bf 1}_{\{ \tau^\dagger > n \}} \right| V_0=v \right]  \le B  t^n |v|^\theta, \quad \mbox{\rm for all} \:\: v \in 
\realsd_+ \setminus {\DD}^\dagger.
\end{equation}
\end{lemma}

Note that \eqref{drift.2} can be viewed as an extension of a standard drift condition from Markov chain theory, typically used to assure that the chain is
geometrically recurrent under the additional assumption of $\psi$-irreducibility (cf. \shortciteN{EN84}, Chapter 5, or \shortciteN{SMRT93}, Section 14.2).

\vspace*{.3cm}

\noindent
{\bf Proof.}   From Lemma \ref{prop:ip} (specifically, \eqref{eq:relation rs and lstar} with $(r_\theta, l_\theta^\ast)$ in place of $(r_\theta^\ast,l_\theta)$),
we have that for some constant $c \in (0,\infty)$,
\begin{equation} \label{lm3.1.0}
|V_n|^\theta r_\theta(\widetilde{V}_n) = c | V_n|^\theta \int_{{\mathbb S}_+^{d-1}} \big<y, \widetilde{V}_n \big>^\theta l_\theta^\ast(dy) 
  = c \int_{{\mathbb S}_+^{d-1}} \big< y, V_n \big>^\theta l^\ast_\theta (dy),
\end{equation}
where $V_n = M_n V_{n-1} + Q_n$.  Then it follows by subadditivity that for any $\theta \in (0,1)$, 
\begin{equation} \label{lm3.1.1}
{\mathbb E} \left[ \left. |V_n|^\theta r_\theta(\wt{V}_n) \: \right| \: {\cal F}_{n-1} \right]  \le 
c \, {\mathbb E}  \left[  \int_{{\mathbb S}^{d-1}_+} \left( \big< y, M_n V_{n-1} \big>^\theta + \big< y, Q_n \big>^\theta \right) l^\ast_\theta(dy) \bigg| V_{n-1} \right].
\end{equation}
To identify the quantity on the right-hand side, first apply \eqref{lm3.1.0} once more to obtain that
\begin{align} \label{lm3.1.2}
c \, {\mathbb E}  \left[ \left. \int_{{\mathbb S}^{d-1}_+} \big< y, M_n V_{n-1} \big>^\theta l^\ast_\theta(dy) \right| V_{n-1} \right]  
& = {\mathbb E} \left[ \left. |M_n V_{n-1} |^\theta r_\theta \big( (M_n V_{n-1})^\sim \big) \right| V_{n-1} \right] \nonumber\\[.2cm]
& = |V_{n-1}|^\theta \int |{\mathfrak m}(V_{n-1})^\sim|^\theta r_\theta \Big( ({\mathfrak m} V_{n-1})^\sim \Big) \mu_M(d{\mathfrak m}) \nonumber\\[.2cm]
& =  |V_{n-1}|^\theta \cdot P_\theta r_\theta(V_{n-1}) =  |V_{n-1}|^\theta \cdot \lambda(\theta) r_\theta(V_{n-1}),
\end{align}
where the operator $P_\theta$ was defined just prior to Lemma \ref{prop:ip}, and, by that lemma, $P_\theta$ has eigenvalue $\lambda(\theta)$ and the right-invariant
function $r_\theta$.    Moreover, since $l_\theta^\ast$ is a probability measure and $Q_n$ is independent of $V_{n-1}$,
\begin{equation} \label{lm3.1.3}
{\mathbb E}  \left[ \int_{{\mathbb S}^{d-1}_+}  \big< y, Q_n \big>^\theta l^\ast_\theta(dy) \bigg| V_{n-1} \right] \le {\mathbb E} \left[ |Q_n|^\theta \right].
\end{equation}
Then substituting \eqref{lm3.1.2} and \eqref{lm3.1.3} in \eqref{lm3.1.1} yields
\begin{equation} \label{lm3.1.4}
{\mathbb E} \left[ \left. |V_n|^\theta r_\theta(\wt{V}_n) \: \right| \: {\cal F}_{n-1} \right]  \le \lambda(\theta) |V_{n-1}|^\theta r_\theta(\widetilde{V}_{n-1}) + c \, {\mathbb E} \left[ |Q_n|^\theta \right].
\end{equation}

Now by Hypothesis $(H_2)$, $\E \big[ \abs{Q_1}^\theta \big] < \infty$ for $0 \le\theta \le \alpha$.
Moreover, $0 < \theta < \min\{\alpha, 1 \} \Longrightarrow \lambda( \theta) < 1$.  Choosing $t \in (\lambda(\theta),1)$
and then choosing $L < \infty$ sufficiently large, we conclude by \eqref{lm3.1.4} that \eqref{drift.2} is satisfied.

To obtain \eqref{drift.2a},  iterate \eqref{drift.2} to deduce that
\begin{equation*}
{\mathbb E} \left[ \left. |V_n|^\theta r_\theta(\widetilde{V}_n)  {\bf 1}_{\{ \tau^\dagger > n \}} \right| V_0 \right]  \le t^n |V_0|^\theta r_\theta(\widetilde{V}_0),
\end{equation*}
and use that the function $r_\theta$ is bounded from above and below, 
by Lemma  \ref{prop:ip}. \halmos

Next, recall that $\pi$ denotes
the stationary measure of $\{ V_n \}$. 

\begin{lemma}\label{lem:geometric return}
Suppose $(H_1)$ and $(H_2)$ are satisfied, and let $\DD=B_r(0) \cap \R^d_+ \setminus\{0\}$ for some $r>0$  such that $\pi(\DD) >0$. Let $\DD^\dagger=\{ v \in \R^d_+ \, : \, |v| \le L\}$, where $L$ is chosen such that \eqref{drift.2a} is satisfied and $\DD^\dagger \supset \DD.$ Then  there exist constants $t \in (0,1)$ and $B< \infty$ such that, for $\tau=\inf \{ n \in \pintegers \, : \, V_n \in \DD\}$,
\begin{equation}\label{eq:geometric return}
\sup_{v \in {\DD}^\dagger} {\mathbb P} \left( \left. \tau > n  \, \right| \,  V_0 = v \right)\le B t^{n},\quad \mbox{\rm for all } n \in \pintegers.
\end{equation}
\end{lemma}

\noindent{\bf Proof.}
As in the previous lemma, let ${\tau}^\dagger$ denote the first return time of ${\DD}^\dagger$. Then \eqref{drift.2a} gives for all $v \in {\DD}^\dagger$ and $n \ge 1 $ that
\begin{align*}
{\mathbb E} & \left[ \left. |V_n|^\theta  {\bf 1}_{\{\tau^\dagger > n \}} \right| V_0=v \right]  ~=~  {\mathbb E} \bigg[ \E \Big[ \left.  |V_n|^\theta  {\bf 1}_{\{\tau^\dagger > n \}} \Big| V_1 \Big] \, \1[{\{V_1 \notin {\DD}^\dagger\}}] \right|  V_0=v \bigg] \\[.2cm]
  & ~\le~ B_1 t^{n-1} \E \bigg[ |V_1|^\theta \1[{\{V_1 \notin {\DD^\dagger}\}}]   \, \bigg| \, V_0 = v \bigg] 
  ~\le~ B_1 t^{n-1} \bigg( \E\|M\|^\theta L^\theta + \E|Q|^\theta \bigg).
\end{align*}
Since 
$|V_n| > L$ on $\{ \tau^\dagger > n \}$, it follows that for some finite constant $B_2$,
\begin{equation}\label{eq:return Dtilde}
\sup_{v \in {\DD}^\dagger} {\mathbb P} \left( \left. \tau^\dagger > n \right| V_0 = v \right) \le B_2 t^{n}.
\end{equation}
Below we shall prove that for some constant $s>0$ and $k \in \pintegers$,
\begin{equation}\label{eq:return D} \sup_{v \in {\DD}^\dagger} \P{ \tau > k \, | \, V_0=v} \le (1-s).
\end{equation}

Now assume that \eqref{eq:return Dtilde} and \eqref{eq:return D} hold.  Without loss of generality, we may further assume that
$k=1$, for if $k > 1$, then we may consider the $k$-step chain $\{ V_{nk}:  n=0,1,\ldots \}$ instead of $\{ V_n\}$, and note that if 
$\{ V_{kn} \}$ returns to $\DD$ at a geometric rate, then so does $\{ V_n \}$.  

Thus, assume that $k=1$ in \eqref{eq:return D}, and observe
from (3.9)  that $\tau^\dagger$ has exponential moments; in particular, there exists a constant $\epsilon >0$ such that 
\[
\sup_{v \in {\DD}^\dagger} {\mathbb E}_v\big[ e^{\epsilon \tau^\dagger} \big] < \frac{1}{1-s}.
\]
Let $\tau_1^\dagger,\tau_2^\dagger,\ldots$ denote the successive
returns of $\{ V_n \}$ to ${\mathbb D}^\dagger$. 
Now if $N$ denotes the random number of returns to ${\DD}^\dagger$ prior to time $\tau$, then
\[
\sup_{v \in {\DD}^\dagger} {\mathbb E}_v \big[ e^{\epsilon \tau} \big] \le \sum_j \sup_{v \in {\DD}^\dagger} {\mathbb E}_v \left[ e^{\epsilon(\tau_1^\dagger + \cdots + \tau^\dagger_j + 1)} \right] \P{N=j} \le e^{\epsilon}
\sum_{j} \left( \sup_{v \in {\DD}^\dagger} \E{\Big[e^{\epsilon \tau^\dagger} \Big]} \right)^j (1-s)^j < \infty,
\]
and \eqref{eq:geometric return} follows.

To establish \eqref{eq:return D}, we use Proposition 4.3.1 of \citeN{BDMbook}, which gives a precise description of ${\rm supp} \, \pi$. Namely, there 
exists a set ${\cal S}$ with $\overline{\cal S}={\rm supp} \, \pi$.
Moreover, for each $v_0 \in {\cal S}$, there exists $l \in \pintegers$ and ${\mathfrak m}_1, \ldots, {\mathfrak m}_l
\in \supp \mu_M$, \,$q_1,\ldots,q_l \in \supp \mu_Q$ such that 
$$ h : v \mapsto {\mathfrak m}_l \cdots {\mathfrak m}_1 v + \sum_{i=1}^l {\mathfrak m}_l \cdots {\mathfrak m}_{i+1} q_i$$ is a contraction on $\R^d_+$ with $v_0$ as the unique fixed point. Hence, using that $\DD^\dagger$ is compact, we obtain that for any $\delta > 0$, there exists $j \in \pintegers$ such that $|h^j(v) - v_0|< \delta/2$ for all $v \in {\DD}^\dagger$. Then, from continuity and the definition of the support, we conclude that
\begin{equation}\label{eq:contraction}
\inf_{v \in {\DD}^\dagger} \P{ |V_{lj} - v_0| < \delta \, | \, V_0=v} > 0. \end{equation}
Since $\DD$ is open and $\pi(\DD)>0$, and hence $\DD \cap {\rm \supp} \ \pi \not= \emptyset$, it follows that
 $\DD \cap {\cal S} \neq \emptyset$ as well. Now let $v_0 \in \DD \cap {\cal S}$ and choose $\delta >0$ such that $B_\delta(v_0) \in \DD$. Then \eqref{eq:return D} follows from \eqref{eq:contraction} with $k=lj$.
\halmos

From the previous result we infer that $\{V_n\}$ returns to $\DD$ at a geometric rate, starting from a state in $\DD^\dagger \supset \DD$. In the next result, we calculate the expected return time, now starting from the stationary distribution restricted to $\DD$, and provide a law of large numbers for the return times. First let
$\kappa_0 = 0$ and
\[
\kappa_i = \inf\left\{ n > \kappa_{i-1}:  V_n \in {\mathbb D} \right\}, \quad i=1,2,\ldots;
\]
and let $\tau_i := \kappa_i - \kappa_{i-1}$ denote the inter-return times, $i=1,2,\ldots$.  Set $
N_{\mathbb D}(n) = \sum_{k=1}^n {\bf 1}_{\mathbb D}(V_k).$

\begin{lemma}\label{lemma:hitting times Vn}
Suppose Hypotheses $(H_1)$ and $(H_2)$ are satisfied and $\pi({\mathbb D}) > 0$.
Then for $\pi$-a.e.\ $v \in \R^d_+$, 
\begin{equation} \label{sec7.erg.1}
\lim_{i \to \infty} \frac{\kappa_i}{i} ~=~
\lim_{n \to \infty}\left( \frac{N_{\mathbb D}(n)}{n} \right)^{-1}
~=~\frac{1}{\pi({\mathbb D})} \qquad {\mathbb P}_v\text{\rm -a.s.}
\end{equation}
and $\pi_{\mathbb D}(\cdot): = \pi(\cdot)/\pi({\mathbb D})$ is invariant with respect to the process $\{ V_{\kappa_i}:  i=0,1,\ldots \}$.  
Moreover, as an alternative representation to \eqref{sec7.erg.1}, we also have that
\begin{equation} \label{sec7.erg.2}
\lim_{i \to \infty} \frac{\kappa_i}{i} ~=~ {\mathbb E}_{\pi_{\mathbb D}}[\tau]  \qquad {\mathbb P}_{\pi_{\mathbb D}}\text{\rm -a.s.}
\end{equation}
\end{lemma}

\noindent
[On the right-hand sides of \eqref{sec7.erg.1} and \eqref{sec7.erg.2}, we recall that ${\mathbb P}_v$-a.s., ${\mathbb P}_{\pi_{\mathbb D}}$-a.s.
mean that these results holds a.s.\ provided that the initial value is $V_0=v$, or the initial distribution is $V_0 \sim \pi_{\mathbb D}$, respectively.]\\[-.3cm]

\noindent
{\bf Proof.}  
If $V_0 \sim \pi$, then the sequence $\{ V_n \}$ is stationary.  Moreover, since each $V_n$ is a function of the ergodic sequence $\{ (M_n,Q_n) \}$, it follows from
Proposition 6.31 of \shortciteN{LB68} that $\{ V_n \}$ is ergodic.
Hence,  by Birkhoff's ergodic theorem, we have for any $\pi$-integrable measurable function $h$ that
\begin{equation}\label{eq:Vn ergodic}
 \lim_{n \to \infty} \frac{1}{n} \sum_{k=0}^{n} h(V_k) ~=~\int_{\realsd_+} h(x) \pi(dx) \qquad {\mathbb P}_v\text{\rm -a.s.},
\end{equation}
for $\pi$-a.e.\  $v \in \R^d_+ \setminus \{ 0 \}$.  Setting
$h=\1[{\mathbb D}]$ then yields
\begin{equation}\label{eq:visits to A} 
\lim_{n \to \infty} \frac{N_{\mathbb D}(n)}{n} ~=~ \lim_{n \to \infty} \frac{1}{n} \sum_{k=0}^{n} \1[{\mathbb D}](V_k) ~=~\pi({\mathbb D}) \qquad {\mathbb P}_v\text{\rm -a.s.}
\end{equation}
As $\kappa_0,\kappa_1,\ldots$ denote the successive return times to ${\mathbb D}$, it follows by definition that
 $N_{\mathbb D}(\kappa_1)=1$, $N_{\mathbb D}(\kappa_2)=2,$ and so on.  Thus
$$ 1 =  \frac{N_{\mathbb D}(\kappa_i)}{\kappa_i} \cdot \frac{\kappa_i }{i}.$$
Noting that $\kappa_{i+1} \ge i$ and applying \eqref{eq:visits to A} along the subsequence  $\{\kappa_i\}$, we then obtain that
\begin{equation} \label{sec7.erg.10}
\lim_{i \to \infty} \frac{\kappa_i}{i}=\lim_{n \to \infty} \left( \frac{N_{\mathbb D}(n)}{n} \right)^{-1}
= \frac{1}{\pi({\mathbb D})} \qquad {\mathbb P}_v \text{\rm -a.s.},
\end{equation}
which is \eqref{sec7.erg.1}.
In particular, this proves
 the recurrence of $\DD$, and hence we may apply Proposition VII.3.4 of \shortciteN{SA03} to infer that $\pi_{\mathbb D}$ is an invariant probability measure with respect to the process $\{ V_{\kappa_i} \}$.

Finally suppose that $V_0 \sim \pi_{\mathbb D}$.  Then $\tau_i \equiv \kappa_i - \kappa_{i-1}$ is stationary and by the ergodic theorem,
\begin{equation}  \label{sec7.erg.12}
\lim_{i \to \infty} \frac{\kappa_i}{i} = {\mathbb E}_{\pi_{\mathbb D}}[\tau] \qquad {\mathbb P}_{\pi_{\mathbb D}}\text{\rm -a.s.},
\end{equation}
using that the left-hand side of this equation converges a.s.\ to a deterministic limit, by \eqref{sec7.erg.10}.
\halmos

Now let $P$ denote the transition kernel of $\{ V_n \}$.
We conclude this section with two results which hold under the following  additional Hypothesis $(H_3)$.  \\[-.1cm] 

\noindent
{\bf Hypothesis \boldmath{$(H_3)$}.} 
(i)  There exists a $\pi$-positive set $F$ such that, for each $v \in F$, $P(v, \cdot)$ has an absolutely continuous component with respect to
some $\sigma$-finite non-null measure $\Phi$.

(ii)  $(\supp \: \pi)^\circ \not= \emptyset$.

\begin{lemma}  \label{georecurrent}
Assume that $(H_1)$, $(H_2)$, and $(H_3)$ are satisfied.  Then $\{V_n\}$ is an aperiodic, positive Harris chain on $\realsd_+$.
Moreover, $\{ V_n \}$ is $\psi$-irreducible, regular, and 
geometrically recurrent. 
\end{lemma}

\noindent
{\bf Proof.}
Under $(H_3)$, it follows from \shortciteN{GA03}, Theorem 2.1 (b) and 
Theorem 2.2 (b) that $\{V_n\}$ is an aperiodic, positive Harris chain on $\realsd_+$. Hence by \shortciteN{SMRT93}, Theorem 13.0.1,
\[
\sup_{E \in {\cal B}(\realsd_+)} \left| P^n(x,E) - \pi(E) \right| \to 0 \quad \mbox{as} \quad n \to \infty, 
\]
for all $x \in \realsd_+$. This implies, in particular, that the chain is $\pi$-irreducible (and hence $\psi$-irreducible for some maximal irreducibility measure $\psi$).

Since $(\supp \, \pi)^\circ \not= \emptyset$,  we also have that every compact set with positive invariant measure is petite (\shortciteN{ENPT82}, Remark 2.7 or \shortciteN{SMRT93}, Proposition 6.2.8).
Let $L$ be chosen sufficiently large such that ${\mathbb D}^\dagger:= \{ |v| \le L \}$ has positive invariant measure and  \eqref{drift.2} holds. 
By Lemma \ref{lem:geometric return}, we have 
$\sup_{v \in {\mathbb D}^\dagger} {\mathbb E} \left[ \tau^\dagger | V_0 = v \right] < \infty$ for the return time $\tau^\dagger$ of $\DD^\dagger$.
Then by \shortciteN{SMRT93}, Theorem 11.3.15, we conclude that $\{ V_n \}$ is regular.   Moreover, from the above calculation
given in the proof of Lemma \ref{lem:geometric return},
\[
\sup_{v \in {\mathbb D}^\dagger} {\mathbb E} \Big[ e^{\epsilon \tau^\dagger} \Big] < \infty, \quad \mbox{\rm some } \epsilon > 0,
\]
and hence $\{ V_n \}$ is geometrically recurrent.
\halmos

Using the $\psi$-irreducibility from the previous lemma, we may observe the following useful result, connecting the stationary distribution of $\{V_n\}$ to 
its average behavior over a given cycle emanating from a $\pi$-positive set ${\mathbb D}$ with initial measure $\pi_{\mathbb D}(\cdot) :=  \pi(\cdot \cap \DD)/\pi(\DD)$.

\begin{lemma}  \label{regn-cycle} Suppose that $(H_1)$\,--\,$(H_3)$ are satisfied, and let $\DD \subset  \R^d_+ \setminus\{0\}$ be 
chosen such that $\pi(\DD) >0$. Let $\tau :=\inf \{n \in \pintegers \, : \, V_n \in \DD\}$ denote the first return time of $\DD$. 
Then for any $\pi$-integrable $\mbox{function }h$, 
\begin{equation} \label{regncycle}
\int h(v) \pi(dv) ~=~ \E\big[ h(V) ] ~=~ \frac{1}{\E_{\pi_\DD} \left[ \tau \right] }
      {\mathbb E}_{\pi_\DD} \left[ \sum_{i=0}^{\tau-1} h(V_i) \right] .
\end{equation}
\end{lemma}

\noindent
{\bf Proof.}  See \shortciteN{EN84}, Proposition 5.9 and the discussion
just prior to Corollary 5.3.  For a closely related result, also see the proof of Theorem 2.1 in \shortciteN{JCGDAV14}.
\halmos

%
%
%
%

\subsection{ Quantifying the discrepancy between \boldmath{$\{V_n\}$}
and \boldmath{$\{e^{S_n}X_n \}$} }\label{sect:convergenceZ}
The objective of this subsection is to precisely quantify the discrepancy between $\{V_n\}$ and $\{e^{S_n}X_n\}$,
which, in essence, will later be shown to determine the constant $C$ in Theorem \ref{thm2A}.
To this end, let
\begin{equation} \label{def-Z}
Z_n ~:=~ \frac{V_n}{\abs{M_n \cdots M_1 X_0}} ~=~ \frac{V_n}{e^{S_n}}, \quad n \in {\mathbb N}
\end{equation}
and
$$
Z_n^{(0)} ~:=~ \frac{V_n - M_n \cdots M_1 V_0}{\abs{M_n \cdots M_1 X_0}} ~=~ \frac{\sum_{i=1}^n M_n \cdots M_{i+1} Q_i}{\abs{M_n \cdots M_1 X_0}}, \quad n \in {\mathbb N}.
$$
Also introduce the shorthand notation
$$ \Pi_n:=M_n \cdots M_1, \quad \text{ and } \quad \Pi_{i}^n:=M_n \cdots M_i.$$
The most important properties of $\{ Z_n \}$, for our purposes,
are summarized in the following.

\medskip

\begin{lemma}  \label{lem:convergence Z}
Assume $(H_1)$ and $(H_2)$. Then{\rm :}

{\rm (i)}   $\sup_{n \in \N} \abs{Z_n} < \infty \quad {\mathbb P}^\alpha \text{\rm -a.s.}$ \:{\rm and}\: $\sup_{n \in \N} \big|Z_n^{(0)}\big| < \infty \quad {\mathbb P}^\alpha \text{\rm -a.s.}$ \\[-.3cm]

{\rm (ii)}  Let $v \in \R_+^d \setminus\{0\}$  .  Then in ${\mathbb P}_{\bfDV}^{\alpha}$-measure, the sequence $\{ Z_n \}$ converges in law to a random variable $Z$, 
 and
$| Z_n| \longrightarrow  |Z|$ a.s., where
\begin{equation} \label{jan16-0}
 \abs{Z} ~=~|v| + \sum_{i=1}^{\infty} \frac{\skalar{Y_i, \widetilde{Q}_i}}{\skalar{Y_i, X_i}} \frac{\abs{Q_i}}{\abs{\Pi_i \widetilde{v}}}  \qquad {\mathbb P}_{\bfDV}^\alpha\text{\rm -a.s.}
 \end{equation}
 Moreover,  $| Z |$ is strictly positive and finite ${\mathbb P}_{\bfDV}^\alpha$-a.s.
Similarly, we have 
\begin{equation}\label{eq:convergence Zn0}  \lim_{n \to \infty} \big|Z_n^{(0)}\big| ~=~ \sum_{i=1}^{\infty} \frac{\skalar{Y_i, \widetilde{Q}_i}}{\skalar{Y_i, X_i}} \frac{\abs{Q_i}}{\abs{\Pi_i \widetilde{v}}}  \qquad {\mathbb P}_{\bfDV}^\alpha\text{\rm -a.s.} 
\end{equation}
 
{\rm (iii)} Let $F \subset \R_+^d\setminus\{0\}$ be a bounded set and  let $\tau^\prime$ be any $\{\F_n\}$-stopping time such that 
\begin{equation} \label{geometric.return.assumption}
\sup_{v \in F} \P{ \tau^\prime > k |V_0=v} \le B t^k, \quad \mbox{for all } k \in \N,
\end{equation}
for some finite constant $B$ and $t \in (0,1)$.  Then for any $v \in \R_+^d \setminus\{0\}$,
\begin{equation} 
\sup_{v \in F} \E_{\bfDV}^{\alpha} \left[ \sup_{n \in \N} \abs{Z_n}^\alpha {\bf 1}_{\left\{ \tau^\prime \ge n \right\}} \right] < \infty \quad \text{\rm and } \quad \sup_{v \in F} \E_{\bfDV}^{\alpha} \big[ |Z|^\alpha \1[{\{ \tau^\prime = \infty\}}] \big] < \infty.
\end{equation}
 
 {\rm (iv)}   For $v \in \R_+^d \setminus\{0\}$, we have the $L^1$-convergence
 \begin{equation} 
\lim_{n \to \infty} \E^\alpha_{\bfDV} \Big[ \big| \, |Z_n|^\alpha \1[{\{ \tau^\prime  \ge n\}}] - |Z|^\alpha \1[{\{ \tau^\prime = \infty\}}]  \big| \Big] ~=~0.
\end{equation}
\end{lemma}

Note that by Lemma \ref{lem:geometric return}, the condition in (iii) holds, in particular, for $\tau^\prime = \tau :=\inf\{n \in \pintegers \, : \, V_n \in \DD\}$, namely the return time of $V_n$ into the set $\DD$. 

\medskip

\noindent {\bf Proof.}  
For any vector  $x \in \reals^d$, let $x^{(i)} = \big< e_i,x \big>$ denote the $i^{\rm th}$ component of $x$, and set ${\vec 1} = (1,\ldots,1)^T$. 
Also, except in part (iii), fix $V_0=v$ throughout the proof.

 First recall that any $\Prob^\alpha_{x,v}$ is absolutely continuous with respect to $\Pstat$ (\shortciteN{Buraczewski.etal:2014}, Lemma 6.2), and hence the convergence of $\{ Z_n \}$ in law, or the convergence of
$\{ |Z_n| \}$ $\Pstat$-a.s., implies the respective convergence under ${\mathbb P}^\alpha_{x,v}$. 
Thus, it is sufficient to prove the convergence results in part (i) and (ii)
with respect to the measure  $\Pstat$, under which the sequence $\{(M_n, Q_n): n =1,2,\ldots\}$ is {\it stationary}; cf.\ the discussion in Section 2 above. 
This will allow us to apply the results of  \shortciteN{Hennion1997}.

(i)  Suppose ${\mathfrak m} \in {\mathfrak M}$, and let $x_{\mathfrak m}$ be chosen such that $\| {\mathfrak m} \| = | {\mathfrak m} x_{\mathfrak m}|.$  Since ${\mathfrak m}$ is
nonnegative, an elementary argument shows that  $x_{\mathfrak m}$ can, in fact, be chosen such that $x_{\mathfrak m}^{(i)} \ge 0$ for all $i$.
Then for any $x \in {\mathbb S}^{d-1}_+$, 
\[
\left|{\mathfrak m} x \right| \ge \big( \min_j x^{(j)} \big) \left|{\mathfrak m}{\vec 1} \right| \ge \big( \min_j x^{(j)} \big) \left| {\mathfrak m} x_{\mathfrak m} \right| = \big( \min_j x^{(j)} \big) \| {\mathfrak m} \|.
\]
Thus
\[
\frac{\| {\mathfrak m} \|}{|{\mathfrak m} x |} \le \frac{1}{ \min_j x^{(j)}}, \quad \mbox{\rm for all } x \in {\mathbb S}^{d-1}_+ \mbox{ and all } {\mathfrak m} \in {\mathfrak M}.
\]
 Recall the stopping time
${\mathfrak N} :=  \inf\big\{ n \in {\mathbb Z}_+ : \, \Pi_n \succ 0 \big\},$
which is finite $\Pstat$-a.s. by $(H_1)$.  [Since $\mu_M$ is equivalent to $\Pstat(M_1 \in \cdot)$, $(H_1)$ holds equally well for $\Pstat(M_1 \in \cdot)$.  Then Lemma 3.1 of \shortciteN{Hennion1997} yields the finiteness of ${\mathfrak N}$.]
Identifying $Q_0:=V_0=v$, we obtain
\begin{align} \label{pf-lemmaZ.1}
\abs{Z_n} & \le  \sum_{i=0}^n \frac{\abs{\Pi^n_{i+1} Q_i}}{\abs{\Pi_n X_0}}  \le \sum_{i=0}^n \frac{\norm{\Pi^n_{i+1}} \abs{Q_i}}{\abs{\Pi^n_{i+1} X_i} \abs{\Pi_i X_0}}
  \nonumber\\[.2cm]
&\le \sum_{i=0}^{\mathfrak N \wedge n} \frac{\norm{\Pi_{i+1}^n} \abs{Q_i}}{\abs{\Pi^n_{i+1} X_i} \abs{\Pi_i X_0}} ~+~ \sum_{i=\mathfrak N \wedge n}^n \frac{1}{ \min_j X_i^{(j)}}  \frac{ \abs{Q_i}}{ \abs{\Pi_i X_0}}.
\end{align}
By \shortciteN{Buraczewski.etal:2014}, Lemma 6.3,
$$ C_i(x):= \inf_{n \in \N} \frac{|\Pi^n_{i+1} x|}{\norm{\Pi^n_{i+1}}} > 0 \qquad \Pstat\text{-a.s.},$$
for all $x \in \SPdp$.  Also observe that
\begin{equation} \label{def-coordX}
X_i^{(j)} =  \frac{(\Pi_i X_0)^{(j)}}{|\Pi_i X_0|},
\end{equation}
which implies that $X_i^{(j)} |\Pi_i X_0| = \left(\Pi_i X_0\right)^{(j)} = \skalar{e_j, \Pi_i X_0} $.   This identifies the denominator in the second sum of \eqref{pf-lemmaZ.1},
and shows that this denominator is positive for $i \ge \mathfrak N$. Hence
\begin{align} \label{pf-lemmaZ.1a}
\sup_{n \in \N} \abs{Z_n} &\le \sum_{k=0}^{\mathfrak N} \frac{\abs{Q_i}}{C_i( X_i) \abs{\Pi_i X_0}} ~+~ \sum_{i= \mathfrak N }^\infty  \frac{ \abs{Q_i}}{\min_{j}  \skalar{e_j, \Pi_i X_0}} \nonumber \\[.2cm]
&\le \sum_{k=0}^{\mathfrak N} \frac{\abs{Q_i}}{C_i( X_i) \abs{\Pi_i X_0}} ~+~  \sum_{i= \mathfrak N }^\infty \, \sum_{j=1}^d \, \frac{ \abs{Q_i}}{ \skalar{e_j, \Pi_i X_0}}.
\end{align}

Since $\mathfrak N < \infty$ $\Pstat$-a.s., it suffices to focus on the second sum. By Lemma \ref{lem:slln},
\begin{equation} \label{pf-lemmaZ.2}
\lim_{n \to \infty} \, \sup \bigg\{ \abs{\frac1n \1[{\{ {\mathfrak N} \le n\}}] \log \skalar{y, \Pi_n x} - \Lambda^\prime(\alpha)} \, : \, x,y \in \SPdp   \bigg\} ~=~0 \qquad \Pstat\text{-a.s.}
\end{equation}
Furthermore, by a Borel-Cantelli argument, $\Pstat \left( \log | Q_i | > \delta i \:\: \mbox{\rm i.o.} \right) = 0$, for all $\delta > 0.$
Thus, given $\epsilon \in (0,\Lambda^\prime(\alpha))$, there exists a finite integer $k_0$ such that, for all $i \ge k_0$ and all $j \in \{1,\ldots,d\}$,
\begin{equation} \label{pf-lemmaZ.3}
\log \abs{Q_i} - \log  \skalar{e_j, \Pi_i X_0}~\le~ -\bigg({\Lambda'(\alpha)} - \epsilon\bigg) i \quad \Pstat\text{-a.s.}
\end{equation}
Since \eqref{pf-lemmaZ.3} holds {\it uniformly} in $j$, substituting \eqref{pf-lemmaZ.3} 
into \eqref{pf-lemmaZ.1a} establishes part (i) of the lemma, where we also use that $\big|Z_n^{(0)}\big| \le |Z_n|$ for all $n \in \N$.

(ii)   Following \shortciteN{Hennion1997}, let $\rho(\Pi_k^n)$ denote the spectral radius of $\Pi^n_k$, and let $R_n^k$ and $L_n^k$ denote the
right and left eigenvectors corresponding to the maximal eigenvalue in modulus; that is,
\[ 
(\Pi^n_k) R_n^k = \rho(\Pi_k^n) R_n^k \quad\mbox{and} \quad \left(\Pi_k^n \right)^T L_n^{k} = \rho(\Pi_k^n) L_n^k, \qquad 1 \le k \le n.
\] Note that the Perron-Frobenius theorem assures that $R_n^k$ and $L_n^k$ have nonnegative entries.
We further assume the following normalization:
$$ \big| L_n^k \big| =1, \qquad \skalar{L_n^k, R_n^k}=1, \qquad  1 \le k \le n.$$
Now let $\{ Y_i \}$ be defined as in \eqref{def-upsilon}.
Then we will show that
\begin{equation} \label{jan16-1}
\lim_{n \to \infty} \left| \big<e_j,Z_n \big> - \big<e_j, \widetilde{R}_n^1 \big> \sum_{i=0}^{n} \frac{\skalar{Y_i, \widetilde{Q}_i}}{\skalar{Y_i, X_i}} \frac{\abs{Q_i}}{\abs{\Pi_i X_0}}  \right| = 0 \quad \Pstat\mbox{-a.s.},
\end{equation}
for all $1 \le j \le d$.   The sequence $\{ \widetilde{R}_n^1\}$ converges {\it in distribution} as $n \to \infty$ (\shortciteN{Hennion1997}, Theorem 1 (ii)\:(b));
hence we obtain the convergence, in distribution, of
$\{ Z_n \}$ to
\[
Z :=  \lim_{n\to \infty}  \widetilde{R}_n^1\cdot  \lim_{n \to \infty} \sum_{i=0}^{n} \frac{\skalar{Y_i, \widetilde{Q}_i}}{\skalar{Y_i, X_i}} \frac{\abs{Q_i}}{\abs{\Pi_i X_0}} .
\]
Moreover, since $\big|\widetilde{R}_n^1 \big|=1$, 
 \eqref{jan16-1}
yields \eqref{jan16-0}, i.e.\ $\lim_{n \to \infty} |Z_n|=|Z|$ $\Pstat$-a.s. In the same way, \eqref{eq:convergence Zn0} is obtained by setting $Q_0=0$.

To establish \eqref{jan16-1}, first recall that (with the identification $Q_0 :=V_0=v$) we have that
\[
Z_n = \sum_{i=0}^n \frac{\Pi^n_{i+1} Q_i}{|\Pi_n X_0|},
\]
and observe that
\begin{align} \label{jan16-2}
 \abs{ \sum_{i= \lfloor n/2 \rfloor +1}^{n} \frac{\skalar{e_j, \Pi^n_{i+1} Q_i}}{ \abs{\Pi_n X_0}  } } 
&\le  \sum_{i= \lfloor n/2 \rfloor +1}^\infty  \frac{\norm{\Pi^n_{i+1}}\abs{Q_i}}{\abs{\Pi^n_{i+1} X_i}\abs{\Pi_i X_0}}, \nonumber
\end{align}
and the right-hand side tends to zero as $n \to \infty$ by the proof of part (i), in particular Eq.\ \eqref{pf-lemmaZ.1a}.
Since $Y_i$ is a unit vector with nonnegative entries, $\skalar{Y_i, X_i} \ge d^{-1} \, {\min_j X_i^{(j)}}$.  Hence we also have
 \begin{align*} 
\abs{ \big<e_j, \widetilde{R}_n^1 \big> \sum_{i= \lfloor n/2 \rfloor +1}^n \frac{\skalar{Y_i, \widetilde{Q}_i}}{\skalar{Y_i, X_i}} \frac{\abs{Q_i}}{\abs{\Pi_i X_0}} }~\le~ \sum_{i= \lfloor n/2 \rfloor +1}^n \frac{d}{ \min_j X_i^{(j)}}  \frac{ \abs{Q_i}}{ \abs{\Pi_i X_0}}.
\end{align*} 
Thus, to establish \eqref{jan16-1} and consequently part (ii) of the lemma, it is enough to show that
\begin{equation} \label{jan16-3}
\lim_{n \to \infty} \left|  \sum_{i= 0}^{\lfloor n/2 \rfloor} \frac{\skalar{e_j, \Pi^n_{i+1} Q_i}}{ \abs{\Pi_n X_0}  }- \big<e_j, \widetilde{R}_n^1 \big> \sum_{i=0}^{
\lfloor n/2 \rfloor} \frac{\skalar{Y_i, \widetilde{Q}_i}}{\skalar{Y_i, X_i}} \frac{\abs{Q_i}}{\abs{\Pi_i X_0}}  \right| = 0 \quad \Pstat\mbox{-a.s.}
\end{equation}
Then by the triangle inequality, it is sufficient to show the following.

\begin{sublemma} 
\begin{align} \label{jan16-3a}
\lim_{n \to \infty} & \sum_{i=0}^{\lfloor n/2 \rfloor} \abs{   \frac{\skalar{e_j, \Pi_{i+1}^n Q_i}}{ \abs{\Pi_n X_0}  }   -  \skalar{e_j, \widetilde{R}_n^{i+1}}  \frac{\skalar{L_n^{i+1}, Q_i}}{\skalar{L_n^{i+1}, X_i}} \frac{1}{\abs{\Pi_i X_0}} }  = 0 \quad \Pstat\mbox{\rm -a.s.};\\[.2cm] \label{jan16-3b}
\lim_{n \to \infty} & \sum_{i=0}^{\lfloor n/2 \rfloor} \abs{  \frac{\skalar{e_j, \widetilde{R}_n^{i+1}} }{\abs{\Pi_i X_0}} \left(
 \frac{\skalar{L_n^{i+1}, Q_i}}{\skalar{L_n^{i+1}, X_i}} -  \frac{\skalar{Y_{i+1}, Q_i}}{\skalar{Y_{i+1}, X_i}} \right) } = 0 \quad \Pstat\mbox{\rm -a.s.};
 \end{align}
and 
\begin{equation} \label{jan16-3c}
\lim_{n \to \infty}   \sum_{i=0}^{\lfloor n/2 \rfloor} \abs{ \frac{1}{\abs{\Pi_i X_0}} \frac{\skalar{Y_{i+1}, Q_i}}{\skalar{Y_{i+1}, X_i}}
\left(  \skalar{e_j, \widetilde{R}_n^{i+1}} -  \skalar{e_j, \widetilde{R}_n^{1}} \right) } = 0 \quad \Pstat\mbox{\rm -a.s.}  \hspace*{1.5cm}
\end{equation}
\end{sublemma}

\noindent
{\bf Proof of the Sublemma.}   First we establish \eqref{jan16-3a}.  To this end,
observe by Corollary 1 of \shortciteN{Hennion1997} that
\begin{equation}\label{eq:convergence matrices} 
\lim_{n \to \infty} \bigg( \frac{\Pi^n_{i+1}}{\norm{\Pi^n_{i+1}}} - \frac{ R_n^{i+1} \otimes L_n^{i+1}}{\norm{R_n^{i+1} \otimes L_n^{i+1}}} \bigg)=0 \quad \Pstat\text{-a.s.},
\end{equation} 
where $a \otimes b$ denotes the rank-one matrix with the property that 
\begin{equation} \label{tensor}
\skalar{e_i, (a \otimes b) e_j}=\skalar{e_i, a}\skalar{b, e_j} .
\end{equation}
From \eqref{eq:convergence matrices} we infer that
\begin{equation}  \label{jan16-5}
\lim_{n \to \infty} \left(\frac{\skalar{e_j, \Pi^n_{i+1} Q_i}}{\norm{\Pi^n_{i+1}}} - \frac{\skalar{e_j, R_n^{i+1}}\skalar{L_n^{i+1},Q_i}}{\norm{R_n^{i+1} \otimes L_n^{i+1}}} \right) = 0 
\end{equation}
and
\begin{equation} \label{jan16-6}
\lim_{n \to \infty} \left(\frac{\abs{\Pi^n_{i+1} X_i}}{\norm{\Pi^n_{i+1}}} - \frac{\abs{R_n^{i+1}}\skalar{L_n^{i+1},X_i}}{\norm{R_n^{i+1} \otimes L_n^{i+1}}} \right)=0. 
\end{equation}
Combining  \eqref{jan16-5} and \eqref{jan16-6}, we conclude that
\begin{align} \label{eq:summands} \lim_{n \to \infty}  \frac{\skalar{e_j,\Pi^n_{i+1} Q_i}}{|\Pi_n X_0|} &= \lim_{n \to \infty}  \frac{\skalar{e_j,\Pi^n_{i+1} Q_i}}{|\Pi^n_{i+1} X_i| |\Pi_i X_0|} \nonumber \\
&=  \lim_{n \to \infty} \skalar{e_j, \widetilde{R}_n^{i+1}} \frac{\skalar{L_n^{i+1}, Q_i}}{\skalar{L_n^{i+1}, X_i}} \frac{1}{\abs{\Pi_i X_0}} \quad \Pstat\text{-a.s.},
\end{align}
showing, in particular, that the individual terms in \eqref{jan16-3a} converge to zero $\Pstat$-a.s.

To prove that the sum in \eqref{jan16-3a} converges to zero, we
now  invoke a dominated convergence argument. Since $\mathfrak N$ is finite a.s., it suffices to focus on summands with $i \ge \mathfrak{N}$, where we can assume that all components of $X_i$ are positive, as the remaining terms form a finite sum.  Observe that
\begin{equation} \label{jan16-6a}
\frac{\skalar{ L_n^{i+1}, {\vec 1}} \max_j Q_i^{(j)}}{\skalar{ L_n^{i+1}, {\vec 1} } \min_j X_i^{(j)}} \le \frac{|Q_i|}{\min_j X_i^{(j)}},
\end{equation}
and therefore
\begin{align} \label{jan16-7}
\sup_n \sum_{i=\mathfrak N}^{\lfloor n/2 \rfloor} & \abs{   \frac{\skalar{e_j, \Pi^n_{i+1} Q_i}}{ \abs{\Pi_n X_0}  }  -  \skalar{e_j, \widetilde{R}_n^{i+1}}  \frac{\skalar{L_n^{i+1}, Q_i}}{\skalar{L_n^{i+1}, X_i}}
 \frac{1}{\abs{\Pi_i X_0}} }
 \nonumber\\[.3cm]
& \hspace*{1cm} \le 2 \sup_n \sum_{i=\mathfrak N}^{\lfloor n/2 \rfloor} \frac{1}{ \min_i X_i^{(j)}}  \frac{ \abs{Q_i}}{ \abs{\Pi_i X_0}} < \infty \quad \Pstat \mbox{-a.s.}
\end{align}
by part (i)  (where we have used the calculation in \eqref{pf-lemmaZ.1} to handle the first term on the left-hand side).
Thus, using a dominated convergence argument to interchange summation and limit (applied pointwise on the space where \eqref{jan16-6} and \eqref{jan16-7} hold), we now conclude that
 \eqref{jan16-3a} follows from \eqref{jan16-6}.

Next we turn to \eqref{jan16-3b}.  It follows by Lemma 3.3 of \shortciteN{Hennion1997} that, under $\Pstat$, the sequence $\{ L_n^{i+1} \}$ converges a.s.\ as $n \to \infty$ to $Y_{i+1}$.
Hence, by a dominated convergence argument (analogous to the one just used to establish \eqref{jan16-3a}), we conclude that \eqref{jan16-3b} holds.

Finally, to establish \eqref{jan16-3c}, 
note by Proposition 3.1 of \shortciteN{Hennion1997} that
\begin{equation} \label{eq:Rn} 
\abs{ \widetilde{R}_{n}^{i+1} - \widetilde{R}_n^1} = \abs{(\Pi^n_{i+1} R_n^{i+1})^\sim - (\Pi_n R_n^1)^\sim}
    \le 2 {\bf c} (\Pi^n_{i+1}),  
\end{equation}
where ${\bf c}(\cdot)$ is bounded above by one and tends to zero $\Pstat$-a.s.\ as $n-i$ tends to infinity (\shortciteN{Hennion1997}, Lemma 3.2).
Then \eqref{jan16-3c} follows, once again, by a dominated convergence argument.  This completes the proof of the sublemma and, consequently, part (ii) of Lemma \ref{lem:convergence Z}.
\halmos

\noindent
{\bf Proof of Lemma \ref{lem:convergence Z} (continued).}
We now return to the proof of main lemma, where it remains to verify that (iii) and (iv) hold.

(iii)  Let $m \in \N$, and set $B_1 = \max_{x,y} \big( r_\alpha(x)/r_\alpha(y) \big) \in (0,\infty)$.
Then for $\alpha >0$,
\begin{align*}
 \E^\alpha_{\bfDV} \bigg[ \Big(  \sup_{n \le m} |Z_n|  & \1[{\{\tau^\prime \ge n\}}] \Big)^{\alpha} \bigg] 
 \le~ \E^\alpha_{\bfDV} \bigg[ \Big( \sup_{n \le m} \Big(  |v| + \sum_{k=1}^n \frac{|\Pi^n_{k+1} Q_k|}{|\Pi_n X_0|} \1[{\{\tau^\prime \ge k-1\}}] \Big) \Big)^{\alpha} \bigg] \\
  &=~ \left(r_\alpha(\widetilde{v})\right)^{-1} \E_{\bfDV} \bigg[ r_\alpha(X_m) |\Pi_m X_0|^\alpha \Big( \sup_{n \le m} \Big(|v| + \sum_{k=1}^n \frac{|\Pi^n_{k+1} Q_k|}{|\Pi_n X_0|} \1[{\{\tau^\prime \ge k-1\}}] \Big) \Big)^{\alpha} \bigg] \\
   &\le~ B_1 \E_{\bfDV} \bigg[   \Big( \sup_{n \le m} \Big(|\Pi_m X_0|  |v| + \sum_{k=1}^n |\Pi^m_{n+1} X_n| \cdot |\Pi^n_{k+1} Q_k| \1[{\{\tau^\prime \ge k-1\}}] \Big) \Big)^{\alpha} \bigg] \\
      &\le~ B_1 \E_v \bigg[   \Big( \sup_{n \le m} \Big( \sum_{k=0}^n \norm{\Pi^m_{n+1}} \cdot \norm{\Pi^n_{k+1}} \cdot  |Q_k| \1[{\{\tau^\prime \ge k-1\}}]  \Big) \Big)^{\alpha} \bigg], \quad \mbox{\rm where } Q_0 := v,\\
            &=~ B_1 \E_v \bigg[  \Big(\sum_{k=0}^m \norm{\Pi^m_{n+1}} \cdot \norm{\Pi^n_{k+1}} \cdot  |Q_k| \1[{\{\tau^\prime \ge k-1\}}] \Big)^{\alpha} \bigg]. 
\end{align*}
Now suppose that $\alpha \ge 1$.  Then by Minkowski's inequality,
\begin{align*}
 \bigg(\E_v \bigg[  &\Big( \sum_{k=0}^m \norm{\Pi^m_{n+1}} \cdot \norm{\Pi^n_{k+1}} \cdot  |Q_k| \1[{\{\tau^\prime \ge k-1\}}] \Big)^{\alpha} \bigg] \bigg)^{1/\alpha} \\
\le&~ \sum_{k=0}^m  \left(\E_v \Big[  \norm{\Pi^m_{n+1}}^\alpha \cdot \norm{\Pi^n_{k+1}}^\alpha \cdot  |Q_k|^\alpha \1[{\{\tau^\prime \ge k-1\}}]  \Big] \right)^{1/\alpha} \\
=&~  \sum_{k=0}^m  \left(\E \Big[  \norm{\Pi^m_{n+1}}^\alpha \Big] \right)^{1/\alpha} \left( \E \Big[ \norm{\Pi^n_{k+1}}^\alpha \Big] \right)^{1/\alpha} 
  \Big(|v| + \E \big[ |Q_1|^\alpha \big] \Big)^{1/\alpha} {\mathbb P}_v \left( \tau^\prime \ge k-1\right)^{1/\alpha}.
\end{align*}
Now by Corollary 4.6 of \shortciteN{Buraczewski.etal:2014}, 
$\E \left[ \norm{\Pi_n}^\alpha \right] \le B_2 \in (0,\infty)$, for all $n$.  Moreover, $ \E [\abs{Q_i}^\alpha] < \infty$ by Hypothesis $(H_2)$,
and by our assumption \eqref{geometric.return.assumption},
$\sup_{v \in F} {\mathbb P} \left( \tau^\prime \ge k \, | \, V_0=v\right) \le
B_3 t^{k}$ for some $t \in (0,1)$.
Hence by monotone convergence,
\begin{align} \label{sec3-ADD100}
\E^\alpha_{\bfDV} \bigg[ \sup_{n \in \N} |Z_n|^\alpha \1[{\{\tau^\prime \ge n\}}] \bigg] 
~&=~ \lim_{m \to \infty} \E^\alpha_{\bfDV} \bigg[ \sup_{n \le m} |Z_n|^\alpha \1[{\{\tau^\prime \ge n\}}] \bigg] \nonumber\\[.2cm]
&\le~ B_1 B_2^2 \Big( |v| + \E[ |Q_1|^\alpha] \Big) \left( \sum_{k=0}^\infty (B_3 t^k)^{1/\alpha} \right)^\alpha  < \infty,
\end{align}
and this bound is uniform over $v \in F$ for any bounded set $F \subset \R_+^d\setminus\{0\}$.

If $\alpha \le 1$, then we use the subadditivity, namely the inequality $\abs{x+y}^\alpha \le \abs{x}^\alpha + \abs{y}^\alpha$ in place of Minkowski's inequality, and then proceed as before, obtaining an analogous estimate to \eqref{sec3-ADD100},
showing again that the left-hand side of \eqref{sec3-ADD100} is finite.

Now it follows from part (ii) that $\abs{Z_n}^\alpha\1[\{ \tau^\prime \ge n\}] \to \abs{Z}^\alpha \1[\{ \tau^\prime=\infty\}]$ ${\mathbb P}^\alpha$-a.s.\ as $n \to \infty$. Consequently,
$$ \sup_{v \in F} \E_{\bfDV}^\alpha \Big[ \abs{Z}^\alpha \1[\{ \tau^\prime = \infty\}] \Big] ~=~ \sup_{v \in F} \E_{\bfDV}^\alpha \left[ \lim_{n \to \infty} \abs{Z_n}^\alpha \1[\{ \tau^\prime\ge n\}] \right] ~\le~ \sup_{v \in F} \E_{\bfDV}^\alpha
  \left[  \sup_{n \in \N} \abs{Z_n}^\alpha \1[\{ \tau^\prime \ge n\}] \right]  < \infty.$$

(iv) The almost sure convergence $\abs{Z_n}^\alpha\1[\{ \tau^\prime \ge n\}] \to \abs{Z}^\alpha \1[\{ \tau^\prime=\infty\}]$ was obtained in part (ii), and it was shown in part (iii) that the sequence $\big\{\abs{Z_n}^\alpha\1[\{ \tau^\prime \ge n\}]\big\}_{n \in \N}$ is uniformly integrable, and the $L^1$-convergence follows.
\halmos

\subsection{Markov nonlinear renewal theory} \label{sect:nonlinearMRT}
We conclude this section by applying Markov nonlinear renewal theory, as developed by Melfi \citeyear{Melfi1992,Melfi1994}, to the processes $\{V_n\}$ and 
$\{ V_n^A \}$, where
\begin{equation} \label{def:VnA}
V_n^A := \frac{V_n}{d_A(\widetilde{V}_n)},  \qquad n \in \N,
\end{equation}
and where we assume here that $d_A$ is bounded and continuous. Melfi's theory allows us to compare the overjump distributions for these 
two processes to those of $\{(X_n, S_n) \}$ and $\{(X_n, S_n^A) \}$, respectively,
where $S_n^A = S_n - \log d_A(X_n)$.

We begin by verifying the conditions of Kesten's renewal theorem for the Markov random walks $\{(X_n, S_n)\}$ and $\{(X_n, S_n^A)\}$ under the $\alpha$-shifted measure.
We start by quoting these conditions as they are stated in \shortciteN{Kesten1974}, with notation adapted for the process $\{(X_n, S_n)\}$.
\begin{itemize}
\item[{\bf I.1}] There exists a measure $\eta_\alpha$ on $\SPdp$ which is invariant for $\{ X_n\}$, and  $\Prob^\alpha_x(X_n \in E, \text{ for some }n )=1$ for all  $x \in \SPdp$ and all {\em open} sets $E \subset \SPdp$ with $\eta_\alpha(E)>0$.
\item[{\bf I.2}] $\E^\alpha_{\eta_\alpha} [S_1-S_0]$ exists and 
$ \lim_{n \to \infty} n^{-1} S_n = \E^\alpha_{\eta_\alpha} [S_1-S_0]>0$ $\Prob^\alpha$-a.s.
\item[{\bf I.3}] There exists a sequence $\{\zeta_i\} \subset \R$ such that the group generated by $\{\zeta_i\}$ is dense in $\R$, and such that for each $\zeta_i$ and each $\delta >0$, there exists $y=y(i,\delta) \in \SPdp$ with the following property: For each $\epsilon >0$, there exists $E \subset \SPdp$ with $\eta_\alpha(E)>0$, $m_1, m_2 \in \pintegers$, and $\gamma \in \R$ such that for any $x \in E$,
\begin{align}
&\Prob_x^\alpha\big( |X_{m_1} - y|<\epsilon, |S_{m_1} - \gamma| \le \delta \big) ~>~0 \quad \text{ and } \\[.2cm]
&\Prob_x^\alpha\big( |X_{m_2} - y|<\epsilon, |S_{m_2} - \gamma - \zeta_i| \le \delta \big)  ~>~0.
\end{align}
\item[{\bf I.4}] For each fixed $x_0 \in \SPdp$ and $\epsilon >0$, there exists $r=r(x_0,\epsilon)$ such that for 
all bounded measurable functions $f: (\SPdp \times \R)^\N \to \R$ and for all $y \in B_r(x_0)$,
\begin{align}
\E^\alpha_y \big[ f(X_0,S_0, X_1, S_1,\dots) \big] ~ \le&~ \E^\alpha_x \big[f^\epsilon(X_0, S_0,X_1, S_1, \dots) \big] + \epsilon |f|_\infty \quad \text{ and } \\[.2cm]
\E^\alpha_x \big[ f(X_0,S_0, X_1, S_1,\dots) \big] ~ \le&~ \E^\alpha_y \big[f^\epsilon(X_0, S_0,X_1, S_1, \dots) \big] + \epsilon |f|_\infty,
\end{align}
where
$ f^\epsilon(x_0,s_0,x_1,s_1,\dots) := \sup \big\{ f(y_0,t_0, y_1, t_1, \dots) \, : \, |x_i-y_i| + |s_i-t_i| < \epsilon \text{ for all } i \in \N \big\}.$
\end{itemize}

Let us give a brief interpretation of these conditions.
Condition I.1 is weaker than $\eta_\alpha$-irreducibility, for only $\eta_\alpha$-positive {\em open} sets are required to be reachable from any initial state. Condition I.2 is the classical assumption of positive drift. Condition I.3 is the implementation of the non-arithmeticity condition;
while Condition I.4 states that if $x \to y$, then $\Prob^\alpha_x\big( (X_n, S_n)_{n \in \N} \in \cdot \big)$ converges to $\Prob^\alpha_y\big( (X_n, S_n)_{n \in \N} \in \cdot \big)$ in Prokhorov distance; for details, see \shortciteN{Melfi1994}, Section 2.2.

\begin{lemma}\label{lem:kesten conditions for XnSnA}
Assume $(H_1)$ and $(H_2)$. Then Conditions I.1 --  I.4 are satisfied by $\{(X_n, S_n):  n \in \N\}$. If $d_A$ is bounded and continuous on $\SPdp$, then these conditions are also satisfied by $\{(X_n, S_n^A): n
\in \N \}$.
\end{lemma}

\noindent {\bf Proof}.  Hypotheses $(H_1)$ and $(H_2)$ imply the assumptions of Proposition 1 in \shortciteN{HK73}, where the validity of {I.1} -- {I.4} is proved for the Markov random walk $\{(X_n, S_n)\}$.  

We now verify {I.1} -- {I.4} for the process $\{(X_n, S_n^A)\}$. 
 {I.1} remains valid, as it only concerns $\{X_n\}$.  Next observe that $S_n^A=S_n - \log d_A(X_n)$ is a measurable function of $(X_n,S_n)$.  Hence any measurable function $f(X_0,S_0^A, X_1, S_1^A, \dots)$ can be transformed into a measurable function of $\{(X_n, S_n)\}$, and thus {I.4} holds as well [using that $d_A$ is bounded from above and below]. Since I.3 holds for $\{(X_n, S_n)\}$, it holds for $\{(X_n, S_n^A)\}$ as well by replacing $\gamma$ with $\gamma^\prime:=\gamma - \log d_A(y)$. Finally, I.2 follows by applying Lemma 2.4, which gives
\begin{align*}
\hspace*{1.6cm}
 \lim_{n \to \infty} \frac{S_n^A}{n} &=  \frac1n \sum_{i=1}^n \big(S_i - S_{i-1}\big) - \big(\log d_A(X_i) - \log d_A(X_{i-1})\big) + \frac1n{\log d_A(X_0)} \\
 &= \E_{\eta_\alpha}^\alpha\big[ S_1 - S_0 - (\log d_A(X_1)-\log d_A(X_0))\big] = \E_{\eta_\alpha}^\alpha[S_1^A- S_0^A] \quad \text{ $\Prob^\alpha$-a.s.}  \hspace*{1.6cm} \halmos
 \end{align*}
 
Write ${\mathfrak T}_u^A=\inf \{n \in \pintegers \, : \, S_n^A > \log u\}$ and ${\mathfrak T}_u=\inf\{n \in \pintegers \, : \, S_n > \log u\}.$ Then Kesten's renewal theorem (\shortciteN{Kesten1974},
Theorem 1) yields the joint asymptotics of the overjump above level $\log\, u$ and the position of $X_n$ at the time of the overjump.

\begin{thm}[Kesten, 1974]\label{prop:asymptotic overjump distribution}
Assume $(H_1)$ and $(H_2)$.   Then:
\begin{enumerate}
\item[{\rm (i)}]  There is a probability measure $\rho$ on $\SPdp \times (0,\infty)$, such that, for all $x \in \SPdp$ and all functions $f \in {\cal C}_b\Big(\SPdp \times (0,\infty)\Big)$, 
$$ \lim_{u \to \infty} \E^\alpha_x\Big[ f(X_{{\mathfrak T}_u}, S_{{\mathfrak T}_u}-\log u) \Big] ~=~ \int_{\SPdp \times \reals_+} f(y,s) \, \rho(dy,ds).$$
\item[{\rm (ii)}]  If the function $d_A$ is bounded and continuous on $\SPdp$, then there is a probability measure
 $\rho^A$ on $\SPdp \times (0,\infty)$, such that for all $x \in \SPdp$ and all $f \in {\cal C}_b\Big(\SPdp \times (0,\infty)\Big)$, 
$$ \lim_{u \to \infty} \E^\alpha_x\Big[ f(X_{{\mathfrak T}_u^A}, S^A_{{\mathfrak T}_u^A}-\log u) \Big] ~=~ \int_{\SPdp \times \reals_+} f(y,s) \, \rho^A(dy,ds).$$
\end{enumerate}
\end{thm}

A representation of $\rho$ in terms of the ascending ladder height process is given in Eq. \eqref{eq:markov delay distribution} below.

\medskip

The essential result we will need is the nonlinear Markov renewal theorem by Melfi (1992, 1994), which extends Kesten's renewal theorem to a wider class of processes which are asymptotically ``close" to the Markov random walk $\{(X_n, S_n):  n \in \N\}$. 

Let $\{(Y_n, W_n):  n \in \N \}$ be a stochastic process on $\SPdp \times \R$, adapted to a filtration $\{ {\cal G}_n\}$, and 
throughout the remainder of this section (with a slight abuse of notation) set
\begin{align*}
T_u ~&:=~ \inf\{n \in \N \, : \, W_n > \log u\}, \\
{\mathfrak T}_u ~&:=~ \inf\{n \in \N \, : \, S_n > \log u\}, \\
(Y_{u,k}, W_{u,k}) ~&:=~ (Y_{T_u + k}, W_{T_u+k} - W_{T_u}) \quad \text{ for } k \ge 0.
\end{align*}
Then we have the following (\shortciteN{Melfi1994}, Theorem 3).

\begin{thm}[Melfi, 1994]  \label{thm:Melfi nonlinear renewal}
Let $\{(X_n, S_n):  n \in \N\}$ be a Markov random walk satisfying the assumptions of Kesten's renewal theorem under $\Prob^\alpha$. Assume the following conditions hold:
\begin{itemize}
\item[{\rm (I)}] For all $m\ge 1$, the Prokhorov distance 
\begin{equation}\label{eq:Melfi condition I} {\mathfrak d} \bigg( \Prob^\alpha \Big( \left. (Y_{u,k}, W_{u,k})_{1 \le k \le m} \in \cdot \, \right| {\cal G}_{T_u} \Big) ~,~ \Prob^\alpha_{Y_{T_u}} \Big(  (X_k, S_k)_{1 \le k \le m} \in \cdot \Big) \bigg)
\end{equation}
 converges to 0 in  $\Prob^\alpha$-probability.
 \item[{\rm (II)}] $\{W_{T_u} - \log u\}_{u \ge 1}$ is tight under  $\Prob^\alpha$.
  \item[{\rm (III)}] $\{Y_{T_u}\}_{u \ge 1}$ is tight  under $\Prob^\alpha$.
\end{itemize}
Let $\rho$ denote the asymptotic overjump distribution of $\{(X_n, S_n): n \in \N\}$ obtained in 
Theorem \ref{prop:asymptotic overjump distribution}.
Then for all $f \in {\cal C}_b \Big( \SPdp \times (0,\infty)\Big)$ and all $x \in \SPdp$,
\begin{equation}\label{eq:Melfi result} \lim_{u \to \infty} \E^\alpha_x\Big[ f(Y_{T_u}, W_{T_u}-\log u) \Big] ~=~ \int_{\SPdp \times \reals_+} f(y,s) \, \rho(dy,ds). 
\end{equation}
\end{thm}

Next, we verify the assumptions of Melfi's theorem for the process $\{(Y_n, W_n)\}$ chosen specifically as $\{ (\widetilde{V}_n, \log \abs{V_n})\}$ and $\big\{\big(\widetilde{V}_n^A, \log \big| V_n^A \big| \big) \big\}$, respectively. 
Note that in this case ${\cal G}_n = {\cal F}_n$, and condition (III) is always satisfied since $Y_{T_u} \in \SPdp$, which is compact.

\begin{lemma} \label{lm-Prokhorov} \label{lem:melfi1}
Assume $(H_1)$ and $(H_2)$. Then Condition (I) is satisfied by the pair
$\{ (\widetilde{V}_n, \log |V_n| ) \}$, $\{ ( X_n, S_n ) \}$.
\end{lemma}

\noindent {\bf Proof.}
Using the Markov property,
\begin{align*}
\Prob^\alpha & \left( \left. (Y_{u,k}, W_{u,k})_{1 \le k \le m} \in \cdot \, \right| \F_{T_u} \right) \\
=&~ \Prob^\alpha \left( \left. (\widetilde{V}_{T_u + k}, \log \abs{V_{T_u + k}} - \log \abs{V_{T_u}})_{1 \le k \le m} \in \cdot \, \right| 
X_{T_u},  V_{T_u} \right) \\
=&~ \Prob^\alpha_{X_{T_u}, V_{T_u}} \left( (\widetilde{V}_{k}, \log \abs{V_{k}} - \log \abs{V_0})_{1 \le k \le m} \in \cdot  \right) .
\end{align*}
By \shortciteN{Guivarch2012}, Lemma 3.5,  the total variation distance between $\Prob^\alpha_{x,v}$ and $\Prob^{\alpha}_{y,v}$ is bounded above by $B \abs{x-y}^{\bar{\alpha}}$ for some finite constant $B$, where $\bar{\alpha}=\min\{\alpha,1\}$. [Their proof is in the setting of invertible matrices, but carries over to the present setting of nonnegative matrices without change.]  Convergence in total variation implies convergence in the Prokhorov metric ${\mathfrak d}$, and thus
\begin{align} & {\mathfrak d}  \Bigl( \Prob^\alpha_{X_{T_u}, V_{T_u}} \left( (\widetilde{V}_{k}, \log \abs{V_{k}} - \log \abs{V_0})_{1 \le k \le m} \in \cdot  \right) \, , \, \Prob^\alpha_{\widetilde{V}_{T_u}} \Big(  (X_k, S_k)_{1 \le k \le m} \in \cdot \Big) \Bigr) \nonumber \\
\le& B \abs{X_{T_u}-\widetilde{V}_{T_u}}^{\bar{\alpha}} \! \!+ \! {\mathfrak d} \Bigl( \Prob^\alpha_{\widetilde{V}_{T_u}, V_{T_u}} \left( (\widetilde{V}_{k}, \log \abs{V_{k}} - \log \abs{V_0})_{1 \le k \le m} \in \cdot  \right)  , \Prob^\alpha_{\widetilde{V}_{T_u}} \left(  (X_k, S_k)_{1 \le k \le m} \in \cdot \right)\Bigr).  \label{eq:nonlinear1} 
\end{align}

We begin by considering the second term. Fix the initial values $(\widetilde{V}_{T_u}, V_{T_u})=(\wt{v},v)$, and introduce the notation $V_k^{(0)} := \sum_{j=1}^k M_k \cdots M_{j+1} Q_j$, $k =2,3,\ldots$ 
and $V_1^{(0)} = Q_1$.  Then  for any $k \in \pintegers$,
\begin{equation} \label{defV_0}
\Big|\log \abs{V_{k}} - \log \abs{V_0} - S_k\Big| ~=~ \log \abs{ \frac{\Pi_k v + V_k^{(0)}}{\abs{ \Pi_k \widetilde{v}}\abs{v}} }   
~\le~ \log  \left( 1 + \frac{|V_k^{(0)}|}{\abs{\Pi_k v}}  \right) ~\le \frac{\big|V_k^{(0)} \big| }{\abs{\Pi_k v}}
\end{equation}
and
\begin{equation}
\abs{\widetilde{V_k} - X_k} ~=~ \abs{\frac{V_k}{\abs{V_k}} - \frac{V_k}{\abs{\Pi_k v}} + \frac{V_k}{\abs{\Pi_k v}} - \frac{\Pi_k v}{\abs{\Pi_k v}} } 
\le~ \abs{V_k}  \abs{\frac{\abs{\Pi_kv}-\abs{V_k}}{\abs{V_k} \abs{\Pi_kv}}} ~+~\frac{\big|V_k^{(0)}\big|}{\abs{\Pi_k v}} 
\le~  \frac{2\big|V_k^{(0)}\big|}{\abs{\Pi_k v}}, \label{estimate Vtilde X}
\end{equation}
using the triangle inequality.
Thus for all $v \in \Rdo$,
\begin{align} \Prob^\alpha_{\tilde{v}} & \Big(  \abs{(\widetilde{V}_{k}, \log \abs{V_{k}}  - \log \abs{V_0})_{1 \le k \le m} - (X_k, S_k)_{1 \le k \le m} }_\infty \ge \epsilon \, \big| \, V_0=v\Big) \nonumber\\[.2cm]
 & \hspace*{2.5cm} \le~ \Prob^\alpha_{\tilde{v}} \left( 2 \sum_{k=1}^m  \frac{|{V_k^{(0)}}|}{\abs{\Pi_k v}} \ge \epsilon \right)
~\le~ \sum_{k=1}^m \Prob^\alpha_{\tilde{v}} \left( \frac{|V_k^{(0)}|}{\abs{\Pi_k v}} \ge \frac{\epsilon}{2m} \right)
\nonumber\\[.2cm]
 & \hspace*{2.5cm} \le~ \frac{(2m)^\alpha}{\epsilon^\alpha} \sum_{k=1}^m \E^\alpha_{\tilde{v}} \left[ \left(
  \frac{1}{|v|} \frac{|V_k^{(0)}|}{\abs{\Pi_k \tilde{v}}} \right)^\alpha \, \right] ~\le~ \frac{B}{\epsilon^\alpha |v|^\alpha} \sum_{k=1}^m \E \Big[ \big|V_k^{(0)} \big|^\alpha \Big], \label{estimate V0Pin}
\end{align}
for some universal constant $B$, where we used Chebyshev's inequality, the definition of the $\alpha$-shifted measure, and the boundedness of $r_\alpha$  in the last identity. 
Hence
$$ \lim_{u \to \infty} \sup_{v \, : \, |v| \ge u} \, \Prob^\alpha_{\tilde{v}} \Big( \abs{(\widetilde{V}_{k}, \log \abs{V_{k}} - \log \abs{V_0})_{1 \le k \le m} - (X_k, S_k)_{1 \le k \le m} }_\infty \ge \epsilon \, \Big| \, V_0=v\Big) ~=~0.$$
Recall that convergence in probability
implies weak convergence, which is equivalent to convergence in the Prokhorov metric.
Since $\Prob^\alpha(T_u < \infty)=1$ and $|V_{T_u}|\ge u$, we conclude that
$$ \lim_{u \to \infty} {\mathfrak d} \Bigl( \Prob^\alpha_{\widetilde{V}_{T_u}, V_{T_u}} \left( (\widetilde{V}_{k}, \log \abs{V_{k}} - \log \abs{V_0})_{1 \le k \le m} \in \cdot  \right)  \, , \, \Prob^\alpha_{\widetilde{V}_{T_u}} \Big(  (X_k, S_k)_{1 \le k \le m} \Big)\Bigr)  ~=~0 \quad \Prob^\alpha\text{-a.s.}$$

To complete the proof, we observe that the first term on the right-hand side 
of \eqref{eq:nonlinear1} tends to zero in $\Prob^\alpha$-probability by Lemma \ref{lem:VTutoXTu},
which is given next. \halmos

\begin{lemma}\label{lem:VTutoXTu}
As $u \to \infty$, we have for any $x \in \SPdp$ and $V_0=v \in \realsd_+ \setminus \{ 0 \}$  that
\begin{equation} \label{eq:VTutoXTu}\abs{X_{T_u}- \widetilde{V}_{T_u}}
\to 0 \quad\text{\rm in $\Prob^\alpha$-probability}.\end{equation}
\end{lemma}

\noindent {\bf Proof.}
Let $w = u/2$, and decompose the process based on its behavior prior and subsequent to the time $T_w$.
This yields
\begin{align*}
 \Prob_x^\alpha\left(\abs{\widetilde{V}_{T_u}-X_{T_u}}> \epsilon\right) 
&~=~ \Prob_x^\alpha\left(\abs{\frac{\Pi_{T_{w}+1}^{T_u}V_{T_{w}} + \sum_{k=T_{w+1}}^{T_u} \Pi_{k+1}^{T_u} Q_k}{\abs{V_{T_u}}} - 
  \Big(\Pi_{T_{w}+1}^{T_u}X_{T_{w}}\Big)^\sim} > \epsilon \right) \\[.2cm]
&~\le~ \E_x^\alpha \left[ \Prob_x^\alpha\left(\left. \abs{\frac{\Pi_{T_{w}+1}^{T_u}V_{T_{w}} + \sum_{k=T_{w}+1}^{T_u} \Pi_{k+1}^{T_u} Q_k}{\abs{V_{T_u}}} - \frac{\Pi_{T_{w}+1}^{T_u}V_{T_{w}}}{\abs{\Pi_{T_{w}+1}^{T_u}V_{T_{w}}}}} > \frac{\epsilon}{2} \, \right| \F_{T_{w}} \,\right) \right]\\[.2cm]
& \hspace*{2cm} +~ \Prob_x^\alpha\left(\abs{ \Big(\Pi_{T_{w}+1}^{T_u}V_{T_{w}}\Big)^\sim -  \Big(\Pi_{T_{w}+1}^{T_u}X_{T_{w}}\Big)^\sim} > \frac{\epsilon}{2} \right) ~=:~ {\mathbb I}_1(u) + {\mathbb I}_2(u)
\end{align*}
by the triangle inequality.

To compute ${\mathbb I}_2(u)$ as $u \to \infty$,  we apply Proposition 3.1 of \shortciteN{Hennion1997}, which gives that
\begin{equation*} 
\sup_{x,y \in \SPdp} \abs{(\Pi^n_{i+1} x)^\sim - (\Pi^n_{i+1} y)^\sim}
    \le 2 {\bf c} (\Pi^n_{i+1})
\end{equation*}
for a  function ${\bf c}(\cdot)$ is bounded above by one and tends to zero $\Pstat$-a.s.\ as $n-i$ tends to infinity (\shortciteN{Hennion1997}, Lemma 3.2). Since $\Prob^\alpha_x$ is absolutely continuous with respect to the measure $\Pstat$ (\shortciteN{Buraczewski.etal:2014}, Lemma 6.2), it follows that ${\bf c}\big(\Pi^{T_u}_{T_{u/2}+1}\big) \to 0$ $\Prob^\alpha_x$-a.s. for all $x \in \SPdp$, and hence
\begin{equation}
{\mathbb I}_2(u)  \le \Prob_x^\alpha\left( 2 {\bf c} \left(\Pi_{T_{w}+1}^{T_u} \right)  > \frac{\epsilon}{2} \right)
   \searrow 0 \quad \mbox{\rm as} \quad u \to \infty.
\end{equation}

Now consider ${\mathbb I}_1(u)$ as $u \to \infty$.  Repeating the calculation in \eqref{estimate Vtilde X} yields
\begin{equation} \label{aug9}
{\mathbb I}_1(u) \le  \E_x^\alpha \left[ \Prob_{\widetilde{V}_{T_{w}}}^\alpha\left( \frac{4\big| V_{T_u}^{(0)} \big| }{u \big|\Pi_{T_u} \widetilde{V}_0 \big| } > \frac{\epsilon}{2} \right) \right].
\end{equation}
Now recall from Lemma \ref{lem:convergence Z} (i) that
\begin{equation} \label{aug9-1}
{\mathcal Z}^0:= \sup_{n \in \N} \frac{ |V_n^{(0)}|}{ |M_n \cdots M_1 X_0| } < \infty \quad \Pstat\text{-a.s.}.
\end{equation}
Using that $\Prob^\alpha_{\tilde{V}_{T_w}} \le B \Pstat$ for some universal constant $B$ by Lemma 6.2 of \shortciteN{Buraczewski.etal:2014}, we obtain
\begin{equation} \label{aug9-2}
{\mathbb I}_1(u)~\le~ \lim_{u \to \infty}  \Pstat\left(  \frac{{\mathcal Z}^0}{u}  > \frac{\epsilon}{8} \right) = 0.
\end{equation}
\halmos

It remains to check condition (II) of Melfi's theorem.

\begin{lemma}\label{lem:tightness}
$\{ W_{T_u} -\log u\}_{u \ge 1}$ is tight under $\Prob^\alpha$.
\end{lemma}

\noindent
{\bf Proof.}
A sufficient condition is given in \shortciteN{Melfi1994}, Section 5.2.  Letting
 $\xi_n := \log|V_{n}| - S_n$
 and supposing that
 $\{\xi_{T_u}\}_{u \ge 1}$ and $\{\xi_{{\mathfrak T}_u}\}_{u \ge 1}$ are tight under $\Prob^\alpha$, then 
 it follows that $\{W_{T_u}-\log u\}_{u \ge 1}$ is tight.  Now by Lemma \ref{lem:convergence Z},
$$ \xi_n =  \log \frac{\abs{V_n}}{\abs{\Pi_n V_0}}  \to \log Z \quad \Prob^\alpha\text{-a.s.},$$
for a finite random variable $Z$.
Since $T_u$ and ${\mathfrak T}_u$ are stopping times with respect to the filtration $\{ \F_n \}$ and tend to infinity as $u \to \infty$, we deduce that
$$ \lim_{u \to \infty} \xi_{T_u} = \log Z \quad \Prob^\alpha\text{-a.s.}\qquad  \mbox{and} \qquad  \lim_{u \to \infty} \xi_{{\mathfrak T}_u} =\log Z \quad \Prob^\alpha\text{-a.s.}$$
Thus, in particular, the families $\{\xi_{T_u}\}_{u \ge 1}$ and $\{\xi_{{\mathfrak T}_u}\}_{u \ge 1}$ converge in distribution under $\Prob^\alpha$ and are consequently tight.
\halmos

 We are now in a position to apply Melfi's theorem.

\begin{thm}\label{thm:Melfi} \label{nonlinear renewal TuA}
Assume $(H_1)$ and $(H_2)$. Then for all $f \in {\cal C}_b\Big( \SPdp \times (0,\infty)\Big)$ and all $x \in \SPdp$ and $v \in \realsd_+ \setminus \{ 0 \}$,
\begin{equation}\label{eq:Melfi result rho} \lim_{u \to \infty} \E^\alpha_{x,v}\left[ f\left(\widetilde{V}_{T_u}, \log \frac{|V_{T_u}|}{u}\right) \right] ~=~ \int_{\SPdp \times \reals_+} f(y,s) \, \rho(dy,ds). 
\end{equation}
Further, if $d_A$ is bounded and continuous on $\SPdp$, then for all $f \in {\cal C}_b\Big( \SPdp \times (0,\infty)\Big)$ and all $x \in \SPdp$ and $v \in \realsd_+ \setminus \{ 0 \}$,
\begin{equation}\label{eq:Melfi result rhoA} \lim_{u \to \infty} \E^\alpha_{x,v}\left[ f\left(\widetilde{V}_{T_u^A}^A , \log \frac{|V_{T_u^A}^A|}{u}\right) \right] ~=~ \int_{\SPdp \times \reals_+} f(y,s) \, \rho^A(dy,ds). 
\end{equation}
\end{thm}

\noindent {\bf Proof.}
Condition (III) is necessarily satisfied since $\SPdp$ is compact, and for the process $\{(\widetilde{V}_n, \log |V_n| )\}$,
the validity of Conditions (I) and (II) has been proved in Lemmas \ref{lem:melfi1} and \ref{lem:tightness}, respectively. Thus \eqref{eq:Melfi result rho} follows from Theorem \ref{thm:Melfi nonlinear renewal}.

Turning to \eqref{eq:Melfi result rhoA}, we need to check the validity of Conditions (I) and (II) for $\{(\widetilde{V}_n^A, \log |V_n^A| )\}$. Recall that $V_n^A:=V_n/d_A(\widetilde{V}_n)$, implying that $\widetilde{V}_n^A=\widetilde{V}_n$ (since these
two quantities have the same direction).
Moreover $S_n^A=S_n - \log d_A(X_n)$.  Thus, for $f(x,s)=(x, s-\log d_A(x))$, we
have that $\{(\widetilde{V}_n^A, \log |V_n^A|)\} = \{f(\widetilde{V}_n, \log |V_n|)\}$ and $\{(X_n, S_n^A)\}=\{f(X_n, S_n)\}$. Hence,
for the process $\{ (\widetilde{V}_n^A, \log |V_n^A| ) \}$,
Condition (I) can be deduced from Lemma \ref{lem:melfi1}. 
Finally, since $d_A$ is bounded, the tightness of $\{\log |V_{T_u^A}^A| - \log u\}=\{\log|V_{T_u^A}| - \log d_A(V_{T_u^A}) - \log u\}$ follows, in the same way, from  Lemma \ref{lem:tightness}. Thus we conclude \eqref{eq:Melfi result rhoA}.
\halmos

We conclude with a result concerning the first passage times in the $\alpha$-shifted measure.

\begin{lemma} \label{lemma-charac-T}
Assume that $d_A$ is bounded and continuous. Then
\begin{equation} \label{charac-T}
\lim_{u \to \infty} \frac{T_u^A}{\log u} = \frac{1}{\lambda'(\alpha)} \quad \mbox{\rm in $\Prob^\alpha$-probability}.
\end{equation}
\end{lemma}

\noindent
{\bf Proof.}
By definition, $V_n = Z_n e^{S_n}$ and $V_n^A = V_n/d_A(\widetilde{V}_n)$, and consequently
\begin{equation} \label{aug9-3}
\log|V_n^A| = S_n + |Z_n| - \log d_A(\wt{V}_n) := S_n + \xi_n.
\end{equation}
Also, it follows by definition that 
$T_u^A := \inf\{ n:  V_n \in A \} = \inf\{n:  |V_n| > d_A(\widetilde{V}_n) u \}
  = \inf\{ n:   \log|V_n^A| > \log u\}.$
Now recall that $\sup_{n \in \N} Z_n$ is finite a.s., by Lemma \ref{lem:convergence Z}. Since $d_A$ is bounded, it
follows that the sequence $\{ \xi_i \}$ in \eqref{aug9-3} satisfies 
$$ \lim_{n \to \infty} \frac{1}{n} \left( \max_{1 \le k \le n} \xi_k \right) = 0 \quad \Prob^\alpha\text{-a.s.},$$
i.e.\ $\{\xi_n\}$ is {\em slowly changing} (as defined in \shortciteN{DS85}, Eq.\ (9.5)).
Moreover, by  Lemma \ref{lem:slln}, $S_n/n \to \lambda'(\alpha)$ a.s., and hence $\log|V_n^A|/n \to \lambda'(\alpha)$ a.s. The
result then follows by reasoning as in \shortciteN{DS85}, Lemma 9.13. \halmos

%
%
%
%

\setcounter{equation}{0}
\section{Characterizing the large exceedances over cycles} \label{sect:4}
Before turning to the proofs of the main theorems of the paper, we first establish a few required results
concerning the behavior of the post-$T_u^A$ process. 
The central results of this section are Proposition \ref{prop1}---which will be used throughout the paper---and Proposition \ref{prop:Kesten result}, which will be the basis for the proof of Theorem \ref{thm2A} in 
the next section.

Recall that $\tau$ denotes the return time to a set $\DD=B_r(0) \cap \Rd_+\setminus\{0\}$ for $\pi(\DD)>0$, and that $T_u^A := \inf\left\{ n :  V_n \in u A \right\} = \inf\left\{ n:  |V_n^A| > u \right\}$,
where
\begin{equation} \label{def-V*}
V_n^A := \frac{V_n}{d_A(\widetilde{V}_n)}, \quad n = 0,1,\ldots.
\end{equation}
Also recall that
\begin{equation} 
r_\alpha^A (x) = r_\alpha (x) \left(d_A(x) \right)^\alpha, \quad x \in {\mathbb S}^{d-1}_+.
\end{equation}
Finally, we say that a function $g: (\Rd_+)^{m+1} \to \R$ is {\em almost $\theta$-H\"older continuous} if 
\begin{equation}\label{def:almost hoelder} g(v_0, \dots, v_m) = \hat{g}(v_0, \dots, v_m) \1[{\{|v_m| \ge \delta\}}] \end{equation}
for some $\delta \ge 0$ and  $\theta$-H\"older continuous  $\hat{g}$.

Much of this section will be devoted to the proof of the following proposition, which can be viewed as a generalization of \shortciteN{JCAV13}, Proposition 6.1, to the setting
of matrix recursions.

\begin{prop} \label{prop1}
Assume Hypotheses $(H_1)$ and $(H_2)$ are satisfied.  Let $m \in \N$ and $g: (\Rd_+)^{m+1} \to \R$ be a bounded almost $\theta$-H\"older continuous function for $\theta \le \min\{1, \alpha\}$, and assume that the
function $d_A$ is bounded and continuous on $\SPdp$.  Then for any $v \in \R^d_+\setminus\{0\}$, 
\begin{align} \label{lm4.1a}
&\lim_{u \to \infty}  u^\alpha {\mathbb E} \left[  \left. g \left( \frac{V_{T_u^A}}{u}, \ldots, \frac{V_{T_u^A +m}}{u} \right) {\bf 1}_{\{ T_u^A < \tau \}} \,\right|
  \, V_0 = v  \right]   \nonumber  \\[.1cm]
&  =  r_\alpha (\widetilde{v})  {\mathbb E}_{\bfDV}^{\alpha} \left[ |Z|^\alpha {\bf 1}_{\{ \tau = \infty \}} \right] \int
    \frac{e^{-\alpha s}}{r_\alpha^A(x)}  {\mathbb E} \Big[ g\big( e^{S_0} X_0, \ldots, e^{S_m} X_m \big) \, \big| \, X_0=x, \, S_0=s+ \log d_A(x)  \Big] \, \rho^A(dx,ds).
\end{align}
\end{prop}

Observe that $\rho^A$ is the asymptotic law, as $u \to \infty$, of $\big(({V_{T_u^A}^A})^\sim
, \log |V_{T_u^A}^A| - \log u\big)$, while on the left-hand side of \eqref{lm4.1a},
we evaluate the function $g$ for the process $\{V_n\}$ (not $\{ V_n^A\}$) 
at a sequence of times commencing  at time $T_u^A$. Taking into account \eqref{def-V*}, this explains the additional summand $\log d_A(x)$ in the expression for $S_0$; namely, it arises when transforming $V_{T_u^A}^A$ into $V_{T_u^A}$.

An important special case occurs when we take $m = 0$ and $g(V_{T_u^A}/u) \equiv 1$, in which case we obtain:

\begin{cor} \label{cor4-1}
Assume Hypotheses $(H_1)$ and $(H_2)$ are satisfied and the function $d_A$ is bounded and continuous on $\SPdp$.
Then for any $v \in \R^d_+\setminus\{0\}$,
\begin{align} 
\lim_{u \to \infty} & u^\alpha {\mathbb P} \left(  \left. T_u^A < \tau\, \right|
\, V_0 = v  \right) =  r_\alpha (\widetilde{v})  {\mathbb E}_{\bfDV}^{\alpha} \left[ |Z|^\alpha {\bf 1}_{\{ \tau = \infty \}} \right]
    \int_{{\mathbb S}^{d-1}_+ \times \reals_+}   \frac{e^{-\alpha s}}{r_\alpha^A(x)}  \, \rho^A(dx,ds). \label{lm4.1b}
\end{align}

\end{cor}

\noindent If $A = \{ v \in \R^d_+ \, : \,  |v| > 1\}$, then $T_u^A$ reduces to the first exceedance time of $|V_n|$ above the level $u$;
that is, $T_u^A = T_u := \inf\{n:  |V_n| > u \}$.
Furthermore, we then have 
$d_A(x) = 1$ and thus  $r^A_\alpha = r_\alpha$.  Consequently, in this case, Proposition \ref{prop1} and Corollary \ref{cor4-1}
hold with $T_u^A$, $r_\alpha^A$, $\rho^A$ replaced with $T_u$, $r_\alpha$, $\rho$, respectively.
Then we can easily apply the definition of the $\alpha$-shifted measure to obtain:

\begin{cor} \label{cor:prop1 for Tu}
Assume Hypotheses $(H_1)$ and $(H_2)$ are satisfied.  Let $m \in \N$ and $g: (\Rd_+)^{m+1} \to \R$ be a bounded almost $\theta$-H\"older continuous function for $\theta \le \min\{1, \alpha\}$.  Then for any $v \in \R^d_+\setminus\{0\}$, 
\begin{align} 
\lim_{u \to \infty} & u^\alpha {\mathbb E} \left[  \left. g \left( \frac{V_{T_u}}{u}, \ldots, \frac{V_{T_u +m}}{u} \right) {\bf 1}_{\{ T_u < \tau \}} \,\right|
  \, V_0 = v  \right] \nonumber \\[.1cm]
&  =  r_\alpha (\widetilde{v})  {\mathbb E}_{\bfDV}^{\alpha} \left[ |Z|^\alpha {\bf 1}_{\{ \tau = \infty \}} \right]
    \int_{{\mathbb S}^{d-1}_+ \times \reals_+} 
      {\mathbb E}^\alpha_x \bigg[ \frac{e^{-\alpha (S_m+s)}}{r_\alpha(X_m)} g\big( e^{s} X_0, \ldots, e^{S_m+s} X_m \big)  \bigg] \, \rho(dx,ds). \label{lm4.1a2}
\end{align}
\end{cor}

To establish Proposition \ref{prop1}, we will rely on the following.

\begin{lemma} \label{sublemma}
Assume the conditions of Proposition \ref{prop1}.  Then:
\begin{enumerate}[{\rm (i)}]
\item For all $v \in \Rd_+\setminus\{0\}$, we have the $L^1$-convergence
\begin{equation} \label{slm4.1.1}
\lim_{n \to \infty} \lim_{u \to \infty} \E_{\bfDV}^\alpha \bigg[\Big| \big|Z_{T_u^A}\big|^\alpha {\bf 1}_{\{ T_u^A < \tau \}} - \abs{Z_n}^\alpha  {\bf 1}_{\{n \le T_u^A\}} {\bf 1}_{\{ n \le \tau\}}  \Big| \bigg] = 0. 
\end{equation}
\item  Let \label{sublemma:part3}
\begin{equation} \label{def-cG}
{\mathfrak G}_u ~:=~ \frac{1}{r_\alpha(X_{T_u^A})} \left( \frac{\big| V_{T_u^A}\big|}{u}\right)^{-\alpha} \E \bigg[ g \left( \frac{V_{T_u^A}}{u}, \ldots, \frac{V_{T_u^A +m}}{u} \right) \bigg| \F_{T_u^A} \bigg], \quad u > 0.
\end{equation}
Then, independent of $n$, we have $\Prob^\alpha$-a.s.\ that
\begin{align} \label{slm4.1.2}
&\lim_{u \to \infty} {\mathbb E}^{\alpha} \left[ \cG_u | {\cal F}_n \right] {\bf 1}_{\{ n \le T_u^A\}}\nonumber\\
&  =\int_{{\mathbb S}^{d-1}_+ \times \reals_+} 
    \frac{e^{-\alpha s}}{r_\alpha^A(x)}  {\mathbb E} \Big[ g\big( e^{S_0} X_0, \ldots, e^{S_m} X_m \big) \, \big| \, X_0=x, \, S_0=s+ \log d_A(x)  \Big] \, \rho^A(dx,ds).   
\end{align}
\end{enumerate}
\end{lemma}

\noindent
{\bf Proof of Lemma \ref{sublemma}.}
 {\rm (i)}  By Lemma \ref{lem:geometric return}, $\tau$ satisfies the assumptions in Lemma \ref{lem:convergence Z} (iii). Thus, this result  is a direct consequence of Lemma \ref{lem:convergence Z} (iv), where the $L^1$-convergence $\abs{Z_n}^\alpha  {\bf 1}_{\{ n \le \tau\}} \to \abs{Z}^\alpha {\bf 1}_{\{ \tau=\infty\}}$ is proved.
It follows that $\abs{Z_n}^\alpha  {\bf 1}_{\{ n \le \tau\}}$ constitutes a Cauchy sequence in $L^1$, which yields the assertion.

{\rm (ii)}  
Let $n \in \pintegers$.  Then by the Markov property, 
$$ \E^\alpha  \left[ \cG_u | {\cal F}_n \right] {\bf 1}_{\{ n \le T_u^A \}} ~=~ \E^\alpha_{X_n,V_n} \big[ \cG_u \big] {\bf 1}_{\{ n \le T_u^A \}} \qquad \Prob^\alpha\text{-a.s.}$$
As $\lim_{u \to \infty}  {\bf 1}_{\{n \le T_u^A \}} =1$ $\Prob^\alpha$-a.s.,  it suffices to determine
 $\lim_{u \to \infty} \E_{x,v}^\alpha \big[ \cG_u]$ and show that this quantity is independent of $x$ and $v$.
For all $v \in \Rd_+$ and $u >0$, set
$$ G_u(v) ~=~ \E \bigg[  g\bigg(v, \Pi_1v + \frac{V_1^{(0)}}{u}, \dots, \Pi_m v + \frac{V_m^{(0)}}{u} \bigg) \bigg], \quad G(v) ~=~ \E \bigg[  g\bigg(v, \Pi_1v, \dots, \Pi_m  v \bigg) \bigg],$$
where $\Pi_k := M_k \cdots M_1$ for $k \ge 1$, and  $V_k^{(0)} := \sum_{i=1}^k M_k \cdots M_{i+1} Q_i$ for $k \ge 2$ and $V_1^{(0)} := Q_1$.
Now consider the decomposition:
\begin{align}
\E_{x,v}^\alpha \big[ \cG_u  \big] ~&=~ \E_{x,v}^\alpha \bigg[ \frac{1}{r_\alpha(X_{T_u^A})} \left( \frac{\big| V_{T_u^A} \big|}{u} \right)^{-\alpha} G_u\bigg( \frac{V_{T_u^A}}{u}  \bigg)  \bigg] \nonumber \\[.2cm]
~&=~ \E_{x,v}^\alpha \bigg[ \frac{1}{r_\alpha(X_{T_u^A})} \left( \frac{\big| V_{T_u^A} \big|}{u} \right)^{-\alpha}  \left( G_u\bigg( \frac{V_{T_u^A}}{u}  \bigg) -  G\bigg( \frac{V_{T_u^A}}{u}  \bigg) \right) \bigg] \nonumber \\[.2cm]
& \hspace*{1.5cm} ~+~ \E_{x,v}^\alpha \bigg[ \frac{r_\alpha(\widetilde{V}_{T_u^A})}{r_\alpha(X_{T_u^A})} \frac{1}{r_\alpha(\widetilde{V}_{T_u^A})} \left( \frac{\big| V_{T_u^A} \big| }{u} \right)^{-\alpha}  G\bigg( \frac{V_{T_u^A}}{u}  \bigg) \bigg] ~:=~{\mathbb I}_1(u) + {\mathbb I}_2(u). \label{eq:G_u converges}  
\end{align}

\Step[1]. We begin by showing that ${\mathbb I}_1(u) \to 0$ as $u \to \infty$. 
Write $g(v_0,\dots,v_m)=\hat{g}(v_0, \dots, v_m)\1[{\{|v_m| \ge \delta\}}]$, where $\hat{g}$ is $\theta$-H\"older continuous. 
Then 
\begin{align*}
&\left|  g \bigg( \bigg(\Pi_k v + \big(V_k^{(0)}/u \big) \bigg)_{k=0}^m \bigg) - g \Big( \big(\Pi_k v \big)_{k=0}^m \Big) \right| \\
& \hspace*{1.5cm} \le~
\left|\hat{g} \bigg( \Big(\Pi_k v + \big(V_k^{(0)}/u \big)\Big)_{k=0}^m \bigg)  - \hat{g} \Big( \big(\Pi_k v \big)_{k=0}^m \Big)  \right|
   \1[{\{|\Pi_m v | \ge \delta\}}] \\
& \hspace*{4cm} + |g|_\infty \left(\1[{[\delta,\infty)}]\bigg( \Big| \Pi_m v + \big(V_m^{(0)}/u \big) \Big| \bigg) - \1[{[\delta,\infty)}]\big( |\Pi_m v| \big) \right) \\[.2cm]
& \hspace*{1.5cm} \le~ \frac{1}{u^\theta}B_1 \Big(\sum_{k=1}^m  \big| V_k^{(0)} \big|^\theta \Big) +  |g|_\infty \left(\1[{[\delta,\infty)}]\bigg( \Big| \Pi_m v + \big(V_m^{(0)}/u \big) \Big| \bigg) - \1[{[\delta,\infty)}]\big( |\Pi_m v| 
\big) \right),
\end{align*}
for some constant $B_1$ arising from the $\theta$-H\"older continuity of $\hat{g}$. Let $(M^*,Q^*)$ be a pair of random variables which are independent of the sequence $\{(M_n, Q_n)\}$,
where the law of $(M^*,Q^*)$ is given by $\P{\big(\Pi_m, V_m^{(0)}\big) \in \cdot}$ under $\Prob^\alpha$. Setting $B_2:=\max_{y \in \SPdp} 
\left( r_\alpha(y) \right)^{-1}$ and using that $\big(|V_{T_u^A}|/u\big)^{-\alpha} <1$,
we obtain that
\begin{align}
{\mathbb I}_1(u) \le \frac{1}{u^\theta} \, B_1 B_2 \,  \E \Big[ \sum_{k=1}^m  \big| V_k^{(0)} \big|^\theta\Big] + B_2 \, |g|_\infty \left( \Prob_{x,v}^\alpha \bigg( \abs{M^* \frac{V_{T_u^A}}{u} + \frac{Q^*}{u}} \ge \delta \bigg) - \Prob_{x,v}^\alpha \bigg( \abs{M^* \frac{V_{T_u^A}}{u} } \ge \delta \bigg)\right). \label{decomposing hoelder}
\end{align}
Since the $\theta$-moment of $V_k^{(0)}$ is finite,
the first term tends to zero as $u \to \infty$.
 For the second term, we use the $\Prob^\alpha$-convergence $(M^*,Q^*/u) \Rightarrow (M^*,0)$ and $\big(\widetilde{V}_{T_u^A}^A, \log |V_{T_u^A}^A| - \log u\big) \Rightarrow \rho^A$ (by Theorem \ref{nonlinear renewal TuA}).
Let $(X,S)\sim \rho^A$ be a random vector independent of $(M^*,Q^*)$ under $\Prob^\alpha$. Using that $V_n=d_A(\widetilde{V}_n^A) V_n^A$ (cf.\ \eqref{def:VnA}), we have
 $V_{T_u^A}/u \Rightarrow d_A(X)e^S X$. [Here $X$ describes the limiting direction of $V_{T_u^A}^A/u$ and $S$ the limiting logarithmic overjump, as $\log |V_{T_u^A}^A| - \log u \Rightarrow S$.] 
  Since the sequences $\{(M^*,Q^*/u)\}$ and $\{V_{T_u^A}/u\}$ are independent, they also converge jointly in distribution. Hence, under $\Prob^\alpha$,
$$ M^* \frac{V_{T_u^A}}{u} + \frac{Q^*}{u}~\Rightarrow~d_A(X)e^S M^* X \quad \text{ and } \quad M^* 
\frac{V_{T_u^A}}{u} ~\Rightarrow~d_A(X)e^S M^*X. $$
Thus, the second term in \eqref{decomposing hoelder}  vanishes if $[\delta, \infty)$ is a continuity set for $d_A(X) e^S | M^*  X|$.

We now show that $[\delta,\infty)$ is a continuity set.
Since $M^*$ is independent of $(X,S)$, it suffices to show that for any allowable matrix ${\mathfrak m}$, $d_A(X) e^S|\mathfrak{m}  X|=\delta$ has probability 0.  Now for  each fixed $y \in \SPdp$, the equation ${\mathfrak h}(s) := d_A(y) e^s | {\mathfrak m} y|=\delta$ has a unique solution $s_y
\in \reals$. Hence
$$ \Prob^\alpha\big(  d_A(X) e^S | {\mathfrak m} X| =\delta \big)~=~\int_{\SPdp \times \R} \1[{\{s=s_y\}}] \, \rho^A(dy,ds) ~=~0,$$
since the radial component of the overjump distribution is absolutely continuous with respect to Lebesgue measure (as can be seen from the representation of $\rho^A$ in Eq.\ (1.16) of \shortciteN{Kesten1974}, or Eq. \eqref{eq:markov delay distribution} below,
which is valid for $\rho^A$ upon replacing $S$ by $S^A$ everywhere).

Thus, having shown that $[\delta, \infty)$ is a continuity set, we conclude by the Portmanteau theorem that for all $x \in \SPdp$ and $v \in \Rdo$,
\begin{align*}
&\Prob_{x,v}^\alpha \bigg( \abs{M^* \frac{V_{T_u^A}}{u} + \frac{Q^*}{u}} \ge \delta \bigg) - \Prob_{x,v}^\alpha \bigg( \abs{M^* \frac{V_{T_u^A}}{u} } \ge \delta \bigg) \\[.2cm]
& \hspace*{2cm} ~\to~\Prob^\alpha\Big(d_A(X) e^S |M^*X|\ge \delta\Big) - \Prob^\alpha \Big(d_A(X) e^S |M^*X|\ge \delta\Big),
\end{align*}
and hence also the second member of \eqref{decomposing hoelder} vanishes as $u \to \infty$. Thus ${\mathbb I}_1(u) \to 0$ as $u \to \infty$.

\medskip

\Step[2].  Now turn to ${\mathbb I}_2(u)$.
Using Theorem \ref{nonlinear renewal TuA}, 
again invoke the convergence $\big(({V_{T_u^A}^A})^\sim, \log |V_{T_u^A}^A| - \log u\big) \Rightarrow \rho^A$ under $\Prob^\alpha$.
Moreover, by Lemma \ref{lem:VTutoXTu}, using the  continuity and boundedness of $r_\alpha$, we have that
$r_\alpha(\widetilde{V}_{T_u})/r_\alpha(X_{T_u})$ tends to one in $\Prob^\alpha$-probability.  Hence by Slutsky's theorem, the quantity inside ${\mathbb I}_2(u)$ 
converges in law, and identifying this limit distribution, we deduce that
\begin{align*} 
\hspace*{.5cm} \lim_{u\to \infty} {\mathbb I}_2(u) =& \lim_{u \to \infty}  \E_{x,v}^\alpha  \bigg[ \frac{r_\alpha(\widetilde{V}_{T_u^A})}{r_\alpha(X_{T_u^A})} \frac{1}{r_\alpha^A(\widetilde{V}_{T_u^A}^A)} \left( \frac{\big| V^A_{T_u^A}\big|}{u} \right)^{-\alpha}   G\bigg( \frac{V_{T_u^A}}{u}  \bigg) \bigg] \\
 ~=&~ \int_{{\mathbb S}^{d-1}_+ \times \reals_+} \frac{e^{-\alpha s}}{r^A_\alpha(y)}    G\Big(d_A(y) e^s y\Big) \rho^A(dy,ds)  \\[.2cm]
=&~ \int_{{\mathbb S}^{d-1}_+ \times \reals_+} 
    \frac{e^{-\alpha s}}{r_\alpha^A(y)}  {\mathbb E} \Big[ g\big( e^{S_0} X_0, \ldots, e^{S_m} X_m \big) \, \big| \, X_0=y, \, S_0=s+ \log d_A(y)  \Big] \, \rho^A(dy,ds) . \hspace*{.5cm} \halmos
\end{align*}\\[-.3cm]

\textbf{\textit{The dual change of measure.}}
Prior to proving  Proposition \ref{prop1}, we  introduce a ``dual" change of measure, where the process $\{ (M_n, Q_n):  n =1,2,\ldots \}$
follows the $\alpha$-shifted measure
until the {\it random} time $T_u^A$, and follows the original measure thereafter.
We shall denote expectation relative to this dual measure by ${\mathbb E}^{\D} [ \cdot ]$.  More formally, for any $n \in {\mathbb N}$, define
\begin{align} \label{dual-def}
&{\mathbb E}^{\D}\left[ h(V_0,M_1,Q_1,\ldots,M_n,Q_n) \big|
 X_0 = x, \: V_0 = v \right] \nonumber\\[.2cm]
& \hspace*{1cm} = {\mathbb E} \left[ \frac{ | M_{n \wedge T_u^A} \cdots M_1 x|^\alpha 
  r_\alpha(X_{T_u^A \wedge n})}{r_\alpha(x)} h(V_0,M_1,Q_1,\ldots,M_n,Q_n) \Big|
 X_0 = x, \: V_0 = v
  \right],
\end{align}
for all measurable functions $h:  \realsd \times  ({\mathfrak M} \times \realsd)^n \to \reals$.   Following the notational conventions
of the previous sections, we  write 
 $\E_\gamma^\mathcal{D} [ \cdot] = \int \E^\mathcal{D}[ \cdot | X_0 = \widetilde{v}, V_0 = v] \, \gamma(dv)$ for any probability measure $\gamma$ on $\Rd_+\setminus\{0\}$.

\medskip

\noindent{\bf Proof of Proposition \ref{prop1}.}
Note that $\{ V_n \}$ is transient in the $\alpha$-shifted measure and thus $T_u^A < \infty$ a.s.;
cf.\ Lemma \ref{lemma-charac-T}.
Hence, employing the dual change of measure in \eqref{dual-def} over the random time interval $[0, T_u^A]$ yields that
\begin{align*} 
&\hspace*{-1cm} u^\alpha {\mathbb E} \left[ \left. g \left( \frac{V_{T_u^A}}{u}, \ldots, \frac{V_{T_u^A +m}}{u} \right) {\bf 1}_{\{ T_u^A < \tau \}} \right| V_0 = v  \right] \\[.2cm]
=&~ u^\alpha r_\alpha(\widetilde{v}) {\mathbb E}_{\bfDV}^{\D}  \left[ \frac{e^{-\alpha S_{T_u^A}}}{r_\alpha(X_{T_u^A})} {\bf 1}_{\{ T_u^A < \tau \}} \E \bigg[ g \left( \frac{V_{T_u^A}}{u}, \ldots, \frac{V_{T_u^A +m}}{u} \right) \bigg| \F_{T_u^A} \bigg]  \right] .
\end{align*}
Now substitute the quantity  ${\mathfrak G}_u$, defined in \eqref{def-cG}, into the previous equation.
Noting that
\[
Z_n := \frac{V_n}{|M_n \cdots M_1 X_0|} := \frac{V_n}{e^{S_n}}, \quad n \in \pintegers, 
\]
we obtain after a little algebra that
\begin{equation} \label{lm4.1.1}
u^\alpha {\mathbb E} \left[ \left. g \left( \frac{V_{T_u^A}}{u}, \ldots, \frac{V_{T_u^A +m}}{u} \right) {\bf 1}_{\{ T_u^A < \tau \}} \right| X_0 = \widetilde{v}, \:V_0 = v  \right] =  r_\alpha(\widetilde{v}) {\mathbb E}_{\bfDV}^{\D}  \Big[ |Z_{T_u^A}|^\alpha \cG_u {\bf 1}_{\{ T_u^A < \tau \}} \Big].
\end{equation}
The right-hand side can be further equated, for $n \in \pintegers$, to
\begin{align} \label{lm4.1.2}
 &~ r_\alpha(\widetilde{v}) {\mathbb E}_{\bfDV}^{\D}  \left[ \bigg( \big| Z_{T_u^A} \big|^\alpha {\bf 1}_{\{ T_u^A < \tau\}} - \abs{Z_n}^\alpha {\bf 1}_{\{n \le T_u^A\}} {\bf 1}_{\{ n \le \tau\}} \bigg) \cG_u \right] \nonumber\\[.2cm]
& \hspace*{2.5cm} + r_\alpha(\widetilde{v}) {\mathbb E}_{\bfDV}^{\alpha}  \Big[ \abs{Z_n}^\alpha {\bf 1}_{\{n \le T_u^A\}} {\bf 1}_{\{ n \le \tau\}} \, \E^{\alpha}  \big[ \cG_u \, | \F_n \big] \Big], 
\end{align}
where we have replaced ${\mathbb E}^\alpha_{\bfDV}[\cdot | {\cal F}_n ]$ with ${\mathbb E}^\alpha [\cdot | {\cal F}_n]$ in the last expectation, since this conditional expectation depends only on $(X_n,V_n)$,
and not on the initial values $(X_0,V_0)$ once $(X_n,V_n)$ has been specified.
Moreover, we have replaced the dual change of measure by the $\alpha$-shifted measure, since they coincide for random variables which are 
$\F_{T_u^A}$-measurable.

To analyze the quantity in \eqref{lm4.1.2}, we first take the limit as $u \to \infty$ and then as $n \to \infty$. 
By part {\rm (i)} of Lemma \ref{sublemma} and the boundedness of $\cG_u$, we deduce from \eqref{lm4.1.1} and \eqref{lm4.1.2} that
\begin{align} 
& \lim_{u \to \infty}  u^\alpha {\mathbb E} \left[  \left. g \left( \frac{V_{T_u^A}}{u}, \ldots, \frac{V_{T_u^A +m}}{u} \right) {\bf 1}_{\{ T_u^A < \tau \}} \right|  V_0 = v  \right] \nonumber \\[.1cm] 
& \hspace*{1.5cm} = \lim_{n \to \infty} \lim_{u \to \infty} 
r_\alpha(\widetilde{v}) {\mathbb E}_{\bfDV}^{\alpha}  \Big[ \left. \abs{Z_n}^\alpha {\bf 1}_{\{n \le T_u^A\}} {\bf 1}_{\{ n \le \tau\}} \, \E^{\alpha} \big[ \cG_u \, \right| \F_n \big] \Big] \label{lm4.1.3}.
\end{align}
Now by Lemma \ref{lem:convergence Z} (iii), $\big\{| Z_n|^\alpha {\bf 1}_{\{ n \le \tau \}}\big\}$  is uniformly integrable.  Denote by $\cG$ the right-hand side of \eqref{slm4.1.2}. Since $\cG_u$ is bounded by $b^{-1}\abs{g}_\infty$
and $T_u^A \uparrow \infty$ ${\mathbb P}^\alpha$-a.s., it follows by Lemma \ref{sublemma} {\rm (ii)} that
\begin{align} 
 \lim_{u \to \infty}  u^\alpha {\mathbb E} & \left[ \left. g \left( \frac{V_{T_u^A}}{u}, \ldots, \frac{V_{T_u^A +m}}{u} \right) {\bf 1}_{\{ T_u^A < \tau \}} \right| V_0 = v  \right] \nonumber \\[.2cm]
~=~& \lim_{n \to \infty}
r_\alpha(\widetilde{v}) {\mathbb E}_{\bfDV}^{\alpha}  \left[ \abs{Z_n}^\alpha {\bf 1}_{\{ n \le \tau\}} \, \lim_{u \to \infty}   {\bf 1}_{\{n \le T_u^A\}} \E^{\alpha} \big[ \cG_u \, | \F_n \big] \right] \nonumber\\[.2cm]
~=~& 
r_\alpha(\widetilde{v}) {\mathbb E}_{\bfDV}^{\alpha}  \left[ \lim_{n \to \infty} \abs{Z_n}^\alpha {\bf 1}_{\{ n \le \tau\}} \, \mathfrak{G} \right]  ~=~ r_\alpha(\widetilde{v}) \E_{\bfDV}^\alpha \big[ \abs{Z}^\alpha {\bf 1}_{\{\tau=\infty\}} \big] \, \cG  \label{lm4.1.4},
\end{align}
where ${\mathfrak G}$ is the limit appearing on the right-hand side of \eqref{slm4.1.2}.
\halmos

In some cases, it is useful to consider functions $g$ which depend on the {\it infinite} path $(V_{T_u^A},V_{T_u^A+1}, \ldots)$,
or to consider functions $g$ which need not be bounded.  Moreover, it is also useful to have {\it uniform} upper bounds.  In these situations, a variant of the above proposition is useful.

\begin{prop}\label{subcorollary}
Suppose that $g : (\Rd_+)^\N \to [0, \infty)$ is a nonnegative measurable function, and set
$$ \bar{\cG}_u ~=~  \frac{1}{r_\alpha \big( X_{T_u^A} \big)} \left( \frac{\big| V_{T_u^A}\big|}{u} \right)^{-\alpha} \E \bigg[ g\left( \bigg(\frac{V_{T_u^A+k}}{u}\bigg)_{k \ge 0} \right) \, \bigg| \, \mathcal{F}_{T_u^A} \bigg]. $$
Further assume that for some finite constant $B$ and some ${\cal U} \ge 0$,
\begin{equation}\label{eq:bound Gu}
\sup_{u \ge {\cal U}} \bar{\mathfrak G}_u \le B \quad {\mathbb P}^\alpha\mbox{\rm -a.s.}
\end{equation}
Then for any bounded set $F \subset \R_+^d\setminus\{0\}$, there exists a finite constant $L$, not depending on $B$, such that
\begin{align} 
0 ~\le~ \sup_{u \ge {\cal U}}\, \sup_{v \in F} \,   u^\alpha {\mathbb E} \left[ g\left( \bigg(\frac{V_{T_u^A+k}}{u}\bigg)_{k \ge 0} \right) {\bf 1}_{\{ T_u^A < \tau \}} \bigg|  V_0 = v  \right] 
  ~\le~  B L.
    \label{cor4.3b}
\end{align}
Moreover, if \eqref{eq:bound Gu} holds and $\limsup_{u \to \infty} \bar{\cG}_u=0$  $\Prob^\alpha$-a.s., then we also have that
\begin{equation}\label{cor4.3c} \lim_{u \to \infty} u^\alpha {\mathbb E} \left[ g\left( \bigg(\frac{V_{T_u^A+k}}{u}\bigg)_{k \ge 0} \right) {\bf 1}_{\{ T_u^A < \tau \}} \bigg|  V_0 = v  \right]  = 0.
\end{equation}
\end{prop}

\noindent {\bf Proof.}
Repeating the argument in the proof of Proposition \ref{prop1} leading to \eqref{lm4.1.3}, we obtain that
\begin{align*}
0 & \le \sup_{u  \ge {\cal U}} \sup_{v \in F}  u^\alpha {\mathbb E} \left[ g\left( \bigg(\frac{V_{T_u^A+k}}{u}\bigg)_{k \ge 0} \right) {\bf 1}_{\{ T_u^A < \tau \}} \bigg| V_0 = v  \right]  ~\le~ B  \es[\alpha](\widetilde{v}) \sup_{v \in F} \E_{\bfDV}^\alpha \left[\sup_{n \in \N} |Z_n|^\alpha  {\bf 1}_{\{ n \le \tau \}} \right] 
\end{align*}
which is finite by Lemma \ref{lem:convergence Z} (iii) and the boundedness of $r_\alpha$. The boundedness of ${\bar{\cG}_u}$ further allows us to use the dominated convergence theorem in order to deduce \eqref{cor4.3c} from \eqref{lm4.1.3}.
\halmos

The remainder of this section is devoted to the proof of several results that are needed in order to establish Theorem \ref{thm2A} in the subsequent section. Therefore, we now restrict our attention to the case
where $A=\{ v \in \R^d_+ \, : \, |v|>1\}$; and thus, $T_u^A=T_u$, $d_A =1$, $r_\alpha^A=r_\alpha$, and $\rho^A=\rho$.

 Recall that in Proposition \ref{prop1}, we studied the behavior of $\{ V_n \}$ over paths of finite length.  Now suppose that we
replace the function $g$ in that proposition with  a function of the form
$\sum_{k=0}^{\tau - 1 - T_u} h\big( V_{T_u+k}/u \big)$.  In the next lemma, we show that it is the path behavior over {\it finite}
time intervals of the form $[T_u,T_u+m]$ which plays the determining role.  Additionally, we establish a technical result, given in 
\eqref{eq:sum tau to m} below, stating that if one computes the path behavior over $[T_u,T_u+m]$ when $T_u - \tau < m$,
then the effect of these additional terms is, roughly speaking, negligible. 

\begin{lemma}\label{sublemma2}
Let $h$ be a bounded measurable function such that $h(x)=0$ for all $x \in B_\delta(0)$, for some $\delta >0$.  Then for all $v \in \Rd_+\setminus\{0\}$ and all $m \in \N$,
\begin{equation} \label{eq:sum m to tau}
\lim_{m \to \infty} \lim_{u \to \infty} u^\alpha {\mathbb E} \left[ \sum_{k=m}^{\tau - 1 - T_u} h\bigg( \frac{V_{T_u+k}}{u} \bigg) {\bf 1}_{\{ T_u + m < \tau \}} \bigg| V_0 = v  \right] ~=~0.
\end{equation}
Moreover, if $F \subset \Rdo$ is bounded, then  by summing over {\em all} terms in the interval $[T_u,\tau)$, we obtain that
\begin{equation}\label{eq:bound for thm2.1}
 \limsup_{u \to \infty} \sup_{v \in F} u^\alpha {\mathbb E} \left[ \sum_{i=0}^{\tau-1} h \left(\frac{V_i}{u} \right) {\bf 1}_{\{ \abs{V_i} > u \}} \bigg| V_0 = v  \right]  < \infty.
\end{equation}
Furthermore,  for all $v \in \Rd_+\setminus\{0\}$,
\begin{equation} \label{eq:sum tau to m}
\lim_{m \to \infty} \lim_{u \to \infty} u^\alpha {\mathbb E} \left[ \sum_{i=\tau}^{T_u + m} h\bigg( \frac{V_{i}}{u} \bigg) 
{\bf 1}_{\{ T_u < \tau \le T_u + m \}} \bigg| V_0 = v  \right] ~=~ 0. 
\end{equation}
\end{lemma}

\noindent
{\bf Proof.}
\Step[1]. First we establish \eqref{eq:sum m to tau}.
By Eq.\ \eqref{cor4.3b} in Proposition \ref{subcorollary}, it suffices to prove that 
$$  \sup_{u \ge \cal U} \bar{\cG}_u ~:=~ \sup_{u \ge \cal U}  \frac{1}{r_\alpha (X_{T_u})} \left( \frac{\big|V_{T_u}\big|}{u} \right)^{-\alpha} \E \bigg[ \sum_{k=m}^{\tau-1-T_u} \abs{h\bigg( \frac{V_{T_u+k}}{u} \bigg)} {\bf 1}_{\{T_u +m < \tau \}} \bigg| \F_{T_u} \bigg] ~\le~B(m,\cal U)$$
 for a sequence $\{B(m,\cal U)\}$ which tends to zero as we first let $\cal U \to \infty$ and then let $m \to \infty$.
By employing  the Markov property and the boundedness of $r_\alpha$, we see that it is enough to show that, for a suitable sequence $B(m, \cal U)$,
\begin{align} \label{eq:uniform bound Huv}
\sup_{u \ge \cal U} \sup_{v : |v| \ge u} H_u(v) ~:=~ \sup_{u \ge \cal U} \sup_{v : |v| \ge u} \E \bigg[ \left( \frac{\abs{V_{0}}}{u} \right)^{-\alpha}  \sum_{k=m}^{\tau-1} h\bigg( \frac{V_{k}}{u} \bigg) {\bf 1}_{\{m < \tau \}} \bigg| V_0=v \bigg] ~\le~B(m, \cal U).
\end{align}
Now  let ${\DD}^\dagger=\{ v \in \R^d_+ \, : \, |v| \le L\}$ be defined as in Lemma \ref{lm-drift}, and set  ${\tau}^\dagger:=\inf\{ n \in \pintegers \, : \, V_n  \in \DD^\dagger\}$.  Recall that $h(x)=0$ for all $x \in B_\delta(0)$.
Hence, for $0 < \theta < \min\{1, \alpha\}$ and $|v| > u$, we have
\begin{align}
H_u(v) ~&\le~ \abs{h}_\infty \E \bigg[ \left( \frac{\abs{V_{0}}}{u} \right)^{-\alpha} \left(  \sum_{k=m}^{\tau^\dagger-1} {\bf 1}_{\{ \abs{V_k} > \delta u \}} + \sum_{k={\tau^\dagger}}^{{\tau}-1} {\bf 1}_{\{ \abs{V_k} > \delta u \}} \right) \bigg| V_0=v \bigg] \nonumber \\
~&\le~ \abs{h}_\infty \left(\frac{\abs{v}}{u}\right)^{-\alpha} \left( \sum_{k=m}^\infty \P{\abs{V_k} > \delta u, \: \tau^\dagger > k \, \big| V_0=v} + \E \bigg[ \E \bigg[  \sum_{k=\tau^\dagger}^{\tau} {\bf 1}_{\{ \abs{V_k} > \delta u \}} \, \bigg| \, {\cal F}_{\tau^\dagger} \bigg] \bigg]\right) \nonumber \\
~&\le~ \abs{h}_\infty \left(\frac{\abs{v}}{u}\right)^{-\alpha} \sum_{k=m}^\infty (\delta u)^{-\theta} \E_v \big[ \abs{V_k}^\theta {\bf 1}_{\{\tau^\dagger >k\}} \big]  +  \abs{h}_\infty \sup_{w \in {\DD}^\dagger} \E_w \bigg[  \sum_{k=0}^{\tau} {\bf 1}_{\{ \abs{V_k} > \delta u \}} \bigg] .  \nonumber \label{eq:estimates F_u(v)}
\end{align}
The first sum can be estimated further by employing Lemma \ref{lm-drift}:
$$  \sup_{v : |v| \ge u} \left(\frac{\abs{v}}{u}\right)^{-\alpha} \sum_{k=m}^\infty (\delta u)^{-\theta} \E_v \big[ \abs{V_k}^\theta {\bf 1}_{\{{\tau}^\dagger >k\}} \big] ~\le~ \frac{B}{\delta^\theta} \left( \sup_{v : |v| \ge u} \left(\frac{\abs{v}}{u}\right)^{\theta-\alpha}\right) \frac{t^{m}}{1-t} ~=~ \frac{B}{\delta^\theta} \frac{t^{m}}{1-t},$$
and the last term tends to $0$ as $m \to \infty$.
For the second term, note that \eqref{eq:geometric return} implies that
$\sup_{w \in {\DD}^\dagger} \E[ \tau \, | \, V_0=w] < \infty$.  Hence we can apply a dominated convergence argument to infer that
$$ \sup_{u \ge \cal U} \sup_{w \in {\DD}^\dagger} \E_w \bigg[  \sum_{k=0}^{\tau} {\bf 1}_{\{ \abs{V_k} > \delta u \}}  \bigg] \le \sup_{w \in {\DD}^\dagger} \E_w \bigg[  \sum_{k=0}^{\tau} {\bf 1}_{\{ \abs{V_k} > \delta \cal U \}} \bigg] ~\to~0  \text{ as } \cal U \to \infty.$$
Combining these estimates, we have established \eqref{eq:uniform bound Huv}, and \eqref{eq:sum m to tau} follows. 

Finally,  \eqref{eq:bound for thm2.1} is a direct consequence of the above estimates for $m=0$ combined with \eqref{cor4.3b}.

\medskip

\Step[2]. Turning to \eqref{eq:sum tau to m}, we now  apply the second part of Proposition \ref{subcorollary}.   
Using that $h = 0$ on $B_\delta(0)$, 
it is now sufficient to show that for any fixed $m \in \N$,
\begin{equation} \label{sec4-ADD100a}
\bar{{\mathfrak G}}_u ~:=~ \frac{1}{r_\alpha (X_{T_u})} \left( \frac{\abs{V_{T_u}}}{u} \right)^{-\alpha} \E \left[ \sum_{i=\tau}^{T_u + m} \abs{h\bigg( \frac{V_{i}}{u} \bigg)} 
  {\bf 1}_{\{ |V_i| > \delta u \}} {\bf 1}_{\{T_u 
  < \tau \le T_u + m \}} \bigg| \F_{T_u} \right] 
\end{equation}
is bounded uniformly in $u$ and tends  to zero $\Prob^\alpha$-a.s.\ as $u \to \infty$.  
The prefactors are bounded, and thus  it suffices to estimate
\begin{align} \label{sec4-ADD101}
\E \bigg[ \sum_{i = \tau}^{T_u + m} & \abs{h\bigg( \frac{V_{i}}{u} \bigg)} {\bf 1}_{\{ |V_i| > \delta u \}}
 {\bf 1}_{\{T_u < \tau \le T_u + m \}} \bigg| \F_{T_u} \bigg] \nonumber\\[.1cm]
  ~&\le~ \abs{h}_\infty \E \left[ \E \left[ \sum_{i=\tau}^{\tau + m} {\bf 1}_{\{ \abs{V_{i}} > \delta u\}} \bigg| \F_{\tau} \right] {\bf 1}_{\{ T_u < \tau\}}  \bigg| \F_{T_u} \right].
\end{align} 
Now let $\theta \in (0,\alpha)$.  Then for all $k=0, \ldots,m$,
\begin{align*}
 {\mathbb P} &\left( \left. |V_{\tau + k}| > \delta u \right| V_{\tau} = v \right) ~\le~  \sup_{v \in \DD} \,  {\mathbb P} \left( \left. |V_k| > \delta u \right| V_0 = v \right)     ~\le~ \sup_{v \in \DD} \,  (\delta u)^{-\theta} {\mathbb E} \left[ \left. |V_k|^\theta \right| V_0 = v \right]  \\
& \hspace*{1.5cm} \le \delta^{-\theta} u^{-\theta} \Big( \sup_{v \in \DD} \E\big[ \norm{M_k \cdots M_1}^\theta \big] \cdot |v|^\theta + \sum_{j=1}^k \E \big[ \abs{M_k \cdots M_{j+1} Q_j}^\theta \big] \Big) ~\le~ B_1 u^{-\theta}
\end{align*}
for some constant $B_1=B_1(m) < \infty$.  Substituting the last estimate into \eqref{sec4-ADD101}
and then \eqref{sec4-ADD100a}, we obtain that for some finite constant $B_2$,
\[
\bar{\mathfrak G}_u \le  B_2 m u^{-\theta} \searrow 0 \quad \mbox{as} \quad u \to \infty.
\]
Thus, by Proposition \ref{subcorollary}, we have for all $m \in \N$ that
$$\lim_{u \to \infty} u^\alpha {\mathbb E} \left[ \sum_{i=\tau}^{T_u + m} h\bigg( \frac{V_{i}}{u} \bigg) 
{\bf 1}_{\{ T_u < \tau \le T_u + m \}} \bigg| V_0 = v  \right] ~=~ 0$$
and \eqref{eq:sum tau to m} follows. 
\halmos

Now by combining Corollary \ref{cor:prop1 for Tu} and Lemma \ref{sublemma2}, we obtain the following result.
 In this lemma, the function $f$ will correspond to that function appearing in the statement of Theorem \ref{thm2A}.

\begin{lemma}\label{lem:Kesten result}
Assume Hypotheses $(H_1)$ and $(H_2)$ are satisfied.
Let $\theta \le \min\{1, \alpha\}$ and let $f$ be a nonnegative bounded $\theta$-H\"older continuous  function. Then for all $v \in \Rdo$,
\begin{align} \label{eq:convergence up to m}
& \lim_{u \to \infty}  u^\alpha {\mathbb E} \left[ \sum_{i=0}^{\tau-1} f \left(\frac{V_i}{u} \right) {\bf 1}_{\{ \abs{V_i} \ge u \}} \bigg| V_0 = v  \right] \nonumber \\ 
&\hspace*{1.25cm} =  \es[\alpha](\widetilde{v}) {\mathbb E}_{\bfDV}^{\alpha} \left[ |Z|^\alpha {\bf 1}_{\{ \tau = \infty \}} \right] 
\lim_{m \to \infty}
    \int_{{\mathbb S}^{d-1}_+ \times \reals_+} 
     {\mathbb E}_{x}^{\alpha} \left[  \sum_{i=0}^m F(X_i,S_i+s) \right] \, \rho(dx,ds),
\end{align}
where $F(x,s):= \big( {e^{-\alpha s}} f(e^sx)/{r_\alpha(x)} \big)  {\bf 1}_{[0,\infty)}(s)$.
\end{lemma}

\noindent{\bf Proof.}
Let $g(v):=f(v) {\bf 1}_{\{|v| \ge 1\}}$, and note that $g$ is an almost $\theta$-H\"older-continuous function. 
For $m \in \pintegers$, consider the decomposition 
\begin{align}
u^\alpha {\mathbb E}_v & \left[ \sum_{i=0}^{\tau-1} f \left(\frac{V_i}{u} \right) {\bf 1}_{\{ \abs{V_i} \ge u \}}   \right] ~=~ u^\alpha {\mathbb E}_v \left[ \sum_{i=T_u}^{\tau-1} f \left(\frac{V_i}{u} \right) {\bf 1}_{\{ \abs{V_i} \ge u \}} {\bf 1}_{\{T_u < \tau\}}   \right] \nonumber \\[.2cm]
 &=~   u^\alpha \sum_{i=0}^m {\mathbb E}_v \left[ g \left(\frac{V_{T_u+i}}{u} \right) {\bf 1}_{\{T_u < \tau\}}   \right]  
+  u^\alpha {\mathbb E}_v \left[ \sum_{i=m+1}^{\tau - 1 - T_u} g \left(\frac{V_{T_u +i}}{u} \right)  {\bf 1}_{\{ T_u + m < \tau\}}   \right] 
    \nonumber\\
 & \hspace*{6.5cm} - u^\alpha {\mathbb E}_v \left[ \sum_{i=\tau}^{T_u+m} g \left(\frac{V_{i}}{u} \right)  {\bf 1}_{\{ T_u < \tau \le T_u+m\}}   \right]. \label{eq:convergent terms}
\end{align} 
On the right-hand side of \eqref{eq:convergent terms}, the last two terms
tend to zero  by Lemma \ref{sublemma2} when taking first the limit $u \to \infty$ and then $m \to \infty$. 
Next, by Corollary \ref{cor:prop1 for Tu}, we obtain that 
\begin{align*} 
 \lim_{u \to \infty} & u^\alpha \sum_{i=0}^m {\mathbb E} \left[ g \left(\frac{V_{T_u+i}}{u} \right) {\bf 1}_{\{T_u < \tau\}} \bigg| V_0 = v
   \right] \\[.2cm]
~=~ &  \es[\alpha](\widetilde{v}) {\mathbb E}_{\bfDV}^{\alpha} \left[ |Z|^\alpha {\bf 1}_{\{ \tau = \infty \}} \right] 
   \sum_{i=0}^m \int_{{\mathbb S}^{d-1}_+ \times \reals_+}  {\mathbb E}_{x}^{\alpha} \left[   \frac{e^{-\alpha (S_i+s)}}{r_\alpha(X_i)} g(e^{S_i +s} X_i)   \right] \, \rho(dx,ds) \\[.2cm]
   ~=~ &  \es[\alpha](\widetilde{v}) {\mathbb E}_{\bfDV}^{\alpha} \left[ |Z|^\alpha {\bf 1}_{\{ \tau = \infty \}} \right] 
   \int_{{\mathbb S}^{d-1}_+ \times \reals_+}  {\mathbb E}_{x}^{\alpha} \left[   \sum_{i=0}^m \frac{e^{-\alpha (S_i+s)}}{r_\alpha(X_i)} f(e^{S_i +s} X_i) {\bf 1}_{\{ S_i +s \ge 0\}}  \right] \, \rho(dx,ds), 
\end{align*} 
and the assertion follows by letting $m \to \infty$. \halmos

To bring the limit (as $m \to \infty$) inside the sum in \eqref{eq:convergence up to m}, we first
need to introduce the definition of a multivariate directly Riemann integrable function.

\begin{defin}\label{def:dRi}
A measurable function $F : \SPdp \times \R \to \R$ is called {\it directly Riemann integrable} if for all $x \in \SPdp$, the function $t \mapsto F(x,t)$ is continuous a.e.\ with respect to Lebesgue measure on $\R$, and
\begin{equation}\label{eq:dRi sum} \sum_{l=-\infty}^\infty \sup \big\{ \abs{F(x,t)} \, : \, x \in \SPdp, t \in [l, l+1] \big\} < \infty. \end{equation}
\end{defin}

In particular, $F$ directly Riemann integrable implies that 
\begin{equation} \label{sec5-201}
 \sup_{x \in \SPdp} \sup_{ s\in \R} \ \E^\alpha_x \left[ \sum_{i=0}^\infty |F(X_i, S_i+s)| \right] < \infty;
 \end{equation} 
cf.\ \shortciteN{Mentemeier2013a}, Section 6.1.  [In fact, the above definition of direct Riemann integrability implies Definition 1 in \shortciteN{Kesten1974}; cf.\ \shortciteN{Buraczewski.etal:2014},
Lemma C.1.  Then \eqref{sec5-201} follows from Lemma 6 of \shortciteN{Kesten1974}.]
Using this definition, we obtain:
\begin{prop}\label{prop:Kesten result}
Under the assumptions of Lemma \ref{lem:Kesten result}, $F(x,s):= \big( e^{-\alpha s} f(e^s x)/r_\alpha(x) \big)
  {\bf 1}_{[0,\infty)}(s)$ is directly Riemann integrable, and thus
\begin{align} \label{conv-limit}
 \lim_{u \to \infty} & u^\alpha {\mathbb E} \left[ \sum_{i=0}^{\tau-1} f \left(\frac{V_i}{u} \right) {\bf 1}_{\{ \abs{V_i} \ge u \}} \bigg| V_0 = v  \right] 
=  \es[\alpha](\widetilde{v}) {\mathbb E}_{\bfDV}^{\alpha} \left[ |Z|^\alpha {\bf 1}_{\{ \tau = \infty \}} \right] 
    \int
     {\mathbb E}_{x}^{\alpha} \left[  \sum_{i=0}^\infty F(X_i,S_i+s) \right] \, \rho(dx,ds).
\end{align}
\end{prop}

\noindent{\bf Proof.} Since $r_\alpha$ is bounded from  below, it follows  that for some positive constant $b$,
\[
\overline{F}(s):= \sup_{x \in {\mathbb S}^{d-1}_+} \left| F(x,s) \right|  \le \frac{1}{b}  \abs{f}_\infty e^{-\alpha s} {\bf 1}_{[0,\infty)}(s).
\]
Since the right-hand side is a decreasing integrable function, we conclude that $\overline{F}$ is (univariate) directly Riemann integrable.  
Since $\overline{F}$ is obtained from $F$ by taking the supremum over all $x \in \SPdp$,
it follows immediately from  Definition \ref{def:dRi}
that $F$ is (multivariate) Riemann integrable.  Then by \eqref{sec5-201}, we can use a dominated convergence
argument to interchange $\lim_{m \to \infty}$ with the integration in \eqref{eq:convergence up to m}.
\halmos

%
%
%
%

\setcounter{equation}{0}
\section{Proof of Theorem \ref{thm2A}}\label{sect:5}
In this section, we provide the proof of Theorem \ref{thm2A}, first under the additional hypothesis $(H_3)$ 
of Section \ref{sect:3}, which is then removed by approximating $\{V_n\}$ from above and below by smoothed processes for which $(H_3)$ is satisfied. 

To establish Theorem \ref{thm2A}, we apply Proposition \ref{prop:Kesten result} directly, except that we must identify the integral in \eqref{conv-limit}.  
Note that if we were in the setting of classical random walk---where $\{ \xi_i \}$ is an i.i.d.\ sequence of random variables and $S_n = \sum_{i=1}^n \xi_i +S_0$---and
if we were to consider $F(s) = {\bf 1}_{[0,t]}(s)$ in \eqref{conv-limit}, then the integral in \eqref{conv-limit} would represent the renewal function, with $S_0$ having the
stationary excess distribution $\rho$.  It is known (see e.g. Proposition 1.0, Eq.\ (1.4), in \shortciteN{Thorisson1987}) that for $S_0 \sim \rho$, the renewal function is {\em equal} to $t/\E[\xi_1]$. In other words, if $S_0 \sim \rho$, then the renewal measure $\sum_{i=0}^\infty \P{S_i \in \cdot}$, restricted to $[0,\infty)$, is equal to $1/\E[\xi]$ 
multiplied by  Lebesgue measure  on $[0,\infty)$.

Our present objective is to extend this identity
into our Markovian framework. Here,  $\E[\xi]$ must be replaced by the drift of $S_1$ under $\Pstat$, which is $\lambda'(\alpha)$, and instead of Lebesgue measure on $\R_+$, we expect to obtain the limiting measure from Kesten's renewal theorem, namely the measure $\eta_\alpha \otimes {\mathfrak l}$ on $\SPdp \times \R_+$, where ${\mathfrak l}$ denotes Lebesgue measure.

\begin{lemma}\label{lem:markovdelay}
Let $g : \SPdp \times \R \to \R$ be a directly Riemann integrable function. 
Then
\begin{equation}\label{eq:markovdelay} 
\int_{{\mathbb S}^{d-1}_+ \times \reals_+}  \E_x^\alpha \left[ \sum_{i=0}^\infty g(X_i, S_i+s) 
 {\bf 1}_{\{ S_i + s  \ge 0 \}}  \right] \, \rho(dx,ds) ~=~
   \frac{1}{\lambda^\prime(\alpha)} \int_{{\mathbb S}^{d-1}_+ \times \reals_+}  g(x,s) \, \eta_\alpha(dx) \, ds .
\end{equation}
\end{lemma}

Before proving this result,
we start by recalling the standard extension of $\{ (X_n, S_n):  n=0,1,\ldots \}$ to a doubly-infinite stationary process, which, in particular,
may be used to identify the measure
 $\rho$ appearing on the left-hand side of \eqref{eq:markovdelay}. 
 For this purpose, let
 $\xi_n:=S_n - S_{n-1}$ and recall that  $\{ (X_n,\xi_n) :  n=1,2,\ldots\}$ is stationary under $\Pstat$; cf.\ \eqref{stationary-theta}.
Now let $\big\{ (X_n^\sharp, \xi_n^\sharp):  -\infty < n < \infty \big\}$ be the two-sided extension of this stationary sequence $\{ (X_n,\xi_n) :  n=1,2,\ldots\}$,
defined on a probability space $(\Omega^\sharp, \F^\sharp, \Prob^\sharp)$. Then
$\big\{ (X_n^\sharp, \xi_n^\sharp):  -\infty < n < \infty \big\}$ is a stationary, doubly-infinite Markov chain; and for each $k \in \Z$, 
\[
\Prob^\sharp\Big( \big(X_{k+n}, \xi_{k+n} \big)_{n \ge 0} \in \cdot \Big) ~=~ \Pstat \Big( \big(X_n, \xi_n\big)_{n \ge 0} \in \cdot \Big).
\]
Further define
$$ S_n^\sharp ~:=~ \begin{cases} \sum_{k=1}^n \xi_k^\sharp, & n > 0, \\ 0, & n=0, \\ - \sum_{k=n+1}^{0} \xi_k^\sharp, & n < 0.\end{cases}$$
Then for all $k \in \Z$, $S_k^\sharp - S_{k-1}^\sharp = \xi_k^\sharp$. Next introduce the ladder indices for $\big\{ (X_n^\sharp, S_n^\sharp) \big\}$, namely,
\begin{align*}
\zeta_0^\sharp ~:=&~ \sup\{n \le 0 \, : \, S_n^\sharp > \sup_{j < n} S_j^\sharp\} ;\\
\zeta_{k+1}^\sharp ~:=&~ \inf\{n > \zeta_k^\sharp \, : \, S_n^\sharp > S_{\zeta_k^\sharp}^\sharp \}, \qquad k=0,1,\ldots;
\end{align*}
and the ladder indices for $\big\{(X_n, S_n) \big\}$, namely,
\begin{equation} 
\zeta_0 ~:=~ 0 ; \qquad
\zeta_{k+1} ~:=~ \inf\{n > \zeta_k \, : \, S_n > S_{\zeta} \}, \qquad k=0,1,\ldots.
\end{equation}
In particular, $\zeta_1=\inf\{n > 0 \, : \, S_n > 0\}$.
Also define the measure $\psi$ on $\SPdp$ by setting
$$ \psi(A) ~=~ \Prob^\sharp\left( \zeta_0^\sharp = 0, X_0^\sharp \in A \right).$$
Then by Kesten's renewal theorem,
\begin{equation}\label{eq:markov delay distribution} 
\rho(E \times \Gamma) = \frac{1}{\lambda'(\alpha)} \int_\Gamma \, {\mathbb P}_{\psi}^\alpha \left( X_{\zeta_1} \in E, \, S_{\zeta_1} > s  \right) \,ds, \qquad E \subset {\mathbb S}^{d-1}_+,\:
\Gamma \in \reals_+;
\end{equation}
see Eqs. (1.16) and (3.10) in \shortciteN{Kesten1974}.


\vspace*{.3cm}

\noindent
{\bf Proof of Lemma \ref{lem:markovdelay}.} 
If the function $g$ takes both positive and negative values, then as the left- and right-hand sides of \eqref{eq:markovdelay} are finite 
(cf.\ \eqref{sec5-201}), we may split $g$ into its positive and negative parts, applying the result separately to each of these parts.
Thus, for the remainder of the proof, we will
assume without loss of generality  that $g$ is nonnegative.

\Step[1].  We start by applying Lemma 3 of \shortciteN{Kesten1974} to obtain an expression analogous to \eqref{eq:markovdelay},
but with respect to the positive ladder heights.   Set $h_t(x,s) = g(x,t-s)$ and
\[
{\cal R}(x,E,t) = \sum_{k=0}^\infty {\mathbb P}_x \Big( X_{\zeta_k} \in E, S_{\zeta_k} \in [0,t] \Big), \qquad x \in {\mathbb S}^{d-1}_+, \: E \in {\cal B}(S^{d-1}_+);
\]
cf.\ \shortciteN{Kesten1974}, Eq.\ (3.26).   Then it follows from these definitions that
\[
 \int_{y \in {\mathbb S}^{d-1}_+,\: 0 \le r \le t-s} h_t(y,t-r-s)  {\cal R}(x,dy,dr) 
~=~ \E_x^\alpha \left[ \sum_{k=0}^\infty g \left(X_{\zeta_k}, S_{\zeta_k}+s\right) {\bf 1}_{\{S_{\zeta_k}+s \in [0,t] \}} \right].
\]
Hence by Lemma 3 of \shortciteN{Kesten1974} and \eqref{eq:markov delay distribution}, we obtain upon letting $t \uparrow \infty$ that
\begin{align} \label{sec5-500}
\int_{{\mathbb S}^{d-1}_+ \times \reals_+}  \E_x^\alpha \left[ \sum_{k=0}^\infty g \left(X_{\zeta_k}, S_{\zeta_k}+s\right) {\bf 1}_{\{S_{\zeta_k}+s \ge 0\}}  \right]\, \rho(dx,ds) ~=&~ \frac{1}{\lambda^\prime(\alpha)}
  \int_{\SPdp} \psi(dx) \int_0^t  h_t(x,s) ds \nonumber\\[.05cm]
~=&~ \frac{1}{\lambda^\prime(\alpha)} \int_{\SPdp \times \reals_+} g(x,s) \psi(dx) ds.
\end{align}
 
\Step[2].  We would now like to extend \eqref{sec5-500} so that the sum on the left-hand side is taken over {\it all} $(X_n, S_n)$, rather than those values corresponding
 to the ladder heights.  Using a cyclic decomposition, observe that for any $(x,s) \in {\mathbb S}^{d-1}_+ \times \reals_+$,
 \begin{equation} \label{sec5-501}
 \E_x^\alpha \left[ \sum_{n=0}^\infty g(X_n, S_n +s) {\bf 1}_{\{S_n + s \ge 0\}} \right] ~=~ \E_x^\alpha \left[ \sum_{k=0}^\infty G(X_{\zeta_k}, S_{\zeta_k}+s) {\bf 1}_{\{S_{\zeta_k}+s \ge 0 \}} \right],
 \end{equation}
where 
$$ G(y,t) ~:=~ \E_y^\alpha \left[ \sum_{n=0}^{\zeta_1-1} g(X_n, S_n+t) {\bf 1}_{\{ S_n + t \ge 0 \}} \right], \quad t \in \reals.$$
Note that $g$ directly Riemann integrable implies that so is $G$.
Now apply \eqref{sec5-500} with $G$ in place of $g$ to obtain that
\begin{align} \label{sec5-502}
\int_{{\mathbb S}^{d-1}_+ \times \reals_+} & \E_x^\alpha \left[ \sum_{n=0}^\infty g \left(X_n, S_n+s\right) {\bf 1}_{\{ S_n + s \ge 0 \}}  \right]\, \rho(dx,ds) 
~=~ \frac{1}{\lambda^\prime(\alpha)} \int_{\SPdp \times \reals_+} G(x,s) \psi(dx) ds \nonumber\\[.2cm]
& \hspace*{3cm} =~  \frac{1}{\lambda'(\alpha)}\int_{\SPdp} \psi(dx) \:  \E_x^\alpha \left[ \sum_{n=0}^{\zeta_1-1} \int_{\reals_+} g(X_n, S_n+s) {\bf 1}_{\{ S_n + s \ge 0 \}}  ds \right] \nonumber\\[.2cm]
& \hspace*{3cm} =~  \frac{1}{\lambda'(\alpha)}\int_{\SPdp \times \reals_+} \psi(dx) \:  \E_x^\alpha \left[ \sum_{n=0}^{\zeta_1-1} g(X_n, s) \right] ds,
\end{align}
where the last step follows by observing that $S_n < 0$ prior to the first ascending ladder height, which implies that
$ \int_{\reals_+} g(X_n, S_n+s) {\bf 1}_{\{ S_n + s \ge 0 \}}  ds =  \int_{\reals_+} g(X_n, s) ds $.

\Step[3].  To establish the lemma, it remains to show that
\begin{equation} \label{sec5-503}
\int_{\SPdp} \psi(dx) \E_x^\alpha \left[ \sum_{n=0}^{\zeta_1 -1} g(X_n,s) \right] = \int_{\SPdp} g(x,s) \eta_\alpha(dx).
\end{equation}
Approximating both the positive and negative parts of $x \mapsto g(x,s)$ by simple functions, it is sufficient to show that
\begin{equation} \label{sec5-504}
\int_{\SPdp} \psi(dx) \E_x^\alpha \left[ \sum_{n=0}^{\zeta_1-1} {\bf 1}_{\{ X_n \in E \}} \right] = \eta_\alpha(E), \qquad E \in {\cal B}(\SPdp).
\end{equation}
To this end, observe using the definition of $\psi$ that
\begin{eqnarray} \label{sec5-505}
 \int_{\SPdp} \psi(dx) \E_x^\alpha \left[ \sum_{n=0}^{\zeta_1-1} {\bf 1}_{\{ X_n \in E \}} \right] 
&=& \sum_{l=0}^\infty \int_{\SPdp} \, \Prob^\sharp \left(  \sup_{j < 0} S_j^\sharp <0, \, X_0^\sharp \in dx \right) \, {\mathbb P}_x^\alpha ( X_l \in E, \, \zeta_1 > l) \nonumber\\[.2cm]
&=& \sum_{l=0}^\infty  \Prob^\sharp \left(  \sup_{j < 0} S_j^\sharp <0, \: \sup_{1 \le j \le l} S_j^\sharp \le 0, \: X_l^\sharp \in E\right) \nonumber\\[.2cm]
&=&  \sum_{l=0}^\infty  \Prob^\sharp \left(   \zeta_0^\sharp = -l , \: X_0^\sharp \in E\right) ~=~ \Prob^\sharp \left(  -\infty < \zeta_0^\sharp , \: X_0^\sharp \in E\right),
\end{eqnarray}
where in the last line, we have used the stationarity of the process $\{ (X_n^\sharp, \xi_n^\sharp):  -\infty < n < \infty \}$
 and the fact that $S_0^\sharp = 0$.
But $\Prob^\sharp  \big(  -\infty < \zeta_0^\sharp \big) =1$ (cf.\ \shortciteN{Kesten1974}, Lemma 2, Eq.\ (3.12)), and hence
\begin{equation} \label{sec5-506}
\Prob^\sharp \left(  -\infty < \zeta_0^\sharp ,\: X_0^\sharp \in E\right)  = \Prob^\sharp( X_0^\sharp \in E) = \Pstat(X_0 \in E) = \eta_\alpha(E).
\end{equation}
Then \eqref{sec5-504} is obtained by substituting \eqref{sec5-506} into \eqref{sec5-505}.
\halmos

We now establish Theorem \ref{thm2A} under the additional Hypothesis $(H_3)$ of Section 3.  This assumption will later be removed using a smoothing argument.

\begin{prop} \label{thm2A-underH3} 
Assume that Hypotheses $(H_1)$, $(H_2)$ and $(H_3)$ are satisfied, and suppose that ${\mathbb D} \in
{\cal B}(\realsd_+\setminus\{0\})$ is bounded and
$\pi({\mathbb D}) > 0$.   Then for any  $f \in {\cal C}_0\big({\R^d_+}\setminus\{0\}\big)$ ,
\begin{equation} \label{sec5-550}
 \lim_{u \to \infty} u^\alpha {\mathbb E} \left[ f \left(\frac{V}{u}\right) \right]  = \frac{C}{\lambda^\prime(\alpha)} 
 \int_{{\mathbb S}^{d-1}_+ \times \reals} e^{-\alpha s} f(e^s x) l_\alpha(dx) ds,
  \end{equation}
 where $C$ is given as in \eqref{thm2A-eq2}.  Equivalently, we have the weak convergence
 \begin{equation} \label{thm2A-eq3a-A}
\lim_{u \to \infty} u^\alpha {\mathbb P} \left( {\abs{V}}>tu, \frac{V}{\abs{V}} \in \cdot \right) ~\Rightarrow~ \frac{C}{\alpha \lambda^\prime(\alpha)} t^{-\alpha} l_\alpha(\cdot), \quad \mbox{for all } t > 0.
\end{equation} 
\end{prop}

\noindent
{\bf Proof.}  We first prove the result under the additional assumptions that $f$ is almost $\theta$-H\"{o}lder continuous (as defined in Eq. \eqref{def:almost hoelder}) for $\theta \le \min\{1, \alpha\}$. 

\Step[1].  First assume that $f(x)=\hat{f}(x){\mathbf 1}_{\{|x|\ge 1\}}$ for a $\theta$-H\"older continuous function $\hat{f}$ (i.e., $\delta=1$ in Eq. \eqref{def:almost hoelder}).  Since $(H_3)$ is satisfied, it follows from Lemma \ref{regn-cycle} that 
\begin{equation} \label{sec5-507}
{\mathbb E} \left[ {\hat f} \left(\frac{V}{u}\right) \1[{\{\abs{V}\ge u\}}] \right] = 
  \frac{1}{{\mathbb E_{\pi_\DD}}\left[ \tau \right]} \int_{{\mathbb D}} {\mathbb E} \left[ \sum_{i=0}^{\tau -1} {\hat f} \left( \frac{V_i}{u} \right){\bf 1}_{\{|V_i| \ge u \}}  \bigg| \, V_0=v\right] \, \pi_{\mathbb D}(dv),
\end{equation}
where $\tau$ denotes the first return time of $\{ V_n \}$ to ${\mathbb D}$.

Now apply Proposition \ref{prop:Kesten result} and \eqref{sec5-507} (separately to the positive and negative parts of ${\hat f}$) to obtain that
\begin{align} \label{sec5-508}
\lim_{u \to \infty} &u^\alpha {\mathbb E} \left[ {\hat f} \left(\frac{V}{u}\right) \1[{\{\abs{V}\ge u\}}] \right] \nonumber\\
 =& ~ \frac{1}{{\mathbb E_{\pi_\DD}} [\tau]} \int_{\mathbb D} \es[\alpha](\widetilde{v}) {\mathbb E}_{{\bfDV}}^{\alpha} \left[ |Z|^\alpha {\bf 1}_{\{ \tau = \infty \}} \right]
   \left( \int_{{\mathbb S}^{d-1}_+ \times \reals_+} {\mathbb E}_{x}^{\alpha} \left[  \sum_{i=0}^\infty F(X_i,S_i+s) \right] \, \rho(dx,ds) \right) \, \pi_{\mathbb D}(dv),
\end{align}  
where  $F(x,s) = \big( e^{-\alpha s} f(e^sx)/r_\alpha(x)\big) \1[{[0,\infty)}]$  (since ${\bf 1}_{\{ |x| \ge 1 \}}$ can be dropped in the last integral in \eqref{sec5-508}).
We remark that Proposition \ref{prop:Kesten result} is conditional on $V_0=v$.  To extend this result so that it holds conditional on $V_0 \sim \pi_{\mathbb D}$, 
we have applied a dominated convergence argument together with the bound provided by \eqref{eq:bound for thm2.1} of  Lemma \ref{sublemma2}.
Next observe by Lemma \ref{lem:markovdelay} that
\begin{align} \label{sec5-509}
&\int_{{\mathbb S}^{d-1}_+ \times \reals_+} {\mathbb E}_{x}^{\alpha} \left[  \sum_{i=0}^\infty F(X_i,S_i +s) \right] \, \rho(dx,ds) 
= \frac{1}{\lambda^\prime(\alpha)} \int_{\SPdp \times \reals_+} F(x,s) \eta_\alpha(dx) ds \nonumber\\[.2cm]
& \hspace*{2cm} = \frac{1}{\lambda^\prime(\alpha)} \int_{\SPdp \times \reals_+} \frac{e^{-\alpha s}}{r_\alpha(x)} f(e^s x) \, \eta_\alpha(dx) \, ds  
= \frac{1}{\lambda^\prime(\alpha)} \int_{\SPdp \times \reals_+} e^{-\alpha s}  f(e^s x) \, l_\alpha(dx) \, ds,
\end{align}
using that $\eta_\alpha(dx) = r_\alpha(x) l_\alpha(dx)$  (cf.\ \eqref{def:eta} and the discussion given there).   

Moreover, by Lemma \ref{lemma:hitting times Vn}, we have that $\pi(\DD)=(\E_{\pi_\DD}[\tau])^{-1}$, and recall 
 that $\pi_\DD(\cdot)=\pi(\cdot \cap \DD)/\pi(\DD)$.  Hence, by applying Lemma \ref{lem:convergence Z} (ii), we obtain that
\begin{align} \label{eq:calculating C}
\frac{1}{{\mathbb E_{\pi_\DD}} [\tau]} \int_{\mathbb D} \es[\alpha](\widetilde{v}) {\mathbb E}_{{\bfDV}}^{\alpha} \left[ |Z|^\alpha {\bf 1}_{\{ \tau = \infty \}} \right] \, \pi_{\mathbb D}(dv) ~=~ \int_{\mathbb D} \es[\alpha](\widetilde{v}) {\mathbb E}_{{\bfDV}}^{\alpha} \left[ |Z|^\alpha {\bf 1}_{\{ \tau = \infty \}} \right] \, \pi(dv) ~=~C,
\end{align}
where $C$ is given as in  \eqref{thm2A-eq2}.
Then \eqref{sec5-508}, \eqref{sec5-509}, and \eqref{eq:calculating C} imply that
\begin{equation}\label{regular variation outside 1}
\lim_{u \to \infty} u^\alpha {\mathbb E} \left[ f \left(\frac{V}{u}\right)  \right]  =   \lim_{u \to \infty} u^\alpha {\mathbb E} \left[ {\hat f} \left(\frac{V}{u}\right) \1[{\{\abs{V}\ge u\}}] \right]  = \frac{C}{\lambda^\prime(\alpha)} 
 \int_{{\mathbb S}^{d-1}_+ \times \reals_+} e^{-\alpha s} f(e^s x) l_\alpha(dx) ds
\end{equation}
for any bounded, almost $\theta$-H\"older continuous function $f$  of the form  $f(x)=\hat{f}(x){\mathbf 1}_{\{|x|\ge 1\}}$, where $\hat{f}$ is $\theta$-H\"{o}lder continuous.

\medskip

\Step[2].  Now suppose that $f$ is a bounded, almost $\theta$-H\"older continuous function with support in $\left( B_\delta(0)\right)^c$ for some $\delta >0$, and define $\check{f}(x)=\delta^{-\alpha}f(\delta v)$. Then $\check{f}$ is almost H\"older-continuous with support  in $\left(B_1(0) \right)^c$, and we infer from \eqref{regular variation outside 1} that 
\begin{align*}
   u^\alpha \E\left[ f\Big(\frac{V}{u}\Big) {\mathbf 1}_{\{|V| \ge  \delta u\}} \right]
  &~=~ u^\alpha \delta^\alpha \E \left[ \check{f}\Big( \frac{V}{\delta u} \Big) {\mathbf 1}_{\{|V| \ge \delta u \}} \right] 
  \stackrel{u \to \infty}{\longrightarrow}~ \frac{C}{\lambda^\prime(\alpha)} 
 \int_{{\mathbb S}^{d-1}_+ \times \reals_+} e^{-\alpha s} \check{f}(e^s x) l_\alpha(dx) ds \\
 &~=~  \frac{C}{\lambda^\prime(\alpha)} 
 \int_0^\infty \int_{{\mathbb S}^{d-1}_+}   e^{-\alpha s} \delta^{-\alpha} f(\delta e^s x) l_\alpha(dx) ds \\
  & ~ =~ \frac{C}{\lambda^\prime(\alpha)} 
\int_{\log \delta}^\infty \int_{{\mathbb S}^{d-1}_+}  e^{-\alpha r}  f( e^r x) l_\alpha(dx) dr, \qquad \mbox{\rm where} \:\: r = s + \log \delta.
\end{align*}
Since $f$ vanishes on $B_\delta(0)$, we may extend the outer integral to range from $0$ to $\infty$.
Thus we have obtained \eqref{sec5-550} for almost $\theta$-H\"older continuous functions on $\Rdo$. 

\medskip

\Step[3].  It remains to remove the assumption that $f$ is almost $\theta$-H\"{o}lder continuous, needed to apply Proposition \ref{prop:Kesten result} in the above argument.
To this end, observe that for all $r>0$,
$$ \Upsilon_u^{(r)} := u^\alpha \P{\frac{V}{u} \in \cdot, \frac{\abs{V}}{u} \ge r}$$
defines a family of uniformly bounded measures on $\Rd_+ \setminus B_r(0)$,
where the boundedness follows by employing \eqref{sec5-550} with $f(x) = \1[{\{|x|\ge r\}}]$, which is an almost $\theta$-H\"older continuous function. The Fourier characters $x \mapsto e^{i \skalar{x,y}}$ are bounded Lipschitz continuous functions for any $y \in \Rd$;  
then $f_y(x):=e^{i \skalar{x,y}}\1[{\{|x|\ge r\}}]$ is almost $\theta$-H\"older continuous for any $\theta \le \min\{1, \alpha\}$.  
 Let ${\mathfrak L}_\alpha$ be the measure on $\R^d_+\setminus\{0\}$ defined by the equation
\[
  \int_{{\mathbb S}^{d-1}_+ \times \reals} e^{-\alpha s} f(e^s x) l_\alpha(dx) ds =  \int_{{\mathbb R}^d_+\backslash \{0\}} f(x) {\mathfrak L}_\alpha(dx),
\]
and let ${\mathfrak L}_\alpha^{(r)}$ denote its restriction to a measure on $\Rd_+ \setminus B_r(0)$.
Then, based on what we have proved so far, we may infer the convergence, for all $y \in \Rd$ (considering real and imaginary part separately), of
$$ \lim_{u \to \infty} \int_{\Rd_+ \setminus B_r(0)} e^{i \skalar{x,y}} \Upsilon_u^{(r)}(dx) ~=~  \lim_{u \to \infty} u^\alpha {\mathbb E} \left[ e^{ i \skalar{u^{-1}V,y} } \1[{\{|V|\ge ru\}}] \right] ~=~ \frac{C}{\lambda'(\alpha)} \int_{\Rd_+ \setminus B_r(0)}  e^{i \skalar{x,y}} {\mathfrak L}_\alpha^{(r)}(dx);
 $$
and thus, by the L\'evy continuity theorem,  the weak convergence $ \Upsilon_u^{(r)} \Rightarrow \frac{C}{\lambda'(\alpha)}{\mathfrak L}_\alpha^{(r)}$, for any $r >0$. Now
if $f \in {\cal C}_0\big( \R_+^d \setminus\{0\}\big)$, then there exists $r>0$ such that $f$ is supported on $\big(B_r(0) \big)^c$. Hence
$$ \lim_{u \to \infty} u^\alpha \E\left[ f\Big(\frac{V}{u}\Big) \right] ~=~ \lim_{u \to \infty} \int f(x) \Upsilon_u^{(r)}(dx) ~=~ \frac{C}{\lambda'(\alpha)} \int f(x) {\mathfrak L}_\alpha^{(r)}(dx) ~=~ \frac{C}{\lambda'(\alpha)} \int f(x) {\mathfrak L}_\alpha(dx),$$ i.e.\ \eqref{sec5-550} holds.
 
Finally, the equivalence of \eqref{sec5-550} to \eqref{thm2A-eq3a-A} follows from Theorem 2 in \shortciteN{Resnick:2004}.
\halmos

\textbf{\textit{Smoothing.}}   To remove Hypothesis $(H_3)$, we employ a lower and upper approximation, where the approximating sequences are smoothed so
that $(H_3)$ is satisfied by these sequences.

We begin by constructing the lower approximating sequence. 
First recall the condition $({\mathfrak K})$ introduced just prior to the statement of Thereom \ref{thm2A}.  Also, from this
discussion in Section \ref{sect:2}, recall the definitions
\begin{equation} \label{MhatQhat}
\widehat{M}_n := M_{kn} \cdots M_{k(n-1) + 1} \quad \mbox{\text and} \quad  \widehat{Q}_n := \sum_{i=k(n-1)+1}^{kn} M_{kn} \cdots M_{i+1} Q_i,
\qquad n \in \pintegers.
\end{equation}
Now let $k \in \pintegers$ be chosen such that \eqref{requirement on k} holds.
Then $\{(\widehat{M}_n, \widehat{Q}_n)\}$ is an i.i.d.\ sequence under $\Prob$, and
with positive probability, $\widehat{Q}_1 - s {\vec 1} \succ 0$ for some $s > 0$.
Let ${\mathbb B}_n = \{ \widehat{Q}_n - s {\vec 1} \succ 0 \}$, and
let $\chi_{n,\epsilon}:=(-\epsilon) \chi_n$ for a sequence $\{\chi_n\}$ of i.i.d.\ random variables which are independent of $\{(\widehat{M}_n, \widehat{Q}_n)\}$  and have a nondegenerate absolutely continuous distribution
concentrated on $[0,1]^d$ (thus, $\chi_{n,\epsilon}$ is concentrated on $[-\epsilon,0]^d$).
Set $\widehat{Q}_{n,\epsilon} := \widehat{Q}_n + {\bf 1}_{{\mathbb B}_n} \chi_{n,\epsilon}$, and note that conditioned on the event ${\mathbb B}_n$, $\widehat{Q}_{n,\epsilon}$ has a 
continuous distribution function.
Since the event ${\mathbb B}_n$ occurs with positive probability, this implies that the distribution function of $\widehat{Q}_{n,\epsilon}$ has an absolutely continuous component
with respect to Lebesgue measure.

Now set 
\begin{equation} \label{sec6-51a}
V_{n,\epsilon} = \widehat{M}_n V_{n-1,\epsilon} + \widehat{Q}_{n,\epsilon}, \quad n=1,2,\ldots; \qquad V_{0,\epsilon} = V_0.
\end{equation}
Then $\{ V_{n,\epsilon} \}$ forms the {\it smoothed lower sequence}.
Let 
\begin{equation} \label{sec6-51aaa}
V_\epsilon := \widehat{Q}_{1,\epsilon} + \sum_{j=1}^\infty \widehat{M}_1 \cdots \widehat{M}_{j} \widehat{Q}_{j+1,\epsilon},
\end{equation}
and note that the law of $V_\epsilon$ is the stationary distribution of the process $\{ V_{n,\epsilon} \}$ with law $\pi_\epsilon$, say. 
 
 A {\it smoothed upper sequence} is constructed analogously, now choosing $\chi^\epsilon_n:=\epsilon \chi_n$, so that
 this random variable is concentrated on the interval $[0,\epsilon]^d$.
 Set $\widehat{Q}^\epsilon_n = \widehat{Q}_n + {\bf 1}_{{\mathbb B}_n} \chi^\epsilon_n$.  Then set
 $V_0^\epsilon = V_0$ and let
\[
V^{\epsilon}_n = \widehat{M}_n V^\epsilon_{n-1} + \widehat{Q}^{\epsilon}_n, \quad n=1,2,\ldots; \qquad
V^\epsilon = \widehat{Q}^{\epsilon}_1 + \sum_{j=1}^\infty \widehat{M}_1 \cdots \widehat{M}_{j} \widehat{Q}^{\epsilon}_{j+1}.
\]
Let  $\pi^\epsilon$ denote the distribution  of $V^\epsilon$.

\begin{rem} \label{drift.factor}
At this stage, it should be emphasized that this smoothing construction only affects the random quantity $\widehat{Q}_n$ and not $\widehat{M}_n$, and so the function $\Lambda$ is unchanged.
Thus, in particular, the solution $\alpha$ to the equation $\Lambda(\alpha) = 0$ and the corresponding invariant function $r_\alpha$ and  invariant measure $l_\alpha$ are the {\it same} as for the unsmoothed process. Moreover,
since $\widehat{M}_1=M_k \cdots M_1$ and $\lambda(\alpha)=1$, the factor $\lambda^\prime(\alpha)$ must now be replaced with $k \lambda^\prime(\alpha)$; cf.\ Lemma \ref{lem:slln}.
\end{rem}

\begin{rem}\label{remark: k skeleton}
Observe that if $k>1$ in \eqref{requirement on k}, then the evolution of the lower and upper smoothed sequence {\em cannot } be compared to the dynamics of the process $\{V_n\}$, but to that of the $k$-step chain $\{V_{kn} \, : \, n \in \N\}$, which at time $n$ is equal to
$$ \widehat{V}_n  := \widehat{M}_n \cdots \widehat{M}_1 V_0 + \sum_{j=1}^n \widehat{M}_n \cdots \widehat{M}_{j+1} \widehat{Q}_j. $$ We then have the sandwich inequality 
$$ V_{n, \epsilon} \le \widehat{V}_n   \le V_n^\epsilon, \quad \mbox{where} \quad \widehat{V}_n = V_{kn}.$$
{\it For the remainder of this section, we consider the $k$-step chain $\{\widehat{V}_n$\}, defined in terms of  $\{\widehat{M}_j, \widehat{Q}_j\}$.}  This $k$-step chain has the same stationary law, but different dynamics, than the 1-step chain $\{V_n\}$.
\end{rem}

For any $x \in \SPdp$ and $F \subset \SPdp$, let $d(x,F) = \inf \left\{ |x-y|:  y \in F \right\}$,
and for a given set $E \subset \SPdp$,  let
$$ E^\epsilon= \left\{x \in \SPdp \, : \, d(x,E) \le \frac{2 \epsilon}{s} \right\}\quad \mbox{\rm and} \quad E_\epsilon= \left\{x \in \SPdp \, : \, d(x,E^c) > \frac{2 \epsilon}{s} \right\}.$$

\begin{lemma}  \label{smoothingLm-1}
Let $\epsilon > 0$.  Then under the assumptions of Theorem \ref{thm2A}:

{\rm (i)}  The approximating sequences $\{ V_{n,\epsilon} \}_{n \in {\mathbb N}}$ and $\{ V_n^\epsilon \}_{n \in {\mathbb N}}$
each satisfy Hypothesis $(H_3)$.

{\rm (ii)}  We have the sandwich inequality
\begin{equation}\label{eq:sandwich P}
\P{ \abs{V_\epsilon} > tu, \ \frac{V_\epsilon}{\abs{V_\epsilon}} \in E_\epsilon}
~\le~ \P{ \abs{V} > tu, \ \frac{V}{\abs{V}} \in E} ~\le~ 
\P{ \abs{V^\epsilon} > tu, \ \frac{V^\epsilon}{\abs{V^\epsilon}} \in E^\epsilon}.
\end{equation}
\end{lemma}

\noindent
{\bf Proof.}
(i)   To verify part (i) of $(H_3)$, let $P_\epsilon$ denote the transition kernel of the process $\{ V_{n,\epsilon} \}$ in  \eqref{sec6-51a}.   Recall that $\chi_\epsilon$ is independent of $\widehat{M}$ and $\widehat{Q}$.  Hence, by construction, we have that
\begin{align*}
P_\epsilon(v,E) ~&=~\P{\widehat{M}v+\widehat{Q} \in E, \mathbb{B}^c} + \int_{[-\epsilon,0]^d} \P{\widehat{M}v+\widehat{Q} + y \in E, \mathbb{B}} \, \P{\chi_\epsilon \in dy}\nonumber\\[.2cm]
 ~&:=~P_{1,\epsilon}(v,E)+P_{2,\epsilon}(v,E).
\end{align*}
The kernel $P_{2,\epsilon}$ is obtained by the convolution of $\P{\widehat{M}v+\widehat{Q} \in \cdot,\, \mathbb{B}}$ with the probability measure $\P{\chi_\epsilon \in \cdot}$, which by 
assumption is smooth; thus $P_{2,\epsilon}(v,\cdot)$ itself has a Lebesgue density for all $v \in \R^d_+$.  Hence part (i) of $(H_3)$
is satisfied  with $\Phi$ taken to be Lebesgue measure and $F = \realsd_+$.

Since $\pi_\epsilon$ is the stationary distribution of the Markov chain with transition kernel $P_\epsilon$, it follows that $\pi_\epsilon$ also has a continuous component with respect to Lebesgue measure. Hence $(\supp \, \pi_\epsilon)^\circ \not= \emptyset$ and part (ii) of $(H_3)$ is satisfied. 

The verification for the process $\{ V_n^\epsilon \}$ is analogous.

(ii)  By construction,
\begin{equation} \label{differenceVepsilon}
V - V_{\epsilon} = - {\bf 1}_{{\mathbb B}_1} \chi_{1,\epsilon} - \sum_{k=1}^\infty \widehat{M}_1 \cdots \widehat{M}_k {\bf 1}_{{\mathbb B}_{k+1}} \chi_{k+1,\epsilon},
\end{equation}
since  $\widehat{Q}_k - \widehat{Q}_{k,\epsilon} = {\bf 1}_{{\mathbb B}_k} \chi_\epsilon$.
[Here we define
${\mathbb B}_k$ in the same way as ${\mathbb B}$, but with respect to the pair $(\widehat{M}_k,\widehat{Q}_k$).]
Consequently, setting $\widehat{M}_0$ to be equal to the identity matrix and recalling that $\chi_\epsilon$ is supported on $[-\epsilon, 0]^d$, we obtain that
\begin{equation}\label{eq:V-Veps1} 
\abs{V - V_\epsilon} ~\le~ \epsilon \left| \sum_{k=0}^\infty \left( \widehat{M}_0 \cdots \widehat{M}_k {\vec 1} \right) \1[{\mathbb B}_{k+1}] \right|.
\end{equation}
Moreover,
\begin{equation} \label{sec6-60}
|V| = \left| \sum_{k=0}^\infty \widehat{M}_0 \cdots \widehat{M}_{k} \widehat{Q}_{k+1} \right| \ge s \left| \sum_{k=0}^\infty  \left( \widehat{M}_0 \cdots \widehat{M}_{k}  {\vec 1} \right)  \1[{\mathbb B}_{k+1}]
\right|.
\end{equation}
This implies that
\begin{align}\label{eq:distVVeps}
\abs{ \frac{V}{\abs{V}} - \frac{V^\epsilon}{\abs{V^\epsilon}} } \le \frac{1}{|V|} \abs{V - V^\epsilon} + |V^\epsilon| \abs{\frac{1}{|V|} - \frac{1}{|V^\epsilon|}} \le 2 \frac{|V-V^\epsilon|}{|V|} \le 
\frac{2 \epsilon}{s}.
\end{align}
Hence
\[
\frac{V_\epsilon}{|V_\epsilon|} \in E_\epsilon \Longrightarrow d \left( \frac{V_\epsilon}{|V_\epsilon|}, E^c \right) > \frac{2\epsilon}{s}
  \Longrightarrow d \left( \frac{V}{|V|}, E^c \right) > 0  \Longrightarrow \frac{V}{|V|} \in E.
\]
Furthermore, by \eqref{differenceVepsilon}, we also have that $|V_\epsilon| > u \Longrightarrow |V| >u$.  Consequently,
\begin{equation} \label{sec6-61}
\P{ \abs{V_\epsilon} > u, \ \frac{V_\epsilon}{\abs{V_\epsilon}} \in E_\epsilon}
~\le~ \P{ \abs{V} > u, \ \frac{V}{\abs{V}} \in E},
\end{equation}
and the remaining inequality in \eqref{eq:sandwich P}
is established by an analogous argument.
\halmos

Since Hypothesis $(H_3)$ is satisfied for the two approximating sequences in Lemma \ref{smoothingLm-1}, it is natural to apply 
Proposition \ref{thm2A-underH3} to these sequences, yielding upper and lower bounds for  $\P{ \abs{V} > u, \, V/\abs{V} \in E}$ as $u \to \infty$.

Let $\big| \widehat Z_\epsilon \big|$ be defined in the same manner as the random variable $|Z|$ in Section \ref{sect:convergenceZ}, but with respect to the process
$\big\{ (\widehat{M}_i, \widehat{Q}_{i,\epsilon}):  i=1,2,\ldots \big\}$; namely,
\begin{equation} \label{sec6-590}
 \big| \widehat Z_\epsilon \big| ~=~|v| + \sum_{i=1}^{\infty} \frac{\skalar{\widehat{Y}_i,(\widehat{Q}_{i,\epsilon})^\sim}}{\skalar{\widehat{Y}_i, \widehat{X}_i}} \frac{\big|\widehat{Q}_{i,\epsilon} \big|}{\big| \widehat{M}_i 
   \cdots \widehat{M}_1 \widetilde{v}\big|}  \qquad {\mathbb P}_{\bfDV}^\alpha\text{\rm -a.s.},
 \end{equation}
where
\[
\widehat{Y}_i   := \lim_{n \to \infty} \left( \widehat{M}_i^{\,\top} \cdots \widehat{M}_n^{\, \top} {\vec 1} \right)^\sim, \quad n = 1,2,\ldots,
\]
and $\widehat{X}_i = \big(\widehat{M}_i 
   \cdots \widehat{M}_1 {v}\big)^\sim$.
Then, with $E = {\mathbb S}^{d-1}_+$, we obtain by Proposition \ref{thm2A-underH3}, Lemma  \ref{smoothingLm-1},
and Remark \ref{drift.factor} that 
\begin{equation}\label{eq:sandwich P-1}
\frac{C_\epsilon}{\alpha k \lambda^\prime(\alpha)} \le u^\alpha \liminf_{u \to \infty} \P{ \abs{V} > u}
~\le~ \limsup_{u \to \infty} u^\alpha \P{ \abs{V} > u} ~\le~ \frac{C^\epsilon}{\alpha k \lambda^\prime(\alpha)},
\end{equation}
where, in view of Lemma  \ref{lem:convergence Z} (ii), 
\[
C_\epsilon =  \int_{\mathbb D} r_{\alpha}(\widetilde{v}) {\mathbb E}^\alpha_{\bfDV} \left[  |\widehat Z_\epsilon|^\alpha {\bf 1}_{\{ \tau_\epsilon = \infty \}} \right] \pi_\epsilon(dv)
\quad \mbox{\rm and} \quad C^\epsilon = \int_{\mathbb D} r_\alpha(\widetilde{v})  {\mathbb E}^\alpha_{\bfDV}  \left[  |\widehat Z^\epsilon|^\alpha {\bf 1}_{\{ \tau^\epsilon = \infty \}} \right] \pi^\epsilon(dv).
\]
In what follows, we will generally write $\tau \equiv \tau({\mathbb D})$ to emphasize the dependence of this quantity on the choice of ${\mathbb D}$.
{\it However, it is important to observe from Proposition \ref{thm2A-underH3}  that $C_\epsilon$ and $C^\epsilon$ are universal constants, not dependent on the choice of ${\mathbb D}$.}

The next lemma shows that these constants converge to the required constant $C$ in \eqref{thm2A-eq2} as $\epsilon \downarrow 0$.

\begin{lemma}  \label{smoothingLm-2}
Assume the conditions of Theorem \ref{thm2A}.  Then for any set ${\mathbb D} = B_r(0) \cap \realsd_+$ with $\pi(\DD)>0$,
\begin{equation} \label{sec6-100}
C_\epsilon \nearrow C \quad \mbox{\rm and} \quad C^\epsilon \searrow C \quad \mbox{\rm as} \quad \epsilon \to 0,
\end{equation}
where $C$ is independent of the choice of $\DD$ and has the representation
$$ C= C({\mathbb D}) =  \int_{\mathbb D} r_{\alpha}(\widetilde{v}) {\mathbb E}^\alpha_{\bfDV} \left[  |\widehat Z|^\alpha {\bf 1}_{\{ \tau({\mathbb D}) = \infty \}} \right] \pi(dv) .$$
\end{lemma} 

\noindent
{\bf Proof.}  
We claim that
\begin{equation} \label{sec6-600}
\lim_{\epsilon \to 0} C_\epsilon \ge C(\bar{\mathbb D}) \qquad \mbox{and} \qquad \lim_{\epsilon \to 0} C^\epsilon \le C({\mathbb D}).
\end{equation}
Note the these limits necessarily exist, since $C_\epsilon$, $-C^\epsilon$ are monotonically increasing (this follows from the
 monotonicity, in $\epsilon$, of $V_\epsilon$ and $V^\epsilon$ and the representation \eqref{thm2A-eq3a-A}).

\Step[1]. We begin by establishing the first inequality in \eqref{sec6-600}.  Set
\begin{equation}\label{def:Hepsilon} H_\epsilon(v) =  r_\alpha(v)  {\mathbb E}^\alpha_{\bfDV} \left[  |\widehat Z_\epsilon|^\alpha {\bf 1}_{\{ \tau_\epsilon(\bar{\mathbb D}) = \infty \}} \right] \quad \text{ and } \quad H(v) = r_\alpha(v)  {\mathbb E}^\alpha_{\bfDV} \left[  |\widehat Z|^\alpha {\bf 1}_{\{ \tau(\bar{\mathbb D}) = \infty \}} \right] .
\end{equation}
Then we need to show that
\begin{equation}\label{eg:Hepsilon to H} \liminf_{\epsilon \to 0} \int_{\bar \DD} H_\epsilon(v) \pi_\epsilon(dv) ~\ge~ \int_{\bar{\DD}} H(v) \pi(dv).
\end{equation}

We will prove below that: (i) $H_\epsilon(v) \uparrow H(v)$ as $\epsilon \to 0$; and (ii) for all $\epsilon \ge 0$, the function $v \mapsto H_\epsilon(v)$ is lower semicontinuous.

Assume that (i) and (ii) hold, and fix $\epsilon_0>0$. Since $H_\epsilon(v)$ is a monotone increasing sequence as $\epsilon \downarrow 0$, we have that
\[
C_\epsilon \ge \int_{\bar{\DD}} H_{\epsilon_0}(v) \pi_\epsilon(dv) \qquad \mbox{for all } \epsilon \le \epsilon_0.
\]
 As the function $v \mapsto H_{\epsilon_0}(v)$ is lower semicontinous and bounded from below by 0,
and $\pi_\epsilon \Rightarrow \pi$ (cf.\ \eqref{differenceVepsilon}),
the Portmanteau theorem (\shortciteN{AVJW96}, Theorem 1.3.4 (iv)) yields that
$$ \liminf_{\epsilon \to 0} C_\epsilon ~\ge~ \liminf_{\epsilon \to 0} \int_{\bar \DD} H_{\epsilon_0}(v) \pi_\epsilon(dv) ~\ge~ \int_{\bar \DD} H_{\epsilon_0}(v) \pi(dv).$$
Now we let $\epsilon_0 \to 0$ and use the monotone convergence $H_{\epsilon_0} \uparrow H$ to infer, using the monotone convergence theorem,
$$ \liminf_{\epsilon \to 0} C_\epsilon ~\ge~ \lim_{\epsilon_0 \to 0} \int_{\bar \DD} H_{\epsilon_0}(v) \pi(dv) ~=~ \int_{\bar \DD} H(v) \pi(dv) ~=~C({\bar \DD}). $$

It remains to prove (i) and (ii). In order to obtain (i), observe that $|V_{n,\epsilon}|$ increases monotonically to $|V_n|$ as $\epsilon \to 0$.  Thus, if the process $\{ V_n \}$ enters $\bar{\mathbb D}$, then so does $\{ V_{n,\epsilon} \}$
for all $\epsilon > 0$.  Hence, we trivially obtain that ${\bf 1}_{\{ \tau(\bar{\mathbb D}) = \infty \}} \ge {\bf 1}_{\{\tau_\epsilon(\bar{\mathbb D}) = \infty\}}$, where
 $\tau(\bar{\mathbb D})$, $\tau_\epsilon(\bar{\mathbb D})$ are the first passage times of $\{ V_n\}$,
$\{ V_{n,\epsilon} \}$ into $\bar{\mathbb D}$, respectively.  Conversely, observe that if $\tau(\bar{\mathbb D}) = \infty$, then
${\bf V}:= (V_1,V_2,\ldots) \in \left(\bar{\mathbb D}^c \right)^{\N}$, which is open.  Now ${\bf V}_\epsilon:= (V_{1,\epsilon},V_{2,\epsilon},\ldots)$ converges to ${\bf V}$ a.s.\
in the product topology (as $\chi_\epsilon$ is supported on $[-\epsilon,0]^d$).  It follows that
${\bf V}_\epsilon \in \left(\bar{\mathbb D}^c \right)^{\N}$ for sufficiently small $\epsilon$.  Consequently, ${\bf 1}_{\{\tau(\bar{\mathbb D}) = \infty\}} \le \liminf_{\epsilon \to 0} 
{\bf 1}_{\{ \tau_\epsilon(\bar{\mathbb D}) = \infty \}}$.
Thus we conclude that ${\bf 1}_{\{\tau(\bar{\mathbb D}) = \infty\}} = \lim_{\epsilon \to 0} {\bf 1}_{\{ \tau_\epsilon(\bar{\mathbb D}) = \infty\}}$, and moreover, the convergence is monotone, i.e. ${\bf 1}_{\{ \tau_\epsilon(\bar{\mathbb D}) = \infty\}} \uparrow {\bf 1}_{\{\tau(\bar{\mathbb D}) = \infty\}}$ as $\epsilon \to 0$.
Furthermore, 
 as $\widehat{Q}_\epsilon$ increases componentwise  to $\widehat{Q}$,
we deduce from \eqref{sec6-590} and
Lemma \ref{lem:convergence Z} (ii) that $|\widehat Z_\epsilon| \uparrow | \widehat Z |$ as $\epsilon \to 0$.
By Lemma \ref{lem:convergence Z} (iii), $|\widehat Z|^\alpha{\bf 1}_{\{\tau(\bar{\mathbb D}) = \infty\}}$ is an integrable upper bound for the family $\big\{ |\widehat Z_\epsilon|^\alpha{\bf 1}_{\{\tau_\epsilon(\bar{\mathbb D}) = \infty\}} \big\}_{\epsilon >0}$, and thus we obtain, for all $v \in \bar{\DD}$, the monotone convergence $H_\epsilon(v) \uparrow H(v)$.

To obtain (ii), observe that if $v \to \hat{v}$, then ${\bf V}(v)$ converges to ${\bf V}(\hat{v})$, where ${\bf V}:= (V_{1,\epsilon},V_{2,\epsilon},\ldots)$
and ${\bf V}(v)$ denotes the dependence of this quantity on its initial state.  Then by repeating the argument given above, we obtain that $\bar{\mathbb D}$ closed 
$\Longrightarrow {\bf 1}_{\{\tau(\bar{\mathbb D},\hat{v}) = \infty\}} \le \liminf_{v \to \hat{v}} 
{\bf 1}_{\{ \tau(\bar{\mathbb D},\hat{v}) = \infty \}}$, where $\tau(\bar{\mathbb D}, \cdot)$ 
again denotes the dependence on the initial state.  From the representation \eqref{sec6-590}, we obtain that $v \mapsto \widehat Z_\epsilon(v)$ is continuous a.s. (the series converges a.s. by Lemma \ref{lem:convergence Z} (i)). Then we can apply Fatou's Lemma and use the continuity of $r_\alpha$ to infer that $H_\epsilon$ is lower semicontinuous.

\Step[2]. To establish the second inequality in \eqref{sec6-600}, we proceed as before, now using 
$${\rm (i^\prime)}: \quad H^\epsilon(v) :=  r_\alpha(v)  {\mathbb E}^\alpha_{\bfDV} \left[  |\widehat Z^\epsilon|^\alpha {\bf 1}_{\{ \tau^\epsilon({\mathbb D}) = \infty \}} \right] ~\downarrow~  r_\alpha(v)  {\mathbb E}^\alpha_{\bfDV} \left[  |\widehat Z|^\alpha {\bf 1}_{\{ \tau({\mathbb D}) = \infty \}} \right] =: H^\circ(v) $$
and (ii$^\prime$) the upper semicontinuity of $H^\epsilon(v)$, which follows since we consider now the hitting time of an open set. Furthermore, Lemma \ref{lem:convergence Z} (iii) gives that $\sup_{v \in \DD} H_{\epsilon_0}(v) \le B$ for some finite bound $B$. Then we can apply the Portmanteau theorem (\shortciteN{AVJW96}, Theorem 1.3.4 (v)) to infer that
$$ \limsup_{\epsilon \to 0} C^\epsilon ~\le~ \limsup_{\epsilon \to 0} \int_{ \DD} H^{\epsilon_0}(v) \pi^\epsilon(dv) ~=~ \int_{ \DD} H^{\epsilon_0}(v) \pi(dv)$$
for all $\epsilon_0 >0$, and thus, letting $\epsilon_0 \to 0$ and using (i$^\prime$),
$$ \limsup_{\epsilon \to 0} C^\epsilon ~\le~ \int_{\DD} H^\circ(v) \pi(dv) ~=~ C(\DD).$$

\Step[3].  Having obtained \eqref{sec6-600},
it remains to show that if ${\mathbb D} = B_r(0) \cap \realsd_+$, where $\pi(\DD)>0$, then, in fact,
\begin{equation} \label{sec6-610}
\lim_{\epsilon \to 0} C_\epsilon \ge C({\mathbb D}).
\end{equation} [This implies immediately that $\lim_{\epsilon \to 0} C_\epsilon = C({\mathbb D})$ and also that
 $\lim_{\epsilon \to 0} C^\epsilon =C(\DD)$, since $C^\epsilon \ge C_\epsilon$ and $\lim_{\epsilon \to 0} C^\epsilon \le C(\DD)$ by  \eqref{sec6-600}.]

To this end, let $\{ r_i \}$ be chosen such that $r_i \uparrow r$ as $i \to \infty$ and set ${\mathbb D}_i = B_{r_i}(0) \cap \realsd_+$.   If $\{ V_n \}$ avoids 
 ${\mathbb D}$, then it also avoids each ${\mathbb D}_i$, so we trivially obtain that ${\bf 1}_{\{ \tau(\bar{\mathbb D}_i) = \infty \}} \ge {\bf 1}_{\{\tau({\mathbb D}) = \infty\}}$.  Conversely,
$V_n \in {\mathbb D} \Longrightarrow V_n \in \bar{\mathbb D}_i$ for sufficiently large $i$.  Thus $\lim_{i \to \infty} {\bf 1}_{\{ \tau(\bar{\mathbb D}_i) = \infty\}} = {\bf 1}_{\{\tau({\mathbb D})
= \infty\}}$.  Now $\lim_{\epsilon \to 0} C_\epsilon$ is a universal constant, independent of the choice of the set $\bar{\mathbb D}$ in \eqref{sec6-600}. 
Consequently,  we conclude by \eqref{sec6-600} that
\begin{eqnarray} \label{convergenceDi}
\lim_{\epsilon \to 0} C_\epsilon &\ge& \lim_{i \to \infty} C(\bar{\mathbb D}_i) 
~=~ \lim_{i \to \infty} \int_{{\bar \DD_i}} r_{\alpha}(\widetilde{v}) {\mathbb E}^\alpha_{\bfDV} \left[  |\widehat{Z}|^\alpha {\bf 1}_{\{ \tau(\bar{\mathbb D}_i) = \infty \}} \right] \pi(dv) \nonumber\\
&=&  \int_{\mathbb D} r_{\alpha}(\widetilde{v}) {\mathbb E}^\alpha_{\bfDV} \left[  |\widehat Z|^\alpha {\bf 1}_{\{ \tau({\mathbb D}) = \infty \}} \right] \pi(dv) = C({\mathbb D}),
\end{eqnarray}
as required. 
\halmos

\noindent
{\bf Proof of Theorem \ref{thm2A}.}  It follows directly from Proposition \ref{thm2A-underH3} and Lemmas \ref{smoothingLm-1}
and \ref{smoothingLm-2} that for any $E \in {\cal B}({\mathbb S}^{d-1}_+)$,
\begin{equation} \label{sec6-650}
\liminf_{u \to \infty} u^\alpha {\mathbb P} \left( \abs{V} > tu, \: \frac{V}{\abs{V}} \in E \right) \ge \frac{C}{\alpha k \lambda^\prime(\alpha)} t^{-\alpha} \,  \limsup_{\epsilon \to 0} l_\alpha(E_\epsilon)
\end{equation}
and
\begin{equation} \label{sec6-651}
\limsup_{u \to \infty} u^\alpha {\mathbb P} \left( \abs{V} > tu, \: \frac{V}{\abs{V}} \in E \right) \le \frac{C}{\alpha k \lambda^\prime(\alpha)} t^{-\alpha} \, \liminf_{\epsilon \to 0} l_\alpha(E^\epsilon).
\end{equation}
Now if $\l_\alpha(\partial E)=0$, then
$$ \limsup_{\epsilon \to 0} l_\alpha(E_\epsilon) = \liminf_{\epsilon \to 0} l_\alpha(E^\epsilon) = l_\alpha(E).$$
Hence the bounds coincide and, thus, for all measurable $E \subset \SPdp$ with $l_\alpha(\partial E)=0$,
$$ \lim_{u \to \infty} u^\alpha\P{|V| > tu, \ \frac{V}{|V|} \in E} ~=~ \frac{C}{\alpha k \lambda'(\alpha)} t^{-\alpha} l_\alpha(E).$$
By the Portmanteau Theorem, this implies the weak convergence
\begin{equation}
 u^\alpha {\mathbb P} \left( |V| > tu, \, \frac{V}{\abs{V}} \in \cdot \right) ~\Rightarrow~ \frac{C}{\alpha \lambda'(\alpha)} t^{-\alpha} l_\alpha(\cdot) \qquad \text{ as } u \to \infty,
\end{equation}
for all $t > 0$,
which is equivalent to \eqref{thm2A-eq1} by Theorem 2 of \shortciteN{Resnick:2004}.
\halmos

%
%
%
%

\section{Proof of Theorem \ref{thm2B}}\label{sect:6}

We now turn to the proof of Theorem \ref{thm2B}.   Note that if $A$ is a semi-cone, then the
 event $\{ V_n \in uA, \:\: \mbox{some } n \le N \}$ corresponds to the event
$\{ \max_{1 \le n \le N} |V_n^A| > u \}$, where $V_n^A := |V_n|/d_A(\tilde{V}_n)$ for $d_A(x) := \inf\{ t>0:  tx \in A\}$.
Thus, as a first step, we study maxima of $\{|V_n^A|\}$ over cycles emanating from a set $\DD$, where 
the cycle is terminated upon the return of the process to $\DD$.  Consequently, 
we obtain the asymptotic distribution of $T_u^A/u^\alpha$, first when $d_A$ is assumed to be bounded, and then for general sets $A$.
Throughout this section, we assume that the set $A$ is a semi-cone.

We begin by establishing a preliminary lemma, where
we identify the constant appearing in the ruin problem for the Markov random
walk $\{ (X_n,S_n^A):  n=0,1,\ldots \}$.  For this purpose, define
\[
{\mathfrak T}^A_u = \inf \Big\{ n \in {\mathbb N}:  M_n \cdots M_1 \widetilde{V}_0 \in uA \Big\}, \quad \mbox{where} \quad V_0 = v \in \realsd_+ \setminus \{ 0 \}.
\]

\begin{lemma} \label{ruin.est}
Suppose that Hypotheses $(H_1)$ and $(H_2)$ are satisfied, and assume that $d_A$  is bounded and continuous on $\SPdp$. Then 
\begin{equation} \label{sec7.r.1}
\lim_{u \to \infty} u^\alpha {\mathbb P} \left( \left. {\mathfrak T}_u^A < \infty \right| V_0 = v \right) =  r_\alpha(\widetilde{v})  \int_{{\mathbb S}^{d-1}_+ \times \reals_+} \: 
\frac{e^{-\alpha s}}{r^A_\alpha(x)} \rho^A(dx, ds) ~:=~ r_\alpha(\widetilde{v}) D_A,
\end{equation}
where $\rho^A$ is given as in Theorem \ref{prop:asymptotic overjump distribution}.
\end{lemma}

\noindent
{\bf Proof.}   Converting to the $\alpha$-shifted measure, we obtain 
for any $v \in \realsd_+ \setminus \{ 0 \}$ that
\begin{eqnarray} \label{sec7.r.2}
u^{ \alpha}{\mathbb P} \big( {\mathfrak T}_u^A < \infty \, | \, V_0 = v \big) & = &  r_\alpha(\widetilde{v}){\mathbb E}^\alpha_{\tilde{v}} \left[ \Big(e^{-\alpha \big(S_{{\mathfrak T}_u^A} -\log u\big)}\big/r_\alpha\big(X_{{\mathfrak T}_u^A}\big) \Big) {\bf 1}_{\{ {\mathfrak T}_u^A < \infty\}} \right]\nonumber\\[.3cm]
& = &  r_\alpha(\widetilde{v}) {\mathbb E}^\alpha_{\tilde{v}} \left[ \Big(e^{-\alpha \big(S^A_{{\mathfrak T}_u^A} -\log u\big)}\big/r_\alpha^A\big(X_{{\mathfrak T}_u^A}\big)\Big) {\bf 1}_{\{ {\mathfrak T}_u^A < \infty\}} \right],
\end{eqnarray}
where the last step follows from the definitions of $S_n^A$ and $r_\alpha^A$.  
To characterize the limit on the right-hand side, apply Kesten's renewal theorem (Theorem \ref{prop:asymptotic overjump distribution}) for the bounded continuous function $g(x,s)= e^{-\alpha s}/r^A_\alpha(x)$.  This yields \eqref{sec7.r.1}. 
\halmos

To establish the weak convergence of $\{ T_u^A \}$, the main idea will be to study the excursions of $\{ V_n \}$ over cycles eminating from the set ${\mathbb D}$. For this purpose, we introduce the random variables
\[
U_i := \max_{\kappa_{i-1} < n \le \kappa_i} V_n^A, \quad i=1,2,\ldots,
\]
where $V_n^A := |V_n|/ d_A(\widetilde{V}_n)$,
and where $\kappa_0=0$ and $\kappa_i=\inf\{n > \kappa_{i-1} \, : \, V_n \in \DD\}$ denote the successive return times to $\DD$.  Also set
\begin{align*}
{\cal M}_n^U &= \max \left\{ U_1, \ldots, U_n \right\}, \quad n =1,2,\ldots;\\[.2cm]
{\cal M}_n &= \max\left\{ V_1^A,\ldots, V_n^A \right\}, \quad n=1,2,\ldots.
\end{align*}
Recall that  $ \{ T_u^A \le N \} = \{ V_n^A > u, \mbox{ some } n \le N\}$.  
Thus, $\{ {\cal M}_n^U  > u \}$ describes the event that $T_u^A$ occurs by the random time $\kappa_n$, 
while $\{ {\cal M}_n > u \}$ describes the event that $T_u^A$ occurs by the deterministic time $n$.

\begin{prop} \label{basicprop}
Suppose that Hypotheses $(H_1)$ and $(H_2)$ are satisfied and there exists $m \in \pintegers$ such that $(H_3)$ holds for the  $m$-skeleton $\{V_{mn} \, : \, n \in \N \}$.  Assume that ${\mathbb D} \in
{\cal B}(\realsd_+)$ is bounded and
$\pi({\mathbb D}) > 0$, and 
suppose that the function $d_A$ is bounded and continuous on $\SPdp$.  Let
$\gamma_0$ be an arbitrary probability distribution on $\Rd_+\setminus\{0\}$.
Then 
\begin{equation} \label{sec7.prop.1}
\lim_{n \to \infty} \P{{\cal M}_n^U \le n^{1/\alpha} u \, \big|  \, V_0 \sim \gamma_0}~=~  \exp\big\{- K_A \, {\mathbb E}_{\pi_{\mathbb D}}[\tau] \, u^{-\alpha} \big\},
\end{equation}
where $K_A = C D_A$ and $C$ is given as in \eqref{thm2A-eq2}.
\end{prop}

Unless explicitly noted, we assume throughout the rest of this section that $V_0 \sim \gamma_0$, i.e., $\Prob=\Prob_{\gamma_0}$.

\medskip

\noindent
{\bf Proof.}  Set $u_n = n^{1/\alpha} u$. Then for any $l \in \pintegers$,
\begin{equation} \label{sec7.prop.3}
\sum_{i=1}^l \P{U_i > u_n} - \sum_{1 \le i < j \le l} \P{U_i > u_n, \, U_j>u_n} \le \P{ {\cal M}_l^U > u_n} \le \sum_{i=1}^l \P{U_i > u_n}.
\end{equation}

We begin by calculating $\sum_{i=1}^{l(n)} \P{U_i > u_n} $ as $n \to \infty$ for the sequence $l(n)=\lfloor n/k \rfloor$ and  fixed $k
\in \pintegers$.  By Corollary \ref{cor4-1} and the Markov property, we have for all $i \in \pintegers$ and $v \in \DD$ that
\begin{align} \label{sec7.prop.4}
\lim_{n \to \infty}  n u^\alpha  {\mathbb P}& \left( U_i > u_n \left| V_{\kappa_{i-1}} = v \right. \right)
= \lim_{n \to \infty} nu^\alpha \P{T_u^A < \tau \, | \, V_0=v}    \nonumber \\
 &= r_\alpha (\widetilde{v})  {\mathbb E}_{\bfDV}^{\alpha} \left[ |Z|^\alpha {\bf 1}_{\{ \tau = \infty \}} \right] \int_{{\mathbb S}^{d-1}_+ \times \reals_+} \: \frac{e^{-\alpha s}}{r^A_\alpha(x)} \rho^A(dx,ds)
 := H(v) \quad \Big(=C(v) D_A\Big).
\end{align}
[Note that this equation also holds if $i=0$ and $v \in \realsd_+ \setminus \{0\}$, as $V_{\kappa_0}$ need not belong to $\DD$.]

Under Hypotheses ($H_3$), $\{V_{mn} \, : \, n \ge 0\}$ is a positive aperiodic Harris chain (Lemma \ref{georecurrent}).
Hence $\{V_n\}$ is a positive, $m$-periodic Harris chain. Then the hitting chain $\{V_{\kappa_i}\}$ is itself
a positive $m$-periodic Harris chain as well (\shortciteN{Alsmeyer1991}, Theorem 8.3.7), 
and the invariant measure of this chain is $\pi_\DD$ (Lemma \ref{lemma:hitting times Vn}). 
If $\gamma_i$ denotes the law of $V_{\kappa_i}$, $i \in \pintegers$, then Harris recurrence gives that $| n^{-1} \sum_{i=1}^n\gamma_i-\pi_{\DD}|_{\text{TV}} \to 0$ as $n \to \infty$, where $|\cdot|_{\text{TV}}$ denotes the total variation distance; see Theorem 13.3.4 in \shortciteN{SMRT93}.
Set
\[
H_n(v) = n u^\alpha \P{ \left. U_1 > u_n \, \right| \, V_0 = v}.
\] 
By  Lemma \ref{sublemma2}
(specifically,  \eqref{eq:bound for thm2.1} with
$h = 1$),
 we have that $\sup \left\{ H_n(v):  v \in \DD,\: n \in \N \right\} \le B$, for some finite constant $B$.
Then $ n u^\alpha \P{U_i > u_n} =\int_{\DD} H_n(v) \gamma_{i-1}(dv)$, and 
\begin{align} \label{sec7.prop.5}
 \bigg| \frac{k}{n} \sum_{i=1}^{\lfloor n/k \rfloor} & nu^\alpha \P{U_i > u_n} - \int_{\DD} H(v) \pi_{\DD}(dv) \bigg|\nonumber\\ 
 \le~& 
\abs{\int_{\DD} H_n(v)  \bigg(\frac{k}{n} \sum_{i=1}^{\lfloor n/k \rfloor} \gamma_{i-1} - \pi_{\DD}\bigg)(dv)  - \int_{\DD} \big(H_n(v)-H(v)\big) \pi_{\DD}(dv)} \nonumber
\\ 
\le~& B \abs{\frac{k}{n} \sum_{i=1}^{\lfloor n/k \rfloor}\gamma_{i-1}- \pi_{\DD}}_{\text{TV}} + \int_{\DD} \abs{H_n(v)-H(v)} \pi_{\DD}(dv).
\end{align}
The second term tends to zero as $n \to \infty$ by dominated convergence and the fact that $H_n(v) \to H(v)$, by \eqref{sec7.prop.4}. 
Thus the left-hand side of \eqref{sec7.prop.5} tends to zero as $n \to \infty$ and hence,  using
 \eqref{sec7.prop.4},
\begin{equation} \label{sec7.prop.6}
\lim_{n \to \infty} \sum_{i=1}^{\lfloor n/k \rfloor} u^\alpha \P{U_i > u_n} = \frac1k \int_{\mathbb D} H(v) \pi_{\mathbb D}(dv) 
= \frac{D_A}{k} 
\int_{\DD} C(v) \, \frac{\pi(dv)}{\pi(\DD)} 
=  \frac{K_A {\mathbb E}_{\pi_{\mathbb D}}[\tau]}{k},  
 \end{equation}
since $C=\int_{\DD} C(v) \pi(dv)$ and ${\mathbb E}_{\pi_{\mathbb D}}[\tau] = \left(\pi({\mathbb D}) \right)^{-1}$, by Lemma \ref{lemma:hitting times Vn}.
Using this equation in \eqref{sec7.prop.3}, we deduce that for any $k \in \pintegers$,
\begin{equation} \label{sec7.prop.8}
\limsup_{n \to \infty}  \P{{\cal M}^U_{\lfloor n/k \rfloor} > u_n} \le \frac{K_A \E_{\pi_{\mathbb D}}{[\tau]}}{k} u^{-\alpha}.
\end{equation}
Note that the right-hand side is {\it independent} of $V_0 \in \Rd_+\setminus\{0\}$.  Since this quantity is asymptotically independent of the initial state, the same calculation yields the asymptotic
behavior of the maximum over any block of length $n/k$; more precisely, for $\limsup_{n \to \infty} {\mathbb P} \left( \left. U_{\lfloor jn/k \rfloor +1},\ldots, U_{\lfloor (j+1)n/k \rfloor} > u_n \, \right| \,
{\cal F}_{\kappa_{\lfloor jn/k \rfloor}} \right)$, $j=0,\ldots,k-1$.
Hence we conclude from \eqref{sec7.prop.8} that
\begin{equation} \label{sec7.prop.9}
\limsup_{n \to \infty} {\mathbb P} \left( {\cal M}^U_n \le u_n \right) \le \left( 1 - \frac{K_A \E_{\pi_{\mathbb D}}[\tau] u^{-\alpha}}{k} \right)^k \to \exp \Big\{-K_A {\mathbb E}_{\pi_{\mathbb D}}[\tau] u^{-\alpha} \Big\} \quad \mbox{as} \quad k \to \infty.
\end{equation}
Moreover,
using once again the uniform upper bound (in the initial state) provided by Lemma \ref{sublemma2}, we obtain that 
for any positive integer $k$,
\begin{equation} \label{sec7.prop.10}
\limsup_{n \to \infty}  \sum_{1 \le i < j \le \lfloor n/k \rfloor} \P{ U_i > u_n, \, U_j > u_n } = o\left( \frac{1}{k} \right) \quad \mbox{as} \quad n \to \infty.
\end{equation}
Finally, using \eqref{sec7.prop.6} and \eqref{sec7.prop.10} in \eqref{sec7.prop.3},  we conclude that
\begin{equation*} 
\hspace*{4.8cm}
\liminf_{n \to \infty} {\mathbb P} \left( {\cal M}^U_n \le u_n \right) \ge \exp \Big\{-K_A {\mathbb E}_{\pi_{\mathbb D}}[\tau] u^{-\alpha} \Big\}. 
\hspace*{4.3cm} \halmos
\end{equation*}

\begin{lemma} \label{cont-in-prob}
Suppose that Hypotheses $(H_1)$--$(H_3)$ are satisfied and 
 the function $d_A$ is bounded and continuous.  Then
for any $\Delta > 0$, there exists a constant $\delta > 0$ such that
\begin{equation}
\limsup_{n \to \infty} {\mathbb P}\left( \max_{|m-n| <  n\delta} \left| {\cal M}_m - {\cal M}_n \right| > n^{1/\alpha} \Delta \right) \le\Delta.
\end{equation}
\end{lemma} 

\noindent
{\bf Proof.}    
Let $k \in \pintegers$.  Then 
\[
{\cal M}_{n+k} = \max \left\{ {\cal M}_n, V_{n+1}^A, \ldots, V_{n+k}^A \right\} \Longrightarrow {\cal M}_n \le {\cal M}_{n+k} \le {\cal M}_n + \max\left\{ V_{n+1}^A,\ldots, V_{n+k}^A \right\}.
\]
Since  $V_n^A := |V_n|/ d_A(\widetilde{V}_n)$, it follows that
\begin{equation} \label{sec7.prop.20}
\max_{|m-n| <  n\delta} \left| {\cal M}_m - {\cal M}_n \right| \le b \max\left\{ |V_{\lfloor n - n\delta \rfloor + 1}|, \ldots, |V_{\lfloor n + n \delta \rfloor}| \right\},
\end{equation}
where $b = \max \big\{ \left(d_A(x) \right)^{-1} \, : \, {x \in \SPdp} \big\} < \infty$.
We now determine the maximum on the right-hand side conditioned on $V_{\lfloor n - n\delta \rfloor} =v$.  Equivalently, we study ${\cal M}_{m_n}$
conditioned on $V_0 =v$, where $m_n = \lfloor n + n \delta \rfloor  - \left( \lfloor n - n \delta \rfloor + 1 \right)  \le 2 n \delta.$

Let $\DD \subset \realsd_+ \setminus \{ 0 \} $ be chosen such that  $\pi({\mathbb D}^c) \le \Delta/2$, and let $v \in \DD$.
Since ${\cal M}_n \le {\cal M}_n^{U}$  and $m_n \le 2 n \delta$, we obtain from Proposition \ref{basicprop} with $A=\{ x \in \Rd_+ \, : \, |x|>1\}$ that
\begin{eqnarray} \label{sec7.prop.21}
\limsup_{n \to \infty} \P{ \left. {\cal M}_{m_n} > n^{1/\alpha}  \Delta \,\right| \, V_0 = v} 
  &\le& 1- \exp\left\{- K_A \, {\mathbb E}_{\pi_\DD}\left[\tau \right] \cdot 2 \delta \Delta^{-\alpha} \right\} \nonumber\\[.1cm]
  &=& 2K_A \, {\mathbb E}_{\pi_\DD} \left[\tau \right] \Delta^{-\alpha} t, \quad \mbox{where} \:\: t \in (0,\delta),
\end{eqnarray}
and the right-hand side is $\le \Delta/2$ when $\delta$ is chosen sufficiently small.  Note that \eqref{sec7.prop.21} holds for all $v \in {\mathbb D}$.
Finally, let $\gamma_n$ denote the distribution function of $V_{\lfloor n - n\delta \rfloor}$.  By the positive Harris recurrence of $\{V_n\}$ (Lemma \ref{georecurrent}), we have that $|\gamma_n - \pi|_{\rm TV} \to 0$ as $n \to \infty$. Then, using Fatou's lemma, we deduce that
\begin{align*} 
\hspace*{2.7cm}
 \limsup_{n \to \infty}{\mathbb P} &\left( \max_{|m-n| <  n\delta} \left| {\cal M}_m - {\cal M}_n \right| > n^{1/\alpha} \Delta \right) \nonumber \\[.2cm] ~\le~&\limsup_{n \to \infty} \left( \gamma_n(\DD^c) + |\gamma_n - \pi|_{\rm TV} + \int_{\DD} \P{ \left. {\cal M}_{m_{{n}}} > {n}^{1/\alpha}  \Delta \,\right| \, V_0 = v} \, \pi(dv) \right)
\nonumber\\[.2cm]
 \le~& \pi(\DD^c)  +  \frac{\Delta}{2} \pi(\DD) ~\le~ \Delta. \hspace*{8.7cm} \halmos
\end{align*}

\medskip

\noindent
{\bf Proof of Theorem \ref{thm2B}.} 
 Assuming that $d_A$ is bounded, the first assertion follows from Corollary \ref{cor4-1} and from the  uniformity provided by Lemma \ref{sublemma2}.    To remove the assumption that $d_A$ is bounded, see Step 4 below.

To establish the remaining assertion, we proceed in four steps.

\Step[1].  First assume that $d_A$ is bounded and continuous and that
$(H_3)$ is satisfied. 
We claim that
\begin{equation} \label{sec7.mainth.0}
\lim_{n \to \infty} \P{{\cal M}_n \le n^{1/\alpha} u} = e^{-K_A u^{-\alpha}};
\end{equation}
that is, we can transfer the result for maxima over cycles (Proposition \ref{basicprop}) to the process of running maxima,
namely to ${\cal M}_n$.

To establish an upper bound for $\limsup_{n \to \infty} \P{M_n \le n^{1/\alpha}u}$ observe that, by definition,
 ${\cal M}_{N_{\mathbb D}(n)}^U$ corresponds to the value of the process $ \{{\cal M}_j \}$ 
evaluated at the time of its last visit to ${\mathbb D}$ within the time interval $[0,n]$. Thus
\begin{equation} \label{sec7.mainth.1}
\P{{\cal M}_n > n^{1/\alpha} u} \ge \P{{\cal M}_{N_{\mathbb D}(n)}^U > n^{1/\alpha}u}.
\end{equation}

To replace the random time $N_\DD(n)$ by a fixed time, observe by Lemma \ref{lemma:hitting times Vn} that for all $\delta>0$,
\begin{equation} 
\P{\left| \frac{N_{\mathbb D}(n)}{n} - \pi(\DD) \right| \ge \delta} \to 0 \quad \mbox{\rm as} \:\: n \to \infty. \label{sec7.mainth.2}
\end{equation}
Set $t_n= n\left( \pi(\DD) - \delta \right)$ and $\Omega_n = \left\{ \big| \left(N_{\mathbb D}(n)/n\right) - \pi(\DD)  \big| < \delta \right\}$,
and note that $N_\DD(n) \ge \lfloor t_n \rfloor$ on $\Omega_n$.  Then
\begin{equation}  \label{sec7.mainth.3}
\P{{\cal M}_{N_{\mathbb D}(n)}^U > n^{1/\alpha}u} \ge \P{ {\cal M}_{N_{\mathbb D}(n)}^U > n^{1/\alpha}u; \, \Omega_n } 
\ge \P{ {\cal M}_{\lfloor t_n \rfloor}^U > n^{1/\alpha}u } - \P{\Omega_n^c}.
\end{equation}  
Then combining \eqref{sec7.mainth.1}, \eqref{sec7.mainth.2}, and \eqref{sec7.mainth.3} and applying Proposition \ref{basicprop}, we conclude that for all $\delta >0$,
\begin{eqnarray*} 
\liminf_{n \to \infty} \P{{\cal M}_n > n^{1/\alpha} u}  \ge \lim_{n \to \infty} \P{ {\cal M}_{\lfloor t_n \rfloor}^U > n^{1/\alpha}u }
 =  1- \exp\Big\{- K_A \E_{\pi_{\mathbb D}}[\tau]\left( \pi(\DD) - \delta \right) u^{-\alpha} \Big\}
\end{eqnarray*}
and hence, letting $\delta \downarrow 0$ and recalling that $\pi(\DD) = \left( \E_{\pi_{\mathbb D}}[\tau] \right)^{-1}$ (Lemma \ref{lemma:hitting times Vn}),
we obtain that 
\begin{equation}\label{sec7.mainth.4a}
\limsup_{n \to \infty} \P{{\cal M}_n \le n^{1/\alpha} u} ~\le~ \exp\{-K_A u^{-\alpha}\}.
\end{equation}

To establish the corresponding lower bound for $\P{{\cal M}_n \le n^{1/\alpha} u} $, observe that for any $\Delta >0$,
\begin{equation} \label{sec7.mainth.8}
\P{{\cal M}_n > n^{1/\alpha} u} \le \P{{\cal M}_{N_{\mathbb D}(n)}^U > n^{1/\alpha}(u-\Delta)} + \P{ \left|{\cal M}_n - {\cal M}_{N_{\mathbb D}(n)}^U \right| > n^{1/\alpha} \Delta}.
\end{equation}
Reasoning as before, we obtain that the first term on the right-hand side of \eqref{sec7.mainth.8} satisfies
\begin{equation} \label{sec7.mainth.8a}
\limsup_{n \to \infty} \P{{\cal M}_{N_{\mathbb D}(n)}^U > n^{1/\alpha}(u-\Delta)} ~\le~1 - \exp\left\{ -K_A(u-\Delta)^{-\alpha} \right\}.
\end{equation}

To quantify the second term on the right-hand side of \eqref{sec7.mainth.8}, first note that ${\cal M}_{N_{\mathbb D}(n)}^U$ corresponds to the value of the process $ \{{\cal M}_j \}$ 
evaluated at the time of its last visit to ${\mathbb D}$ in the interval $[0,n]$.  Since $\kappa_i$ denotes the time of the $i^{\footnotesize {\rm th}}$ visit to ${\mathbb D}$, this gives
 ${\cal M}_{N_{\mathbb D}(n)}^U = {\cal M}_{\kappa_{N_{\mathbb D}(n)}}$.  Moreover, for any $\delta > 0$, it follows from
Lemma \ref{lemma:hitting times Vn} that
\begin{equation*} 
\P{\left|  \frac{\kappa_{N_{\mathbb D}(n)}}{n}  - 1 \right| \ge \delta} \to 0 \quad \mbox{\rm as} \:\: n \to \infty.
\end{equation*}
Set $\widehat{\Omega}_n = \left\{  \left| \left( \kappa_{N_{\mathbb D}(n)} / n \right) - 1 \right|  < \delta \right\}$. 
Then
\begin{eqnarray*}  
\P{ \left|{\cal M}_n - {\cal M}_{N_{\mathbb D}(n)}^U \right| > n^{1/\alpha} \Delta} 
& \le & \P{ \left|{\cal M}_n - {\cal M}_{N_{\mathbb D}(n)}^U \right| > n^{1/\alpha} \Delta; \, \widehat{\Omega}_n } 
                             +  {\mathbb P} \big( \widehat{\Omega}_n^c \big)  \\[.2cm]
& \le & \P{ \max_{|m-n| < n\delta} \left| {\cal M}_m - {\cal M}_n \right| > n^{1/\alpha} \Delta}  + o(1) 
\end{eqnarray*}
as $n \to \infty$.
Hence by Lemma \ref{cont-in-prob},
\begin{equation} \label{sec7.mainth.8b}
\limsup_{n \to \infty} \P{ \left|{\cal M}_n - {\cal M}_{N_{\mathbb D}(n)}^U \right| > n^{1/\alpha} \Delta} \le \Delta.
\end{equation}
Finally, substituting \eqref{sec7.mainth.8a} and \eqref{sec7.mainth.8b} into \eqref{sec7.mainth.8} and letting $\Delta \to 0$, we obtain that
$$ \liminf_{n \to \infty} \P{ {\cal M}_n \le n^{1/\alpha} u} ~\ge~\exp \big\{ -K_A u^{-\alpha} \big\}.$$
Together with \eqref{sec7.mainth.4a}, this implies the assertion.

\medskip

\Step[2]. Next we remove the additional assumption ($H_3$), but still assume that the function $d_A$ is bounded and continuous.

To remove $(H_3)$, we employ the smoothing argument  introduced in Section \ref{sect:5}.   Let $\{(\widehat{M}_n,\widehat{Q}_n:
n=0,1,\ldots \}$ be defined as in \eqref{MhatQhat}.
Then, since we have assumed that $k=1$ in $({\mathfrak K})$, it follows
that $(\widehat{M}_n, \widehat{Q}_n) = (M_n,Q_n)$ for all $n \in \pintegers$. 
This gives that $V_{n,\epsilon}\le V_n \le V_n^\epsilon$ for all $n \in \pintegers$. 

By repeating the computation  leading to \eqref{eq:distVVeps}, we obtain that
$$ \abs{\widetilde{V}_{n,\epsilon} - \widetilde{V}_n} \le \frac{2\epsilon}{s} \quad \text{ and } \quad \abs{\widetilde{V}_{n}^\epsilon - \widetilde{V}_n} \le \frac{2\epsilon}{s} \qquad \text{ for all } n \in \pintegers.$$ Since the function $d_A$ is assumed to be continuous on the compact set $\SPdp$, it is, in fact, equicontinous.  Hence there is a sequence $\delta(\epsilon)$, tending to zero as $\epsilon \to 0$, such that
\begin{equation}
\abs{d_A\big(\widetilde{V}_{n,\epsilon}\big) - d_A\big(\widetilde{V}_n)} \le \delta(\epsilon) \quad \text{ and } \quad \abs{d_A\big(\widetilde{V}_{n}^\epsilon\big) - d_A\big(\widetilde{V}_n\big)} \le \delta(\epsilon) \qquad \text{ for all } n \in \pintegers.
\end{equation}
Since $A \subset \{ v: |v| > 1 \} \Longrightarrow d_A > 1$, it follows that for all $n \in \pintegers$,
\begin{align*}
|V_n| \le u \cdot d_A\big(\widetilde{V}_n\big) ~\Longrightarrow~|V_{n,\epsilon}| \le u \cdot d_A\big(\widetilde{V}_n\big) \le u \cdot d_A\big(\widetilde{V}_{n,\epsilon}\big) + u \cdot \delta(\epsilon) \le u \cdot d_A\big(\widetilde{V}_{n,\epsilon}) \big( 1+ \delta(\epsilon)\big).
\end{align*}
Consequently, $\P{{\cal M}_n \le u} \le \P{{\cal M}_{n,\epsilon} \le u(1+\delta(\epsilon))}$.  Similarly, for all $n \in \pintegers$,
$$ 
|V_n^\epsilon| \le u \cdot d_A\big(\widetilde{V}_n^\epsilon\big) \big(1-{\delta(\epsilon)}\big) ~\Longrightarrow~|V_{n}| \le u \cdot d_A\big(\widetilde{V}_n^\epsilon\big) - u\cdot \delta(\epsilon) \le u \cdot d_A\big(\widetilde{V}_{n}\big),  $$
and we obtain that $\P{{\cal M}_n \le u} \ge \P{{\cal M}_{n}^\epsilon \le u(1-\delta(\epsilon))}$. Using these upper and lower bounds 
together with Step 1, we conclude that
\begin{align} \label{eq:smoothing hitting time}
\exp &\left\{-K^\epsilon(u-\delta(\epsilon))^{-\alpha}\right\}  \le \liminf_{n \to \infty} \P{{\cal M}_n \le n^{1/\alpha} u} \nonumber \\[.2cm]
& \hspace*{.5cm} \le \limsup_{n \to \infty} \P{{\cal M}_n \le n^{1/\alpha} u} \le \exp \left\{-K_\epsilon(u+\delta(\epsilon))^{-\alpha}\right\},
\end{align}
for constants $K_\epsilon := C_\epsilon D_A$ and $K^\epsilon := C^\epsilon D_A$, 
where $C_\epsilon$ and $C^\epsilon$ are given as in Section \ref{sect:5}.
  Now by Lemma \ref{smoothingLm-2}, $\limsup_{\epsilon \to 0} C_\epsilon \le C$ and $\liminf_{\epsilon \to 0} C^\epsilon \ge C$. Thus, letting $\epsilon \downarrow 0$, we obtain \eqref{sec7.mainth.0}.

\Step[3].  We now relate the behavior of the maxima to the behavior of the first passage times.

Recall that $V_n \in u A \Longleftrightarrow |V_n| > u \cdot d_A(\widetilde{V}_n) \Longleftrightarrow V_n^A >u$.  Hence for all $n \in \pintegers$ and all $u > 0$,
\begin{equation} \label{sec7.18}
{\mathbb P} \left( T_u^A \le n \right) = {\mathbb P} \left( V_i^A > u, \mbox{ some } 1 \le i \le n \right) =  {\mathbb P} \left( {\cal M}_n > u \right).
\end{equation}
Then by \eqref{sec7.mainth.0} and \eqref{sec7.18},
\[
\lim_{n \to \infty} {\mathbb P} \Big( T_{n^{1/\alpha}v}^A \le n \Big) = 1 - e^{-K_A \, v^{-\alpha}}.
\]
Finally, setting $u = n^{1/\alpha} v$ and $z = v^{-\alpha}$ yields
\begin{equation} \label{sec7.19}
\lim_{v \to \infty} {\mathbb P} \left( \frac{T_{u}^A}{u^\alpha} \le z \right) = 1 - e^{-K_A \, z}, \quad z \ge 0.
\end{equation}

\Step[4].  Finally suppose that
${\mathfrak P}_A := \{ x \in \SPdp \, : \, d_A(x) < \infty\} \neq {\mathbb S}_+^{d-1}$.  For any $L \ge 1$, set 
\[
{\cal K}_L = \{ v \in \realsd_+:  |v| \ge L \} \quad \mbox{and} \quad A_L = A \cup {\cal K}_L.
\] 

First observe that $d_{{\cal K}_L}(x) := \inf\{ t:  tx \in {\cal K}_L \} = L$, for all $x \in \SPdp$.  Hence,
letting $r^{{\cal K}_L}_\alpha$ be defined as in \eqref{def-S*},
we have that $r^{{\cal K}_L}_\alpha (x) = L^\alpha r_\alpha(x)  \uparrow \infty$ as $L \to \infty$, uniformly in $x$, since
$r_\alpha$ is uniformly bounded from below by a positive constant, by Lemma \ref{prop:ip}.
Now in general, the constant $K_A$ is proportional to $D_A$, where the latter constant  was characterized in Lemma \ref{ruin.est}.  Using this characterization,
we see that $r_\alpha^{{\cal K}_L}(x) \uparrow \infty \: \forall x \Longrightarrow D_A^{{\cal K}_L} \downarrow 0$ as $L \to \infty$.  Consequently,
\begin{equation} \label{sec7.22a}
\Delta(L):= \lim_{u \to \infty} {\mathbb P} \left( \frac{T_u^{{\cal K}_L}}{u^\alpha} \le z \right) \searrow 0 \quad \mbox{as} \quad L \to \infty.
\end{equation}
Since
\begin{equation*}
\left| {\mathbb P} \left( \frac{T_u^{A_L}}{u^\alpha} \le z \right) -  {\mathbb P} \left( \frac{T_u^{A}}{u^\alpha} \le z \right) \right| 
   \le  {\mathbb P} \left( \frac{T_u^{{\cal K}_L}}{u^\alpha} \le z \right),
\end{equation*}
we conclude that for all $z \ge 0$,
\begin{align} \label{sec7.24}
- \Delta(L) + \lim_{u \to \infty} {\mathbb P} \left( \frac{T_u^{A_L}}{u^\alpha} \le z \right)& \le  \liminf_{u \to \infty} {\mathbb P} \left( \frac{T_u^{A}}{u^\alpha} \le z \right) \le \limsup_{u \to \infty} {\mathbb P} \left( \frac{T_u^{A}}{u^\alpha} \le z \right) \nonumber \\[.2cm]
& \le  \Delta(L) + \lim_{u \to \infty} {\mathbb P} \left( \frac{T_u^{A_L}}{u^\alpha} \le z \right).
\end{align}
Thus, by \eqref{sec7.22a} and Step 3,
\begin{equation} \lim_{u \to \infty} {\mathbb P} \left( \frac{T_u^{A}}{u^\alpha} \le z \right) ~=~1 - \lim_{L \to \infty} \exp\big\{-C D_{A_L} z \big\} ~:=~ 
  1- \exp\big\{-C D_A z \big\}.\label{sec7.24a} \end{equation}
Observe that $D_A := \lim_{L \to \infty} D_{A_L}$ exists, since 
$D_{A_L} =u^\alpha\P{\mathfrak{T}_u^{A_L} < \infty \, | \, V_0 \sim \pi_\DD}$ is a decreasing sequence;
that is to say, the hitting probability of a decreasing sequence of sets.

It remains to identify $D_A$ as the ruin constant in this case.
Arguing as before, we have for all $u >0$ that
\begin{equation*} 
 \Big| u^\alpha\P{\mathfrak{T}_u^{A_L} < \infty \, | \, V_0 =v} - u^\alpha\P{\mathfrak{T}_u^A < \infty \, | \, V_0 =v }\Big| 
   \le  u^\alpha\P{\mathfrak{T}_u^{{\cal K}_L} < \infty \, | \, V_0 =v } ~\searrow 0~ \text{as } L \to \infty.
\end{equation*}
Thus, we deduce by another sandwich argument that
$$ \lim_{u \to \infty} u^\alpha\P{\mathfrak{T}_u^A < \infty \, | \, V_0 = v} = \lim_{L \to \infty} \bigg( \lim_{u \to \infty} u^\alpha\P{\mathfrak{T}_u^{A_L} < \infty \, | \, V_0  = v } \bigg) = r_\alpha(\widetilde{v})  \lim_{L \to \infty} D_{A_L} = r_\alpha (\widetilde{v}) D_A,$$
which gives the required identification of $D_A$ as the constant in the ruin problem for the Markov random walk; cf.\ Lemma \ref{ruin.est}.

To conclude the proof, observe that the same reasoning yields \eqref{thm2B-eq1} for unbounded functions $d_A$; namely, one
can again introduce the family $A_L = A \cup {\cal K}_L$ for $L \ge 1$, and argue that the hitting probability
of the set ${\cal K}_L$---now prior to the return time $\tau$---becomes
asymptotically negligible as $L \to \infty$.  The argument is entirely identical, so we omit the details.
 \halmos

%
%
%
%

\setcounter{equation}{0}
\section{Determining the path of large exceedance}\label{sect:7}
We conclude by studying the path of large exceedance conditioned on $\{ T_u^A < \tau \}$.  
In particular, we provide the proofs of Theorem  \ref{thm:empirical law}
and a stronger version of Theorem \ref{thm:excursion}, where we also allow for paths of infinite length. 

To study paths of infinite length in the context of Theorem \ref{thm:excursion}, 
first introduce the normalized process
\[
 \bar{S}_n = S_n - n \log a, \quad n = 0,1,\ldots,
\]
where $a \ge 0$.  A natural choice is $\log a = \Lambda^\prime(\alpha) = \Estat[S_1]$, in which case 
 $n^{-1} \bar{S}_n \to 0$ a.s.\ as $n \to \infty$, by Lemma \ref{lem:slln}.
Similarly, we introduce a normalization for $\{V_n\}$, describing the behavior after this process has exceeded an initial barrier $u_\epsilon$, where $u_\epsilon = o(u)$ 
and $u_\epsilon \uparrow \infty$ as $u \to \infty$.  Conditioned on $\{ T_u^A < \tau \}$,
our objective is to characterize $\{ V_n \}$ over the time interval $[T_{u_\epsilon}, T_u^A]$, and to show that this process
resembles the process $\{e^{S_n} X_n\}$ under ${\mathbb P}^\alpha$.

Let 
\begin{equation} \label{def-IuJu}
\gamma_u = u - \epsilon_u, \quad \mbox{and set} \quad  I_u = T_{\epsilon_u} \quad  \mbox{and} \quad J_u = T_{\gamma_u}^A.
\end{equation}

Now by the nonlinear renewal theory in Subsection 3.3  (cf.\ Eq. \eqref{aug9-3} in Lemma \ref{lemma-charac-T}), 
$\log|V_{I_u+n}| - \log |V_{I_u}|$ grows under ${\mathbb P}^\alpha$ at roughly the rate $n \Lambda^\prime(\alpha)$ as $n \to \infty$,
where $\Lambda^\prime(\alpha)$ represents the mean of the Markov random walk $\{ (X_n,S_n)\}$.  Thus, to describe the 
post-$I_u$ behavior of $\{ V_n \}$, it is natural to consider the normalized process
\[
\bar{V}_{I_u + n}^{(u)} := \frac{1}{a^n}  \frac{V_{I_u + n}}{|V_{I_u}|} , \quad n=0,1,\ldots,
\]
where $\log a = \Lambda^\prime(\alpha)$ and
and $\bar{V}_k^{(u)} := 0$ for all $k < I_u$.  
In practice, when studying the behavior of $\{ V_n \}$ over a path of {\it infinite} length, we generally need to consider other choices of $a$, where $\log a> \Lambda'(\alpha)$, and thus $a$ should be
viewed here as an arbitrary parameter subject to the constraint that $\log a \ge 0$.

Finally, given a measurable function $g:  (\reals^d)^{\mathbb N} \to \reals$ and $m \in \pintegers$, define ${\cal P}_{[m]} g:  (\reals^d)^m \to \reals$ by setting
${\cal P}_{[m]} g (x_1, \ldots, x_m) = g(x_1, \ldots, x_m,0,0,\ldots)$.  Thus, ${\cal P}_{[m]} g$ is determined by $g$ under the {\it projection} of $(x_1,x_2,\ldots)$ onto $(x_1,\ldots,x_m)$.  
Also, extending the standard finite-dimensional definition, we say that $g:  (\reals^d)^{\mathbb N} \to \reals$ is $\theta$-H\"older continuous if for all sequences $\{x_i\}_{i\in \N_+}$, $\{y_i\}_{i \in \N_+} \in  (\reals^d)^{\mathbb N}$,
$$ \Big| g (x_1, x_2, \dots,) - g(y_1, y_2, \dots) \Big| \le \sum_{i \ge 1} \abs{x_i-y_i}^\theta. $$


In the next theorem, note that ${\cal P}_{[k]} g = {\cal P}_{[m]} g$ for all $k \ge m \Longleftrightarrow
g(x_1,x_2,\ldots) = g(x_1,\ldots,x_m,0,$ $0,\ldots) := g^\ast(x_1, \ldots, x_m)$.  Recall that ${\mathfrak D}$ denotes the intersection
of the domain of $\Lambda$ with $\reals_+$.

\begin{thm} \label{prop8.1}
Suppose that Hypotheses $(H_1)$ and $(H_2)$ are satisfied, and assume that the function $d_A$ is finite and continuous on $\SPdp$.
Let $\{ \epsilon_u:  u \in \reals_+ \}$ be any sequence such that $\epsilon_u = o(u)$ and $\epsilon_u \uparrow \infty$ as $u \to \infty$.
Let $I_u$ and $J_u$ be defined as in \eqref{def-IuJu}, and set $g_u = {\cal P}_{[J_u-I_u]}g$,
where $g:  \reals^{\mathbb N} \to \reals$ is a bounded 
measurable function.  Assume that either: 
\begin{itemize}
\item[{\rm (i)}] ${\cal P}_{[k]} g = {\cal P}_{[m]} g$ for all $k \ge m,$ some $m \in \pintegers$, and that $g$ is $\theta$-H\"{o}lder continuous for some $\theta \le \min\{1,\alpha\}${\rm ;} or
\item[{\rm (ii)}] $g$ is $\theta$-H\"{o}lder continuous for some $0 <\theta \le 1$ satisfying $\alpha + \theta \in {\mathfrak D}$, 
where $\E[\norm{M}^{\alpha} |Q|^\theta] < \infty$ and
 $\log a> \Lambda(\alpha + \theta)/\theta$.
 \end{itemize}
Then for all $v \in \R_+^d$,
\begin{align}\label{conditional limit law}
 \lim_{u \to \infty} &\, \E\bigg[ g_u \bigg( \bar{V}_{I_u}^{(u)}, \ldots, \bar{V}_{J_u}^{(u)} \bigg) \, \bigg| \, 
    V_0=v,\, T_u^A < \tau \bigg] \nonumber\\[.1cm]
 & \hspace*{1.5cm} =~ \int_{{\mathbb S}^{d-1}_+ \times \reals_+}
 \E_x^\alpha \bigg[ g \big(X_0, e^{\bar{S}_1} X_1, e^{\bar{S}_2}X_2, \ldots)\bigg] \rho(dx,ds).
 \end{align} 
\end{thm}

Note by  definition that $|\bar{V}_{I_u}^{(u)}|=1$, and consequently the initial value on the right-hand side is $X_0 \in \SPdp$,
not $e^s X_0$ [as would appear if, instead, we had normalized $V_{I_u}$ by dividing by $\epsilon_u$].  Further, since $\Lambda$ is a convex function, we have that $\Lambda(\alpha + \theta) \ge \Lambda'(\alpha) \theta$ for $\Lambda(\alpha)=0$. Thus, in (ii),
we require $\log a> \Lambda'(\alpha)$.

\medskip

\noindent
{\bf Proof.}
To establish the result, it suffices to prove that 
\begin{align}\label{eq: behavior during excursion}
&\lim_{u \to \infty}  u^\alpha \E\bigg[ g_u \bigg( \bar{V}_{I_u}^{(u)}, \ldots, \bar{V}_{J_u}^{(u)} \bigg)
  \mathbf{1}_{\{T_u^A < \tau\}} \bigg|  \, V_0=v  \bigg] 
  =  r_\alpha(\widetilde{v}) \E_{\bfDV}^\alpha \left[ |Z|^\alpha \mathbf{1}_{\{ \tau=\infty\}} \right]  \nonumber\\
& \hspace*{2cm} \times  \int_{{\mathbb S}^{d-1}_+ \times \reals_+} 
       \E_x^\alpha \bigg[ g \big( X_0, e^{\bar{S}_1} X_1, e^{\bar{S}_2} X_2, \ldots)\bigg] \rho(dx,ds) 
        \left( \int_{{\mathbb S}^{d-1}_+ \times \reals_+} \frac{e^{-\alpha s}}{r^A_\alpha(x)} \, \rho^A(dx,ds) \right).
\end{align} 
Once \eqref{eq: behavior during excursion} is established,
then the assertion follows by using \eqref{eq: behavior during excursion} twice (once with $g = 1$), observing that
\begin{equation} \label{sec6A.1}
\E\bigg[ g_u \bigg( \bar{V}_{I_u}^{(u)}, \ldots, \bar{V}_{J_u}^{(u)} \bigg) \, \bigg| \, 
    V_0=v,\, T_u^A < \tau \bigg] = \frac{{\cal A}_1}{{\cal A}_2},
\end{equation}
where
\begin{equation*}
{\cal A}_1 = u^\alpha \E\bigg[ g_u \bigg(   \bar{V}_{I_u}^{(u)}, \ldots, \bar{V}_{J_u}^{(u)} \bigg)
  \mathbf{1}_{\{T_u^A < \tau\}} \bigg| \, V_0=v  \bigg]  \quad \mbox{and} \quad
{\cal A}_2 = u^\alpha \E\bigg[  \mathbf{1}_{\{T_u^A < \tau\}} \bigg| \, V_0=v  \bigg].
\end{equation*}

To verify \eqref{eq: behavior during excursion}, we proceed as in the proof of Proposition \ref{prop1}, 
converting to the 
$\alpha$-shifted measure to
obtain that
\begin{align}  \label{eq:decompose excursion behavior}
u^\alpha & \E\bigg[ g_u \bigg(  \bar{V}_{I_u}^{(u)}, \ldots, \bar{V}_{J_u}^{(u)} \bigg)
  \mathbf{1}_{\{T_u^A < \tau\}} \bigg| \,  V_0=v  \bigg] \nonumber\\[.2cm]
&=~ u^\alpha r_\alpha(\widetilde{v}) \, 
\E_{\bfDV}^\alpha\bigg[ \frac{e^{-\alpha S_{T_u^A}}}{r_\alpha(X_{T_u^A})} g_u \bigg( \bar{V}_{I_u}^{(u)}, \ldots, \bar{V}_{J_u}^{(u)} \bigg)
  \mathbf{1}_{\{T_u^A < \tau\}}   \bigg] \nonumber\\[.2cm]
&=~ r_\alpha(\widetilde{v}) \, \E_{\bfDV}^\alpha \Big[ |Z_{T_u^A}|^\alpha {\mathfrak G}_u {\bf 1}_{\{T_u^A < \tau \}} \Big],
\end{align}
where $Z_n:= |V_n|/e^{S_n}$, $n=0,1,\ldots,$ and
\begin{equation} \label{def-G-U}
{\mathfrak G}_u ~:=~ \frac{1}{r^A_\alpha(X_{T_u^A}) } \left( \frac{\big|V^A_{T_u^A}\big|}{u}\right)^{-\alpha}
 \left(  \frac{d_A(X_{T_u^A})}{d_A(\widetilde{V}_{T_u^A})} \right)^\alpha
  \: g_u \bigg(  \bar{V}_{I_u}^{(u)}, \ldots, \bar{V}_{J_u}^{(u)} \bigg).
   \end{equation}
[The term $\left( d_A(X_{T_u^A})/d_A(\widetilde{V}_{T_u^A}) \right)^\alpha$ arises when we replace $|V_{T_u^A}|^{-\alpha}/r_\alpha(X_{T_u^A})$ with
$|V^A_{T_u^A}|^{-\alpha}/r^A_\alpha(X_{T_u^A})$.]
The right-hand side of \eqref{eq:decompose excursion behavior} can be written as the difference of two terms, namely 
\begin{align} \label{eq:decompose excursion behavior II}
&r_\alpha(\widetilde{v}) {\mathbb E}_{\bfDV}^{\alpha}  \left[ \bigg( \big| Z_{T_u^A} \big|^\alpha {\bf 1}_{\{ T_u^A < \tau\}} - \abs{Z_n}^\alpha {\bf 1}_{\{n \le T_u^A\}} 
   {\bf 1}_{\{ n \le \tau\}} \bigg) {\mathfrak G}_u \right] \nonumber \\[.2cm]
& \hspace*{3cm} +~ r_\alpha(\widetilde{v}) {\mathbb E}_{\bfDV}^{\alpha}  \bigg[ \abs{Z_n}^\alpha {\bf 1}_{\{n \le T_u^A\}} {\bf 1}_{\{ n \le \tau\}}  \, \E^{\alpha} \big[ {\mathfrak G}_u \, | \F_n \big] \bigg]. 
\end{align} 
As in the proof of Proposition \ref{prop1}, we may then apply Lemma \ref{sublemma} (i) and use the uniform boundedness of $\{ {\mathfrak G}_u \}$
to conclude that the first term in \eqref{eq:decompose excursion behavior II} tends to zero as $u \to \infty$ and then $n \to \infty$.

Thus, it suffices to analyze the second term in \eqref{eq:decompose excursion behavior II}.  Indeed, the proof of the theorem will be complete 
once we have established the following.

\begin{lemma} \label{lemma8.2}  Under the conditions of Theorem \ref{prop8.1},
\begin{align} \label{sec6A.99}
\lim_{n \to \infty} \lim_{u \to \infty} & {\mathbb E}_{\bfDV}^{\alpha} \bigg[ \abs{Z_n}^\alpha {\bf 1}_{\{n \le T_u^A\}} {\bf 1}_{\{ n \le \tau\}}  \, \E^{\alpha} \big[ {\mathfrak G}_u \, | \F_n \big] \bigg] \nonumber\\[.2cm]
& = H r_\alpha(\tilde{v}) \E_{\bfDV}^\alpha \left[ |Z|^\alpha \mathbf{1}_{\{ \tau=\infty\}} \right]  \int_{{\mathbb S}^{d-1}_+ \times \reals_+} 
      \E_x^\alpha \bigg[ g \big( X_0, e^{\bar{S}_1} X_1, e^{\bar{S}_2} X_2, \ldots)\bigg] \rho(dx,ds),
\end{align}
where
$H :=   \int_{{\mathbb S}^{d-1}_+ \times \reals_+} \big(e^{-\alpha s}/r^A_\alpha(x)\big) \, \rho^A(dx,ds).$
\end{lemma}

\noindent
{\bf Proof of Lemma \ref{lemma8.2}.}
\Step{1}. We begin by analyzing $\E^\alpha[{\mathfrak G}_u|{\cal F}_n]$,
showing that there is asymptotic independence between $\{ V_n:  n \le J_u \}$ and $V_{T_u^A}^A/u$.
Set
\begin{equation} \label{eq:hu}
h(u) ~=~  {\mathbb E}^\alpha \left[ \left. \frac{1}{r^A_\alpha(\widetilde{V}_{T_u^A})} \left( \frac{ \big| V^A_{T_u^A} \big|}{u} \right)^{-\alpha}   \left( \frac{d_A(X_{T_u^A})}{d_A(\widetilde{V}_{T_u^A})} \right)^\alpha \frac{r^A_\alpha(\widetilde{V}_{T_u^A})}{r^A_\alpha(X_{T_u^A})}\:\right| {\cal F}_{J_u} \right].
\end{equation}
Recall that $J_u=T_{u-\epsilon_u}^A$, where $u - \epsilon_u \to \infty$. Now decompose $h$ into two parts, namely
$$h(u) = h(u) {\mathbf 1}_{\left\{|V_{J_u}^A| \le u - \frac{\epsilon_u}{2} \right\}} + h(u) {\mathbf 1}_{\left\{|V_{J_u}^A| > u - \frac{\epsilon_u}{2} \right\}},$$ 
and observe that the second term tends to zero in $\Prob^\alpha$-probability, since the sequence $(u - \epsilon_u)^{-1} |V_{T_u}^A|$ is tight, by Lemma \ref{lem:tightness}.  

To analyze the first term, we employ nonlinear renewal theory. Since $d_A^\alpha$ and $r_\alpha^A$ are continuous functions and bounded 
from below, Lemma \ref{lem:VTutoXTu} implies that $\left(d_A(X_{T_u^A})/d_A(\widetilde{V}_{T_u^A}) \right)^\alpha$ and $r^A_\alpha(\widetilde{V}_{T_u^A})/r^A_\alpha(X_{T_u^A})$  tend to one in ${\mathbb P}^\alpha$-probability. Moreover, by the continuous mapping theorem and Theorem \ref{thm:Melfi nonlinear renewal}, 
$\big( r^A_\alpha(\widetilde{V}_{T_u^A}) \big)^{-1} \exp\{-\alpha(\log |V^A_{T_u^A}| - \log u)\}$ converges in law, independent of the initial distribution. Then by the Markov property and Slutsky's theorem, 
\begin{align} \label{sec6A.10}
 h(u)  {\mathbf 1}_{ \left\{|V_{J_u}^A| \le u - \frac{\epsilon_u}{2} \right\}}  & =   {\mathbb E}^\alpha_{X_{J_u},
  V_{J_u}} \left[  \frac{ \exp\big\{-\alpha\big( \log |V^A_{T_u^A}| - \log u\big) \big\} }{r^A_\alpha(\widetilde{V}_{T_u^A})}  \left( \frac{d_A(X_{T_u^A})}{d_A(\widetilde{V}_{T_u^A})} \right)^\alpha \frac{r^A_\alpha(\widetilde{V}_{T_u^A})}{r^A_\alpha(X_{T_u^A})}  \right] {\mathbf 1}_{\left\{|V_{J_u}^A| \le u - \frac{\epsilon_u}{2} \right\}} \nonumber\\[.2cm]
&  \longrightarrow  \int_{[0,\infty) \times {\mathbb S}^{d-1}_+} \: \frac{e^{-\alpha s}}{r^A_\alpha(x)} \rho^A(dx,ds) := H \quad \mbox{in }{\mathbb P}^\alpha \mbox{-probability as }u \to \infty.
\end{align}
Thus we obtain the convergence $h(u) \to H$ in $\Prob^\alpha$-probability. 
Next observe that $h(u)$ is uniformly bounded, since we have assumed that $d_A$ is bounded and continuous.
Hence, since $g$ bounded $\Longrightarrow g_u$ bounded, 
it follows that
$(h(u)-H)g_u$ tends to zero in $\Prob^\alpha$-probability and also in $L^1$.
 Using the Markov property, we then deduce that on $\{n \le T_u^A\}$,
 \begin{align} \label{sec6A.11} 
\lim_{u \to \infty} {\mathbb E}^\alpha \left[ {\mathfrak G}_u | {\cal F}_n \right] & =  \lim_{u \to \infty} {\mathbb E}^\alpha_{X_n,V_n} \left[ \left(h(u)-H\right) \, g_u \Big(   \bar{V}_{I_u}^{(u)}, \ldots, \bar{V}_{J_u}^{(u)} \Big) \right] 
+ \lim_{u \to \infty} H {\mathbb E}^\alpha_{X_n,V_n} \left[  g_u \Big(   \bar{V}_{I_u}^{(u)}, \ldots, \bar{V}_{J_u}^{(u)} \Big) \right] \nonumber\\[.2cm]
& = H \cdot  \lim_{u \to \infty}  {\mathbb E}^\alpha_{X_n,V_n} \left[  g_u \Big(   \bar{V}_{I_u}^{(u)}, \ldots, \bar{V}_{J_u}^{(u)} \Big) \right]  \qquad \Prob^\alpha\text{-a.s.} 
\end{align}

\Step{2}.  Next assume that either (i) or (ii) holds.  Then we claim that for all $x \in \SPdp,$ $v \in \reals_+^d$,
\begin{equation} \label{eq:guVtogX} 
\lim_{u \to \infty}  {\mathbb E}^\alpha_{x,v} \left[  g_u \bigg(   \bar{V}_{I_u}^{(u)}, \ldots, \bar{V}_{J_u}^{(u)} \bigg) \right] ~=~ \int_{{\mathbb S}^{d-1}_+ \times \reals_+} 
       \E_y^\alpha \bigg[ g \big( X_0, e^{\bar{S}_1} X_1, e^{\bar{S}_2} X_2, \ldots)\bigg] \rho(dy,ds)
\end{equation}

We focus on the proof under the set of assumptions (ii), as the calculations needed under (i) are essentially identical, except
simpler.

Thus assume (ii) holds and, for $m \in \N_+$, consider the decomposition 
\begin{align} \label{july10a}
 {\mathbb E}^\alpha_{x,v} & \left[  g_u \Big(   \bar{V}_{I_u}^{(u)}, \ldots, \bar{V}_{J_u}^{(u)} \Big) \right]  - \int_{{\mathbb S}^{d-1}_+ \times \reals_+} 
       \E_y^\alpha \bigg[ g \big( X_0, e^{\bar{S}_1} X_1, e^{\bar{S}_2} X_2, \ldots)\bigg] \rho(dy,ds)\nonumber\\
=~& {\mathbb E}^\alpha_{x,v} \left[  (g_u - {\cal P}_{[m]} g)\Big(   \bar{V}_{I_u}^{(u)}, \ldots, \bar{V}_{J_u}^{(u)} \Big) \,  \right] 
\nonumber\\
& + {\mathbb E}^\alpha_{x,v} \left[  {\cal P}_{[m]} g\Big(   \bar{V}_{I_u}^{(u)}, \ldots, \bar{V}_{I_u+m}^{(u)} \Big) \right] - \int_{{\mathbb S}^{d-1}_+ \times \reals_+} \E_y^\alpha \bigg[ {\cal P}_{[m]} g \big( X_0, e^{\bar{S}_1} X_1, \dots,  e^{\bar{S}_m} X_m)\bigg] \rho(dy,ds)    \nonumber\\
&  + \int_{{\mathbb S}^{d-1}_+ \times \reals_+} \E_y^\alpha \bigg[ (g-{\cal P}_{[m]} g) \big( X_0, e^{\bar{S}_1} X_1, \dots,  e^{\bar{S}_m} X_m)\bigg] \rho(dy,ds) \nonumber\\
:=~& {\mathbb I}_1(u,m) + {\mathbb I}_2(u,m) + {\mathbb I}_3(u,m).
\end{align} 

Next we show that $\abs{ {\mathbb I}_1(u,m)} + \abs{ {\mathbb I}_3(u,m)} \le \Delta(m)$ for a sequence $\left\{ \Delta(m) \right\}_{m \in \pintegers}$ which is independent of $u$ and tends to zero as $m \to \infty$. 

First consider ${\mathbb I}_1(u,m)$.  Note by Lemma \ref{lemma-charac-T} that $\Prob^\alpha(J_u-I_u>m) \to 0$ as $u \to \infty$.
Thus it is sufficient to study the restriction to the set $\{J_u - I_u \le m\}$.
 Now let $\{ (M_n^\ast, Q_n^\ast):  n=0,1,\ldots \}$ be a process which is independent of $\{ (M_n,Q_n):  n=0,1,\ldots \}$
but sharing the same distribution function, and let $V_0^\ast$ denote the initial value corresponding to $\{(M_n^\ast,Q_n^\ast)\}$. 
Then using the $\theta$-H\"older continuity of $g$ and the Markov property, together with the subadditivity of $|\cdot|^\theta$, $\theta \le 1$, we obtain  
\begin{align*}
| {\mathbb I}_1(u,m)| ~\le~&\E_{x,v}^\alpha \bigg[ \sum_{j >m} |\bar{V}_{I_u+j}^{(u)}|^\theta \,\bigg] 
~\le~\E_{x,v}^\alpha \bigg[ \sum_{j>m} \E^\alpha_{X_{I_u}, \tilde{V}_{I_u}} \bigg[ a^{-\theta j} \Big|M_j^\ast \cdots M_1^\ast V_0^\ast  + \sum_{k=1}^j M_j^\ast \cdots M_{k+1}^\ast Q_k^\ast  \Big|^\theta \bigg]  \bigg] \\
~\le~&\E_{x,v}^\alpha \bigg[ \sum_{j > m} a^{-\theta j}\, B_1 \E \bigg[ \norm{M_j \cdots M_1}^\alpha \Big( \norm{M_j \cdots M_1 }^\theta + \sum_{k=1}^j |M_j \cdots M_{k+1} Q_k|^\theta\Big)  \bigg]  \bigg] .
\end{align*} 
In the last line, we used the definition of the $\alpha$-shifted measure and set $B_1=\max_{x,y \in \SPdp} \big( r_\alpha(x)/r_\alpha(y) \big)$. Now Corollary 4.6 in \shortciteN{Buraczewski.etal:2014} gives that $\E\big[\norm{M_j \cdots M_{1}}^\beta \big] \le \overline{B}_\beta \left( \lambda(\beta) \right)^{j}$ for all $\beta \in \mathfrak{D}$, where $\overline{B}_\beta$ is a finite constant
which is independent of $j$. Thus we deduce that $\abs{{\mathbb I}_1(u,m)}$ is bounded above by 
\begin{align*}
B_1 \sum_{j > m} & a^{-\theta j} \left( \E\big[\norm{M_j \cdots M_1}^{\alpha + \theta}\big] + \sum_{k=1}^j \E\big[ \norm{M_j \cdots M_{k+1}}^{\alpha + \theta}  \norm{M_k}^{\alpha} |Q_k|^\theta \norm{M_{k-1} \cdots M_1}^\alpha\big]\right) \\
\le~& B_1 \sum_{j > m} \left( \overline{B}_{\alpha+\theta}\frac{\left(\lambda(\alpha+\theta) \right)^j}{a^{\theta j}} + \sum_{k=1}^j \frac{\lambda(\alpha+\theta)^{j-k} \lambda(\alpha)^{k-1}}{a^{\theta j}} \overline{B}_{\alpha+\theta} \overline{B}_{\alpha} \E\left[\norm{M}^{\alpha} |Q|^\theta\right] \right) \\
\le~& B \sum_{j > m} (j+1)\left(\frac{\lambda(\alpha+\theta)}{a^\theta}\right)^j ~:=~\Delta_1(m)
\end{align*} 
for some finite constant $B$. Recall that by assumption,  $\E[\norm{M}^{\alpha+\theta} |Q|^\theta]$ is finite and $\theta \log a > \Lambda(\alpha + \theta)$ $\Longleftrightarrow$
 $a^\theta > \lambda(\alpha + \theta)$. Thus $\Delta_1(m) \to 0$ as $m \to \infty$. 
 
Repeating the same argument for the term ${\mathbb I}_3(u,m)$ yields that
$\abs{ {\mathbb I}_3(u,m) } \le \Delta_2(m)$ where $|\Delta_2(m)| \to 0$ as $m \to \infty$.
[In contrast to the previous calculation, the terms involving $Q_k$ now vanish.] 

Next we show that $\lim_{u \to \infty} {\mathbb I}_2(u,m)=0$ for all $m \in \pintegers$. 
To this end, define
\begin{eqnarray*}
 G_u(x,v)& = &\E^\alpha_x\left[ {\cal P}_{[m]} g\left(\tilde{v},a^{-1}\left(M_1\tilde{v}+\frac{V_1^{(0)}}{v\epsilon_u}\right), \dots,a^{-m}\left(M_m \cdots M_1\tilde{v}+\frac{V_m^{(0)}}{v\epsilon_u} \right) \right) \right];\\[.2cm]
 G(x,v)& =&\E^\alpha_x\bigg[ {\cal P}_{[m]} g\bigg(\tilde{v},a^{-1}M_1\tilde{v}, \dots, a^{-m}M_m \cdots M_1\tilde{v} \bigg) \bigg],
 \end{eqnarray*}
  where $V_k^{(0)}$ is defined as in \eqref{defV_0} for $k \in \pintegers.$
Then the $\theta$-H\"older continuity of $g$ gives that
\begin{equation}\label{eq:commonbound} \lim_{u \to \infty} \sup \bigg\{ \, \abs{G_u(x,v) - G(x,v)} \, : x \in \SPdp, \, v \in \reals^d_+ \setminus \overline{B_1(0)} \bigg\} ~=~0. 
\end{equation}

Now by Lemma \ref{lem:VTutoXTu}, $\lim_{u\to \infty} {X_{I_u}}/{\tilde{V}_{I_u}} \to 1$ in $\Prob^\alpha$-probability. 
Moreover,  under $\Prob^\alpha_{x,v}$, we have by Melfi's nonlinear renewal theorem (Theorem \ref{thm:Melfi nonlinear renewal})
that
$\epsilon_u^{-1}\big( \widetilde{V}_{I_u}, \log | V_{I_u}| \big)$ converges in law to a random variable $(X,S)$, say, having the distribution $\rho$. Then by the continuous mapping theorem and Slutsky's theorem,
\begin{equation} \label{eq:YY} (X_{I_u}, \epsilon_u^{-1}V_{I_u}) ~=~ \left(\frac{X_{I_u}}{\tilde{V}_{I_u}} \tilde{V}_{I_u}, \frac{V_{I_u}}{\epsilon_u} \right) ~\Rightarrow~(X,e^SX) \quad \mbox{\rm as} \:\: u \to\infty. \end{equation} 
In this notation, $\E_{x,v}\big[G(X,e^SX)\big]=\int_{{\mathbb S}^{d-1}_+ \times \reals_+} \E_y^\alpha \big[ {\cal P}_{[m]} g \big( X_0, e^{\bar{S}_1} X_1, \dots,  e^{\bar{S}_m} X_m)\big] \rho(dy,ds)$.
Then, using the Markov property and the boundedness and continuity of $G$, 
\begin{align}
\lim_{u \to \infty} {\mathbb I}_2(u,m) ~=~& \lim_{u \to \infty} \E_{x,v}^\alpha \bigg[ G_u\left({X}_{I_u},\frac{V_{I_u}}{\epsilon_u}\right) - G\left({X}_{I_u},\frac{V_{I_u}}{\epsilon_u}\right)\bigg] \nonumber\\
&\hspace*{3cm}+ \lim_{u \to \infty} \E_{x,v}^\alpha \bigg[ G\left({X}_{I_u},\frac{V_{I_u}}{\epsilon_u}\right)\bigg] - \E_{x}^\alpha \left[G(X,e^SX) \right]   ~=~0.\end{align} 
Note that the same calculation proves \eqref{eq:guVtogX} under the set of assumptions (i).

In conclusion, we have shown that for all $m \in \N_+$,
$$ \lim_{u \to \infty} \abs{{\mathbb E}^\alpha_{x,v} \left[  g_u \bigg(   \bar{V}_{I_u}^{(u)}, \ldots, \bar{V}_{J_u}^{(u)} \bigg) \right] - \int_{{\mathbb S}^{d-1}_+ \times \reals_+} 
       \E_y^\alpha \bigg[ g \big( X_0, e^{\bar{S}_1} X_1, e^{\bar{S}_2} X_2, \ldots)\bigg] \rho(dy,ds)} ~\le~\Delta(m).$$ 
Since $\Delta(m) \to 0$ as $m \to \infty$, we conclude that \eqref{eq:guVtogX} holds.

\medskip

\Step[3].  By combining \eqref{sec6A.11} and \eqref{eq:guVtogX}, we obtain that
$$ \lim_{u \to \infty} {\mathbf 1}_{\{n \le T_u^A\}} \E^\alpha [ {\mathfrak G}_u \, | {\cal F}_n ] ~=~H \int_{{\mathbb S}^{d-1}_+ \times \reals_+} 
       \E_y^\alpha \bigg[ g \big( X_0, e^{\bar{S}_1} X_1, e^{\bar{S}_2} X_2, \ldots)\bigg] \rho(dy,ds) \qquad \Prob^\alpha\text{-a.s.}$$
Finally, \eqref{sec6A.99} is obtained by reasoning as in \eqref{lm4.1.4}.
\halmos

\medskip

\noindent
{\bf Proof of Theorem \ref{thm:empirical law}.} It suffices to show that 
\begin{equation}\label{eq:8.20}
\limsup_{u \to \infty}\, u^\alpha \E_v \left[ \bigg| \frac{1}{T_u^A} \sum_{k=1}^{T_u^A} g\left( \log \bigg(  \frac{|V_k|}{|V_{k-1}|} \bigg) \right)  ~-~ \Estat \left[g( S_1) \right] \bigg| {\mathbf 1}_{\{ T_u^A < \tau \}} \right] ~=~0.
\end{equation} 
For simplicity, we introduce the shorthand notation 
$$\mu_g:=\Estat[g(S_1)] \quad \mbox{and} \quad \Sigma_j^n=\sum_{k=j}^n g\big( \log|V_k| - \log |V_{k-1}|\big).$$ 

Let $\{ \epsilon_u:  u > 0 \}$ be a sequence such that $\epsilon_u = o(u)$ and $\epsilon_u \uparrow \infty$ as $u \to \infty$.
Set $\gamma_u = u - \epsilon_u$ and $J_u = T_{\gamma_u}^A$, and set $B_1=\max_{x,y} \left( r_\alpha(x)/r_\alpha(y) \right)$.
Then from a change of measure, we infer that
\begin{align}  
u^\alpha & \E \bigg[  \abs{ \frac{1}{T_u^A} \Sigma_1^{T_u^A} - \mu_g }
  \mathbf{1}_{\{T_u^A < \tau\}} \:\Big| \: V_0 = v \bigg] 
~=~ r_\alpha(\widetilde{v}) \, 
\E_{\bfDV}^\alpha\bigg[ \frac{e^{-\alpha (S_{T_u^A}-\log u)}}{r_\alpha(X_{T_u^A})} \abs{ \frac{1}{T_u^A} \Sigma_1^{T_u^A} - \mu_g }
  \mathbf{1}_{\{T_u^A < \tau\}}  \bigg] \nonumber\\[.2cm]
& \hspace*{2cm} \le~ B_1 \, 
\E_{\bfDV}^\alpha\bigg[ \abs{ \frac{1}{J_u} \Sigma_1^{J_u} - \mu_g }
   \bigg] ~+~  B_1 \, 
\E_{\bfDV}^\alpha\bigg[ \abs{ \frac{1}{J_u} \Sigma_1^{J_u} - \frac{1}{T_u^A} \Sigma_1^{T_u^A}}   \bigg]~:=~ {\mathbb I}_1(u) 
  + {\mathbb I}_2(u).
\nonumber
\end{align}

First consider ${\mathbb I}_2(u)$.  Note that
$$ \abs{ \frac{1}{J_u}\Sigma_1^{J_u}  - \frac{1}{T_u^A}\Sigma_1^{T_u^A} } ~\le~\abs{\frac{1}{J_u} \left( \Sigma_1^{J_u} - \Sigma_1^{T_u^A} \right)} + \abs{ \left( \frac{1}{J_u} - \frac{1}{T_u^A}\right) \Sigma_{1}^{T_u^A}} ~\le~ 2 \abs{\frac{T_u^A - J_u}{J_u}} \cdot |g|_\infty. $$
By Lemma \ref{lemma-charac-T},  $\left( J_u/\log u \right) \to \left(\lambda^\prime(\alpha) \right)^{-1}$ in $\Prob^\alpha$-probability and $\left((T_u^A - J_u)/\log u \right) \to 0$ in $\Prob^\alpha$-probability. Hence the term inside the expectation of ${\mathbb I}_2(u)$
 tends to zero in $\Prob^\alpha$-probability and, furthermore, is bounded above by $2 |g|_\infty$. Hence, $\limsup_{u \to \infty} {\mathbb I}_2(u) =0$.

Next consider ${\mathbb I}_1(u)$.
Here, our objective is to apply Fatou's lemma and to observe that
 $$\limsup_{u \to \infty}  \abs{ \frac{1}{J_u} \Sigma_1^{J_u} - \mu_g }=0 \quad \Prob^\alpha\mbox{\rm -a.s.}$$  Consider the decomposition
$$ \abs{ \frac{1}{J_u}\Sigma_1^{J_u}  - \mu_g } ~\le~   \frac{1}{J_u} \sum_{k=1}^{J_u} \bigg| g\left( \log \bigg(  \frac{|V_k|}{|V_{k-1}|} \bigg) \right)   - g \big( S_k - S_{k-1} \big)   \bigg|  
 ~+~  \left| \frac{1}{J_u} \sum_{k=1}^{J_u} g \left(S_k- S_{k-1} \right)  - \mu_g  \right|. $$
By Lemma \ref{lemma-charac-T}, $J_u \uparrow \infty$ a.s.\ as $u \to \infty$.  The second term tends to zero $\Prob^\alpha$-a.s., by Lemma \ref{lem:ergodic theorem}.  For the first term,  
use the Lipschitz continuity of $g$ to infer that for some finite constant $B_g$,
\begin{align*} 
 \frac{1}{J_u} \sum_{k=1}^{J_u} & \bigg| g\left( \log \bigg(  \frac{|V_k|}{|V_{k-1}|} \bigg) \right)  - g \big( S_k - S_{k-1} \big)   \bigg|   ~\le~ \frac{B_g}{J_u} \sum_{k=1}^{J_u} \bigg( \Big| \log |V_k| - \log |V_{k-1}|  - (S_k - S_{k-1})  \Big| \bigg). 
\end{align*}  
Also, it follows directly from the definitions (as given in \eqref{def-S} and \eqref{def-Z}) that
\begin{equation} \label{sec8.3}
\Big| \log |V_k| - \log |V_{k-1}|  - (S_k - S_{k-1})  \Big| :=  \Big| \log |Z_k| - \log |Z_{k-1}| \Big|.
\end{equation}
Now by Lemma \ref{lem:convergence Z}, $|Z_k|$ converges a.s.\ to the proper random variable $|Z|$ (and thus forms  a Cauchy sequence).  
Then by C\'esaro's theorem,
\begin{equation} \label{sec8.4}
\limsup_{u \to \infty} \frac{1}{J_u} \sum_{k=1}^{J_u} \Big| \log |Z_k| - \log |Z_{k-1}| \Big| = 0 \quad \mbox{a.s.},
\end{equation}
and we conclude that $\limsup_{u \to \infty} {\mathbb I}_1(u) =0$ $\Prob^\alpha$-a.s.
This establishes \eqref{empirical distribution function}. 
\halmos

\nocite{JCAV13b, ADSK91a, ADSK91b, Melfi1992, Melfi1994, RP94}

\vspace{.3cm}
\noindent
{\small
{\sc J.\ F.\ Collamore},
Department of Mathematical Sciences,
University of Copenhagen,
Universitetsparken 5,\\
\hspace*{.4cm} DK-2100 Copenhagen $\emptyset$,
Denmark\\
\hspace*{.4cm} E-mail: collamore$@$math.ku.dk

\vspace{.25cm}
\noindent
{\small
{\sc S.\ Mentemeier},
Technische Universit\"{a}t Dortmund,
Fakult\"{a}t f\"{u}r Mathematik, Lehrstuhl IV,\\
\hspace*{.4cm} Vogelpothsweg 87, 44227 Dortmund\\
\hspace*{.4cm} E-mail:  sebastian.mentemeier$@$tu-dortmund.de}
 
\end{document}